\documentclass[a4paper]{article}

\usepackage{geometry}
\usepackage{amsmath}
\usepackage{amssymb}
\usepackage{amsthm}
\usepackage{bm}
\usepackage[backend=biber,bibstyle=alphabetic,citestyle=alphabetic,sorting=nyt,maxnames=6,giveninits=false]{biblatex}
	\AtEveryBibitem{\clearfield{doi}}
	\AtEveryBibitem{\clearfield{isbn}}
	\AtEveryBibitem{\clearfield{issn}}
	\AtEveryBibitem{\clearfield{url}}
	\renewbibmacro{in:}{\ifentrytype{article}{}{\printtext{\bibstring{in}\intitlepunct}}}
	\AtEveryBibitem{\ifentrytype{book}{}{\ifentrytype{incollection}{}{\clearfield{number}}}}
	\AtEveryBibitem{\ifentrytype{book}{\clearfield{pages}}{}}
	\DeclareFieldFormat[article,incollection,inbook,unpublished]{title}{#1}
	\DeclareFieldFormat[article]{volume}{\textrm{#1}}
	\DeclareFieldFormat[article]{pages}{#1}
	\DeclareFieldFormat[book]{edition}{#1}
\usepackage{calc}
\usepackage{csquotes}
\usepackage[british]{babel}
\usepackage[UKenglish]{isodate}
\usepackage{dsfont}
\usepackage{upgreek}
\usepackage{enumitem}
\usepackage[OT2,OT1]{fontenc}
\usepackage{ifthen}
\usepackage[cal=cm,scr=boondoxo]{mathalfa}
\usepackage{pifont}
\usepackage{stmaryrd}
\usepackage{textcomp}
\usepackage{tikz}
	\usetikzlibrary{arrows}
	\usetikzlibrary{patterns}
\usepackage{tikz-cd}
\usepackage[textwidth=3cm,colorinlistoftodos]{todonotes}
\usepackage{mathtools}
\usepackage[capitalise]{cleveref}
	\crefformat{equation}{\textup{#2(#1)#3}}

\numberwithin{equation}{section}			
\swapnumbers						

\newcommand\cyr{
\renewcommand\rmdefault{wncyr}%
\renewcommand\sfdefault{wncyss}%
\renewcommand\encodingdefault{OT2}%
\normalfont
\selectfont}
\DeclareTextFontCommand{\textcyr}{\cyr}

\newcounter{Enum}					
\newenvironment{Enumerate}{\begin{enumerate}[label={\rm({\roman*})}]}{\end{enumerate}}
\newcommand{\Enumref}[1]{{\setcounter{Enum}{#1}{\rm(\roman{Enum})}}}

\newcounter{Enumalph}					
\newenvironment{Enumeratealph}{\begin{enumerate}[label={\rm({\alph*})}]}{\end{enumerate}}

\newcommand{\descriptionlabelsave}{}			
\newenvironment{Itemize}{%
	\renewcommand{\descriptionlabelsave}{\descriptionlabel}\renewcommand{\descriptionlabel}{$\triangleright$}%
	\begin{description}[leftmargin=15pt,itemindent=-5.2pt]}{%
	\end{description}\renewcommand{\descriptionlabel}{\descriptionlabelsave}}

\newcounter{StepsCount}					
\newenvironment{Steps}{%
	\begin{list}{\ding{\value{StepsCount}}}{\usecounter{StepsCount} \leftmargin=0pt \labelwidth=12pt \itemindent=\labelwidth%
	\itemsep=5pt\listparindent=\parindent} \setcounter{StepsCount}{191}}{\end{list}}

\theoremstyle{plain}
	\newtheorem{lemma}{Lemma}[section]
	\newtheorem{proposition}[lemma]{Proposition}
	\newtheorem{theorem}[lemma]{Theorem}
	\newtheorem{corollary}[lemma]{Corollary}
\theoremstyle{definition}
	\newtheorem{definitioN}[lemma]{Definition}
	\newtheorem{assumptioN}[lemma]{Assumption}
	\newtheorem{notatioN}[lemma]{Notation}
\theoremstyle{remark}
	\newtheorem{remarK}[lemma]{Remark}
	\newtheorem{examplE}[lemma]{Example}
\renewcommand{\qedsymbol}{\raisebox{-2pt}{\large\ding{113}}}
\newcommand{\defendsymbol}{$\lozenge$}
\newcommand{\qedsymbolsave}{\qedsymbol}
\newenvironment{definition}{\begin{definitioN}}{
	\renewcommand{\qedsymbolsave}{\qedsymbol}\renewcommand{\qedsymbol}{\defendsymbol}
	\popQED{\qed}\renewcommand{\qedsymbol}{\qedsymbolsave}\end{definitioN}}
\newenvironment{assumption}{\begin{assumptioN}}{
	\renewcommand{\qedsymbolsave}{\qedsymbol}\renewcommand{\qedsymbol}{\defendsymbol}
	\popQED{\qed}\renewcommand{\qedsymbol}{\qedsymbolsave}\end{assumptioN}}
\newenvironment{notation}{\begin{notatioN}}{
	\renewcommand{\qedsymbolsave}{\qedsymbol}\renewcommand{\qedsymbol}{\defendsymbol}
	\popQED{\qed}\renewcommand{\qedsymbol}{\qedsymbolsave}\end{notatioN}}
\newenvironment{remark}{\begin{remarK}}{
	\renewcommand{\qedsymbolsave}{\qedsymbol}\renewcommand{\qedsymbol}{\defendsymbol}
	\popQED{\qed}\renewcommand{\qedsymbol}{\qedsymbolsave}\end{remarK}}
\newenvironment{example}{\begin{examplE}}{
	\renewcommand{\qedsymbolsave}{\qedsymbol}\renewcommand{\qedsymbol}{\defendsymbol}
	\popQED{\qed}\renewcommand{\qedsymbol}{\qedsymbolsave}\end{examplE}}

\crefname{definitioN}{Definition}{Definitions}
\Crefname{definitioN}{Definition}{Definitions}
\crefname{assumptioN}{Assumption}{Assumptions}
\Crefname{assumptioN}{Assumption}{Assumptions}
\crefname{notatioN}{Notation}{Notations}
\Crefname{notatioN}{Notation}{Notations}
\crefname{remarK}{Remark}{Remarks}
\Crefname{remarK}{Remark}{Remarks}
\crefname{examplE}{Example}{Examples}
\Crefname{examplE}{Example}{Examples}

\newcommand{\mc}[1]{{\mathcal{#1}}}			
\newcommand{\ms}[1]{{\mathscr{#1}}}			
\newcommand{\mf}[1]{{\mathfrak{#1}}}			
\newcommand{\bb}[1]{{\mathbb{#1}}}			
\newcommand{\qu}{\overline}				
\newcommand{\mr}{\mathring}				

\DeclareMathOperator{\RE}{Re}				
\renewcommand{\Re}{\RE}
\DeclareMathOperator{\IM}{Im}				
\renewcommand{\Im}{\IM}

\newcommand{\smmatrix}[4]{\Bigl(			
\begin{smallmatrix}
\hspace*{-0.2ex} #1 \hspace*{0.2ex} & \hspace*{0.2ex} #2 \hspace*{-0.2ex}
\\[0.5ex]
\hspace*{-0.2ex} #3 \hspace*{0.2ex} & \hspace*{0.2ex} #4 \hspace*{-0.2ex}
\end{smallmatrix}
\Bigr)}

\DeclareMathOperator{\tr}{tr}				
\DeclareMathOperator{\ran}{ran}				
\DeclareMathOperator{\ind}{ind}				
\newcommand{\Dummy}{\text{\textvisiblespace\kern1pt}}	
\newcommand{\BigO}{{\rm O}}				

\newcommand{\DS}{\colon\mkern3mu}			
\newcommand{\D}{\mathrm{d}}			
\newcommand{\DD}{\mkern4mu\mathrm{d}}			
\newcommand{\RD}{\mkern4mu\mathrm{d}}			
\newcommand{\DP}{{:\kern5pt}}				
\newcommand{\DF}{\colon}				
\newcommand{\DE}{\mathrel{\mathop:}=}			
\newcommand{\DEalign}{\mathrel{\mathop:}\hspace*{-0.72ex}&=}	
\newcommand{\ED}{=\mathrel{\mathop:}}			

\newcommand{\NHam}{{\bb H}^1}			
\newcommand{\Ham}{\bb H}				

\DeclareMathOperator\sgn{sgn}

\DeclareMathOperator\artanh{artanh}
\newcommand\Fhyperg[2]{\prescript{}{#1}{F}_{#2}^{}} 
\DeclareMathOperator\supp{supp}

\newcommand\tn[1]{\widehat{#1}} 

\newcommand{\hatM}{\hspace{0.42ex}\widehat{\rule[1.5ex]{1.5ex}{0ex}}\hspace{-2.2ex}M}
\newcommand{\tildeM}{\hspace{0.42ex}\widetilde{\rule[1.5ex]{1.5ex}{0ex}}\hspace{-2.2ex}M}
\newcommand{\qSL}{q_{\mathsf{SL}}}
\newcommand{\strS}{\mathsf{S}}
\newcommand{\qstr}{q_{\mathsf{S}}}
\newcommand{\qigs}{q_{\mathsf{igs}}}
\newcommand{\sfw}{\mathsf{w}}

\begin{document}

\begin{flushleft}
	{\Large\textbf{Canonical systems whose Weyl coefficients have \\[2mm] regularly varying asymptotics}}
	\\[5mm]
	\textsc{
	Matthias Langer
	\,\ $\ast$\,\
	Raphael Pruckner
	\,\ $\ast$\,\
	Harald Woracek
		\hspace*{-14pt}
		\renewcommand{\thefootnote}{\fnsymbol{footnote}}
		\setcounter{footnote}{2}
		\footnote{The second and third authors were supported by the project P~30715-N35
			of the Austrian Science Fund (FWF). The third author was supported by the
			project I~4600 of the Austrian Science Fund (FWF).}
		\renewcommand{\thefootnote}{\arabic{footnote}}
		\setcounter{footnote}{0}
	}
	\\[6mm]
\end{flushleft}
	{\small
	\textbf{Abstract:}
		For a two-dimensional canonical system $y'(t)=zJH(t)y(t)$ on an interval $(0,L)$ with $0<L\le\infty$
		whose Hamiltonian $H$ is a.e.\ positive semidefinite, denote by $q_H$ its Weyl coefficient.
		De~Branges' inverse spectral theorem states that the assignment
		$H\mapsto q_H$ is a bijection between trace-normalised Hamiltonians and Nevanlinna functions.
		We prove that $q_H$ has an asymptotics towards $i\infty$ whose leading term
		is some (complex) multiple of a regularly varying function
		if and only if the primitive $M$ of $H$ is regularly or rapidly varying at $0$
		and its off-diagonal entries do not oscillate too much.
		The leading term in the asymptotics of $q_H$ towards $i\infty$ is related
		to the behaviour of $M$ towards $0$ by explicit formulae.
		The speed of growth in absolute value depends only on the diagonal entries of $M$,
		while the argument of the leading coefficient corresponds to the relative size
		of the off-diagonal entries.
		Translated to the spectral measure $\mu_H$ and the Hamiltonian $H$,
		this means that the diagonal of $H$ determines the growth of the
		symmetrised distribution function of $\mu_H$, and the relative size and
		sign distribution of its off-diagonal is a measure for the asymmetry of $\mu_H$.
		The results are applied to Sturm--Liouville equations, Krein strings and generalised indefinite strings
		to prove similar characterisations for the asymptotics of the corresponding Weyl coefficients.
	}
	\\
\begin{flushleft}
	\small{
	}
\end{flushleft}


%
%
%
\section{Introduction}

We consider two-dimensional \emph{canonical systems} on an interval $(0,L)$:
\begin{equation}\label{Z1}
	y'(t) = zJH(t)y(t),\qquad t\in(0,L),
\end{equation}
where $L>0$ or $L=\infty$, $J\DE\smmatrix 0{-1}10$, $z\in\bb C$,
and the \emph{Hamiltonian} $H$ of the system satisfies
\begin{Enumeratealph}
\item 
$H\DF(0,L)\to\bb R^{2\times 2}$ is measurable and locally integrable on $[0,L)$;
\item 
$H(t)\geq 0$ and $\tr H(t)>0$, $t\in(0,L)$ a.e.;
\item 
$\int_0^L \tr H(t)\RD t=\infty$.
\end{Enumeratealph}
Many differential and difference equations can be written as canonical systems:
one-dimensional Schr\"odinger equations, Sturm--Liouville equations, Dirac systems,
Krein strings, generalised indefinite strings, and the eigenvalue equation
for Jacobi operators.
A crucial construction to obtain a spectral measure is Weyl's nested disc method,
originally invented by H.~Weyl in the context of Sturm--Liouville equations;
see \cite{weyl:1910}.
Given a Hamiltonian $H$, this method produces a function $q_H$, called its
\emph{Weyl coefficient}, via the fundamental solution of \eqref{Z1} and a
geometric observation based on the limit point condition (c) above; 
see \eqref{Z111} in \cref{Z146} for the definition of $q_H$.
The Weyl coefficient is a \emph{Nevanlinna function}, i.e.\ it is analytic in the
open upper half-plane $\bb C^+$ with $\Im q_H(z)\ge 0$, $z\in\bb C^+$,
or identically equal to $\infty$.
All essential properties of $H$ are encoded in $q_H$: by de~Branges' theory
the set of all Hamiltonians---up to reparameterisation---corresponds bijectively
to the set of all Nevanlinna functions;
see \cite{debranges:1968}; an explicit deduction can be found in \cite{winkler:1995}.
Some general reference for the theory of canonical systems are
\cite{atkinson:1964,gohberg.krein:1967,hassi.snoo.winkler:2000,romanov:1408.6022v1,remling:2018,behrndt.hassi.snoo:2020}.

The system \eqref{Z1} has an operator model, which consists of a Hilbert space $L^2(H)$
and a self-adjoint operator\footnote{%
	In some cases, $A_H$ is a multi-valued operator (or linear relation).
	For simplicity of presentation we systematically neglect these cases throughout
	the present introduction.
	Of course, therefore, some of the following statements have to be understood appropriately.
	In fact, this situation corresponds to `$b_H>0$' in \eqref{Z4}, which is excluded
	in our main theorem; see the discussion in \cref{Z90}\,(ii).}%
$A_H$ in $L^2(H)$.
This construction goes back to B.\,C.~Orcutt and I.\,S.~Kac; see \cite{orcutt:1969,kac:1984,kac:1986a};
a more accessible reference is \cite{hassi.snoo.winkler:2000}.
The operator $A_H$ has simple spectrum, and a \emph{spectral measure} $\mu_H$ is obtained
from the Herglotz integral representation of $q_H$:
\begin{equation}\label{Z4}
	q_H(z) = a_H+b_Hz+\int_{\bb R} \biggl(\frac{1}{t-z}-\frac{t}{1+t^2}\biggr)\RD\mu_H(t),
	\qquad z\in\bb C^+;
\end{equation}
here $a_H\in\bb R$, $b_H\ge0$ and $\mu_H$ is a positive measure on $\bb R$ that
satisfies $\int_{\bb R} \frac{\RD\mu_H(t)}{1+t^2}<\infty$.
The operator $A_H$ is unitarily equivalent to the operator of multiplication 
by the independent variable in the space $L^2(\mu_H)$ via a generalised Fourier transform.

The main objective of this paper is to answer the question how the asymptotic behaviour 
of $\mu_H$ towards $\pm\infty$ relates to $H$.  
The behaviour of $\mu_H$ towards infinity is often called the high-energy asymptotics.
In order to study this question, it turns out to be appropriate to work
with the Weyl coefficient $q_H$ rather than with the measure $\mu_H$ directly
and consider the asymptotic behaviour of $q_H(ri)$ as $r\to\infty$.
In view of \eqref{Z4}, Abelian and Tauberian theorems can be used to translate
between $\mu_H$ and $q_H$.  

The asymptotic behaviour of Weyl coefficients and spectral measures has been studied
for a long time, especially in the context of Sturm--Liouville equations.
Maybe the paper \cite{marchenko:1952} by V.\,A.~Marchenko can be viewed as a
starting point, where the asymptotics of the spectral measure for one-dimensional
Schr\"odinger equations is determined.
A selection of other references dealing mainly with the Sturm--Liouville case are
\cite{everitt:1972,kac:1973a,atkinson:1981,kaper.kwong:1986,bennewitz:1989,clark.gesztesy:2002,rybkin:2002,luger.teschl.woehrer:2016,
sakhnovich:2025}.
Our present work grew out of a line of research which was started
by Y.~Kasahara in \cite{kasahara:1975} for Krein strings,
and continued in \cite{kasahara.watanabe:2010} for so-called Kotani strings,
in \cite{bennewitz.wood:1997} for higher-order equations 
and in \cite{eckhardt.kostenko.teschl:2018} for canonical systems
of the form \eqref{Z1}.
In the latter paper \cite{eckhardt.kostenko.teschl:2018} it is assumed
in most results that the diagonal entries of $H$ dominate the off-diagonal ones,
which leads to situations that are very close to Krein strings and Sturm--Liouville
equations.  Further, in that paper some more delicate cases where the growth of $q_H$
is close to a constant or close to maximal or minimal growth are treated only partially.

The main result of the current paper, \cref{Z10}, gives a full characterisation when the 
Weyl coefficient $q_H$, along the imaginary axis, is asymptotically equal to a constant times 
a regularly varying function in terms of the asymptotic behaviour of the
primitive of $H$ at the left endpoint $0$.
Regularly varying functions are often used as comparison functions in asymptotic analysis,
and they have many properties that are similar to power functions;
for the definition of regularly varying functions see \cref{Z151,Z152}.

In his celebrated inverse spectral theorem L.~de Branges proved that the Weyl coefficient 
determines the Hamiltonian uniquely up to reparameterisation.
In \cref{Z236}  we consider an asymptotic counterpart.
In particular, the asymptotic behaviour of the Weyl coefficient determines the Hamiltonian 
up to asymptotic equivalence of the primitive and reparameterisation, almost every regularly varying asymptotics
can be achieved, and the connection is constructive.

Similar to \cite{kasahara:1975,bennewitz.wood:1997,kasahara.watanabe:2010,eckhardt.kostenko.teschl:2018},
the starting point of the proof of \cref{Z10} is the rescaling trick invented by Y.~Kasahara,
which makes it possible to relate the asymptotics of $q_H$ to the convergence of
certain rescalings of $H$; see \cref{Z128}.  The essential steps in our proof then are to explicitly solve
the inverse problem for Weyl coefficients being an arbitrary complex multiple of a power
(\cref{Z132} and \cref{Z133}), and to relate regular or rapid variation
of the entries of the primitive of $H$ to convergence of rescalings of $H$, 
which is done in \cref{Z63}.

In \cref{Z64} we translate the asymptotic behaviour of the Weyl coefficient $q_H$
to the asymptotic behaviour of the distribution function of the spectral measure $\mu_H$
using Abelian and Tauberian theorems from \cite{langer.woracek:kara}.
In the generic case there is a relation between the asymmetry of the spectral measure $\mu_H$
and the value of a limit that relates off-diagonal to diagonal entries of $H$.
When the growth of $q_H$ is very slow, then a full translation between the
behaviour of $q_H$ and $\mu_H$ is not possible.
Finally, in \cref{Z65} we apply our main theorem to three different scalar equations
to obtain new characterisations for the asymptotic behaviour of the corresponding Weyl coefficients:
Sturm--Liouville equations, Krein strings and generalised indefinite strings; 
the latter were introduced in \cite{eckhardt.kostenko:2016} in connection with
the Camassa--Holm equation.

\subsection{The main theorem}
\label{Z149}

Before we formulate our main theorem (\cref{Z10} below), we recall the definition of regularly varying functions
(since we use this notion to describe the asymptotic behaviours of $q_H$ and the primitive of $H$), and fix some further notation.

\begin{definition}\label{Z151}
Let $\mc I\subseteq(0,\infty)$ be an interval with $\sup\mc I=\infty$.
A measurable function $\ms f\DF\mc I\to(0,\infty)$
is called \emph{regularly varying at $\infty$ with index $\alpha\in\bb R$} if
\begin{equation}\label{Z32}
	\lim_{r\to\infty}\frac{\ms f(\lambda r)}{\ms f(r)} = \lambda^\alpha
	\qquad\text{for all} \ \lambda>0.
\end{equation}
We write $\ind\ms f\DE\alpha$ and denote the set of regularly varying functions with index $\alpha\in\bb R$ 
at $\infty$ by $R_\alpha(\infty)$.

\end{definition}

\medskip

\noindent
Examples include functions $\ms f$ behaving, for large $r$, like
\begin{equation}\label{Z193}
	\ms f(r) = r^\alpha\cdot\bigl(\log r\bigr)^{\beta_1}\cdot\bigl(\log\log r\bigr)^{\beta_2}
	\cdot\ldots\cdot
	\bigl(\underbrace{\log\cdots\log}_{\text{\footnotesize$m$\textsuperscript{th} iterate}}r\bigr)^{\beta_m},
\end{equation}
where $\alpha,\beta_1,\ldots,\beta_m\in\bb R$.
Other examples are $\ms f(r) = r^\alpha e^{(\log r)^\beta}$ with $\beta\in(0,1)$,
$\ms f(r) = r^\alpha e^{\frac{\log r}{\log\log r}}$,
and $\ms f(r)=r^\alpha e^{(\log r)^\beta \cos((\log r)^\beta)}$ with $\beta\in(0,\frac12)$;
note that the last function oscillates.
All these functions have index $\alpha$;
see \cite[Section~1.3]{bingham.goldie.teugels:1989}.
A regularly varying function with index $0$ is also called \emph{slowly varying}.

The property of having regular variation can equally well be considered for $r$ tending to $0$
instead of $\infty$.  At $0$ we also need the notion of rapid variation.

\begin{definition}\label{Z152}
	Let $\mc I\subseteq(0,\infty)$ be an interval with $\inf\mc I=0$.
	Further, let $\ms f\DF\mc I\to(0,\infty)$ be a measurable function.
	\begin{Enumerate}
	\item
	$\ms f$ is called \emph{regularly varying at $0$ with index $\alpha\in\bb R$} if
	\begin{equation}\label{Z33}
		\lim_{t\to 0}\frac{\ms f(\lambda t)}{\ms f(t)} = \lambda^\alpha
		\qquad\text{for all} \ \lambda>0.
	\end{equation}
	\item
	$\ms f$ is called \emph{rapidly varying at $0$ with index $\infty$}
	if \eqref{Z33} holds for $\alpha=\infty$, where we set
	\[
		\lambda^{\infty}\DE
		\begin{cases}
			0, & \lambda\in(0,1),
			\\[0.5ex]
			1, & \lambda=1,
			\\[0.5ex]
			\infty, & \lambda\in(1,\infty).
		\end{cases}
	\]
	\end{Enumerate}
	In both cases we write $\ind\ms f\DE\alpha$.

	We denote the set of regularly or rapidly varying functions with index $\alpha\in(-\infty,\infty]$ 
	at $0$ by $R_\alpha(0)$.
\end{definition}

\noindent
Note that $\ms f$ is regularly varying at $0$ with index $\alpha$ if and only if
the function $\ms g(r)\DE \ms f(r^{-1})^{-1}$ has the respective property at $\infty$.
An example for a rapidly varying function at 0 with index $\infty$
is $\ms f(t)\DE e^{-\frac 1t}$.  

\begin{notation}\label{Z205}
\rule{1ex}{0ex}
\begin{Itemize}
\item 
	We always use the branches of the logarithm and complex powers
	that are analytic on $\bb C\setminus(-\infty,0]$ and take the value $0$ or $1$,
	respectively, at the point $1$.
\item 
	We set $\bb C^+\DE\{z\in\bb C\DS \Im z>0\}$ and $\bb C^-\DE\{z\in\bb C\DS \Im z<0\}$.
\item
	We use the notation $f\sim g$ to express that $\frac fg\to 1$.
	Often we use this notation also locally uniformly with respect to a parameter, e.g.\
	$f(r,z)\sim g(r,z)$ as $r\to\infty$ locally uniformly in $z$ means $\lim_{r\to\infty}\frac{f(r,z)}{g(r,z)}=1$
	locally uniformly in $z$.
	Further, we write \\
	\hspace*{1ex} $f\ll g$ if $\frac fg\to 0$; \\[0.5ex]
	\hspace*{1ex} $f\lesssim g$ if there exists a constant $c>0$ such that $f\le cg$; \\[0.5ex]
	\hspace*{1ex} $f\asymp g$ if $f\lesssim g$ and $g\lesssim f$. \\[0.5ex]
	The domain of validity of these inequalities will be stated or will be clear from the context.
\item
	We denote by $\mathds{1}_M$ the indicator function of a set $M$.
\vspace*{-4ex}
\end{Itemize}
\end{notation}

\begin{notation}\label{Z206}
	Let $H=\smmatrix{h_1}{h_3}{h_3}{h_2}$ be a Hamiltonian on an interval $(0,L)$
	that satisfies the assumptions \textup{(a)--(c)} stated at the beginning of the Introduction.
	For $t\in[0,L)$ we set
	\begin{align*}
		& m_j(t) \DE \int_0^t h_j(s)\RD s,\;\; j=1,2,3,\qquad
		M(t) \DE \begin{pmatrix} m_1(t) & m_3(t) \\ m_3(t) & m_2(t) \end{pmatrix},
		\\[1ex]
		& \mf t(t) \DE \int_0^t\tr H(s)\RD s = \tr M(t) = m_1(t)+m_2(t).
	\end{align*}
\end{notation}

\medskip

\noindent
Now we are ready to state the main result of the paper.
It gives a characterisation when the Weyl coefficient of a Hamiltonian $H$
is asymptotically equal to a constant times a regularly varying function
along the imaginary axis.
The characterisation is given in terms of the behaviour
of the primitives of the entries of the Hamiltonian towards the left endpoint.

\begin{theorem}\label{Z10}
	Let $H=\smmatrix{h_1}{h_3}{h_3}{h_2}$ be a Hamiltonian on an interval $(0,L)$
	that satisfies the assumptions \textup{(a)--(c)} stated at the beginning of the Introduction,
	let $m_j$ and $\mf t$ be as in \cref{Z206},
	and let $q_H$ be the Weyl coefficient of $H$.
	Further, assume that neither $h_1$ nor $h_2$ vanishes a.e.\ on any neighbourhood of $0$
	and that $\mf t$ is regularly varying at $0$ with positive index.

	Then the following two statements are equivalent.
	\begin{Enumerate}
	\item
		There exist a function $\ms a\DF(0,\infty)\to(0,\infty)$ that is
		regularly varying at $\infty$ and a constant $\omega\in\bb C\setminus\{0\}$, such that
		\begin{equation}\label{Z110}
			q_H(ri) \sim i\omega\ms a(r)
			\qquad \text{as} \ r\to\infty.
		\end{equation}
	\item
		The functions $m_1$ and $m_2$ are regularly or rapidly varying at $0$, 
		and, in the case when both $m_1$ and $m_2$ are regularly varying, the limit
		\begin{equation}\label{Z56}
			\delta\DE\lim\limits_{t\to0}\frac{m_3(t)}{\sqrt{m_1(t)m_2(t)}\,}
		\end{equation}
		exists.
	\end{Enumerate}
	Assume that \textup{(i)} and \textup{(ii)} hold. Then the following statements hold.
	\begin{Itemize}
	\item At most one of $m_1$ and $m_2$ is rapidly varying. 
		If $m_1$ or $m_2$ is rapidly varying, then the limit in \eqref{Z56} exists and $\delta=0$.
	\item Let $\mr t$ be the strictly decreasing bijection from $(0,\infty)$ onto $(0,L)$
		such that $\mr t(r)$ is the unique number that satisfies
		\begin{equation}\label{Z38}
			(m_1m_2)\bigl(\mr t(r)\bigr) = \frac{1}{r^2}
		\end{equation}
		for $r>0$, set $\rho_i \DE \ind m_i$, $i\in\{1,2\}$,
		and define the function
		\begin{equation}\label{Z109}
			\ms a_H(r)\DE\sqrt{\frac{m_1(\mr t(r))}{m_2(\mr t(r))}\,},\qquad r>0.
		\end{equation}
		Then $\ms a_H$ is regularly varying at $\infty$ with index 
		\begin{equation}\label{Z78}
			\alpha = \frac{\rho_2-\rho_1}{\rho_2+\rho_1}
		\end{equation}
		which is interpreted as $+1$ when $\rho_2=\infty$ and as $-1$ when $\rho_1=\infty$.
		Moreover, $\frac{1}{r}\ll\ms a_H(r)\ll r$ as $r\to\infty$.
	\item The Weyl coefficient has the asymptotics
		\begin{equation}\label{Z57}
			q_H(rz) \sim i\omega_{\alpha,\delta}\Bigl(\frac zi\Bigr)^\alpha\ms a_H(r) \qquad \text{as} \ r\to\infty
		\end{equation}
		locally uniformly for $z\in\bb C^+$, where \textup{(}explicit formulae are given in \eqref{Z98}\textup{)}
		\[
			\omega_{\alpha,\delta}\in\bb C\setminus\{0\},\quad 
			|\arg\omega_{\alpha,\delta}| \le \frac{\pi}{2}\bigl(1-|\alpha|\bigr).
		\]
	\end{Itemize}
\end{theorem}

\medskip

\noindent
The proof of \cref{Z10} is given in \cref{Z61,Z63}.
Let us add here a few remarks concerning this theorem.

\begin{remark}\label{Z90}
\rule{1ex}{0ex}
\begin{Enumerate}
\item
	Additional relations among $\alpha,\omega_{\alpha,\delta},\delta,\rho_1,\rho_2$ are given in \cref{Z100} below.
\item
	Since $q_H$ is a Nevanlinna function, we have
	\[
		\frac{1}{r} \lesssim |q_H(ri)| \lesssim r, \qquad r>0.
	\]
	The following equivalences show that the cases of vanishing $h_1$ or $h_2$
	in a neighbourhood of $0$, which are excluded in \cref{Z10}, correspond to the 
	extreme cases in the growth behaviour of $q_H$, and the asymptotic behaviour is
	known exactly in these cases (see, e.g.\ \cite{kac.krein:1968a,winkler:1995}).
	For maximal growth we have
	\begin{alignat*}{3}
		&\exists\varepsilon>0\DP h_2|_{(0,\varepsilon)}=0\text{ a.e.}
		\quad&&\Longleftrightarrow\quad
		\exists b_1>0\DP q_H(ri) \sim ib_1r
		\quad&&\Longleftrightarrow\quad
		\limsup_{r\to\infty}\frac 1r|q_H(ri)|>0
		\\
		&&&\Longleftrightarrow\quad
		b_H>0 \quad \text{(with $b_H$ as in \eqref{Z4})};
		\hspace*{-40ex}
	\end{alignat*}
	if one of these equivalent statements holds, then $b_1=\int_0^{c_1} h_1(s)\RD s$ where
	$c_1\DE\max\{\varepsilon>0\DS h_2|_{(0,\varepsilon)}=0\}$.

	For minimal growth we have:
	\begin{alignat*}{3}
		&\exists\varepsilon>0\DP h_1|_{(0,\varepsilon)}=0\text{ a.e.}
		\quad&&\Longleftrightarrow\quad
		\exists b_2>0\DP q_H(ri) \sim \frac{i}{b_2}\cdot\frac{1}{r}
		\quad&&\Longleftrightarrow\quad
		\liminf_{r\to\infty}r|q_H(ri)|<\infty
		\\
		&&&\Longleftrightarrow\quad
		q_H(z) = \int_{\bb R} \frac{1}{t-z}\RD\mu_H(t)
		\quad \text{with a finite measure} \ \mu_H;
		\hspace*{-40ex}
	\end{alignat*}
	if one of these equivalent statements holds, then $b_2=\int_0^{c_2} h_2(s)\RD s$ where
	$c_2\DE\max\{\varepsilon>0\DS h_1|_{(0,\varepsilon)}=0\}$,
	and $\mu_H(\bb R)=\frac{1}{b_2}$.
\item
	The assumption that $\mf t$ is regularly varying at $0$ with positive index
	is no loss of generality since we could even assume $H$ to be
	trace-normalised, i.e.\ $\tr H(t)=1$ a.e.\ and hence $\mf t(t)=t$; this can be achieved by
	reparameterisation (i.e.\ a change of the independent variable); see \cref{Z146}.
\item
	The expression $\ms a_H$ has been used recently in connection with uniform
	bounds for the Weyl coefficient: it follows from \cite[Theorem~1.1]{langer.pruckner.woracek:heniest},
	together with a simple rescaling, that $|q_H(ri)|\asymp\ms a_H(r)$, $r>0$,
	with constants in $\asymp$ independent of $H$.
	This implies, in particular, that $q_H(ri)\to0$ as $r\to\infty$ if and only if $\frac{m_1(t)}{m_2(t)}\to0$
	as $t\to0$; see \cref{Z174}\,(i).
\item
	Instead of an exact solution of \eqref{Z38} it is sufficient fo find a function that is asymptotically
	equal to $\mr t(r)$.  For $\alpha\in(-1,1)$ this follows from \cite[Theorem~1.8.7]{bingham.goldie.teugels:1989}.
	When $\alpha\in\{1,-1\}$, i.e.\ $m_1$ or $m_2$ is rapidly varying, one has to use \eqref{Z17} below 
	instead of \eqref{Z109}.
\item
	We chose the form $i\omega_{\alpha,\delta}(\frac{z}{i})^\alpha$ as the coefficient of $\ms a_H(r)$ on the
	right-hand side of \eqref{Z57} because $\arg\omega_{\alpha,\delta}$ reflects 
	asymptotic asymmetry of $q_H$.  Note that $q_H$ is a symmetric Nevanlinna function 
	(i.e.\ $q_H(-\overline z)=-\overline{q_H(z)}$ for $z\in\bb C^+$)
	if and only if $q_H(ri)\in i\bb R$ for $r>0$.
	A relation between $\arg\omega_{\alpha,\delta}$ and $\delta$ is given in \eqref{Z6}.
\item
	When $\alpha\ne0$, then the validity of an asymptotic behaviour
	as in \eqref{Z57} can be translated to an asymptotics of the distribution function
	of the spectral measure $\mu_H$ towards $\pm\infty$.
	The diagonal of $H$ determines the growth of the symmetrised distribution function,
	and $\delta$ determines the asymmetry of $\mu_H$.
	In the exceptional case $\alpha=0$ a full translation is not possible as counterexamples show,
	but partial results covering many cases exist.
	A detailed discussion is contained in \cref{Z64}.
\item
	Let us relate \cref{Z10} to \cite{eckhardt.kostenko.teschl:2018},
	which contains partial results:
	(a) the case when $q_H(ri)\to\zeta_0$ for some $\zeta_0\in\bb C\setminus\{0\}$;
	(b) the situation when $\alpha\in(-1,1)\setminus\{0\}$ and $\delta=0$;
	(c) one implication when $\alpha\in\{1,-1\}$.
	In particular, the restriction $\delta=0$ in (b) means that those results
	are very close to the situation of a Krein string considered in \cite{kasahara:1975}.
	Also the exclusion of $\alpha=0$ apart from the situation in (a) is
	quite restrictive as it excludes any slowly varying function for $\ms a_H$
	apart from those converging to a positive constant.
\end{Enumerate}
\end{remark}

\subsection{Explicit formulae relating constants}
\label{Z100}

Assume that $H$ is given as in \cref{Z10} and that the equivalent properties (i), (ii) in \cref{Z10} hold.
Further, denote by $\sigma$ the index of $\mf t$, and let $\delta,\alpha,\omega_{\alpha,\delta},\rho_1,\rho_2$ 
be the quantities appearing in \cref{Z10}.
We list some facts that will be seen in the course of the proof of \cref{Z10}.
\begin{Steps}
\item
	The following formulae hold:
	\begin{align}
		& \min\{\rho_1,\rho_2\}=\sigma,
		\label{Z86}
		\\[1ex]
		& \frac{\rho_2}{\rho_1} = \frac{1+\alpha}{1-\alpha},
		\label{Z194}
	\end{align}
	where both sides of the latter relation are understood as $\infty$ when $\alpha=1$, $\rho_1\in(0,\infty)$, $\rho_2=\infty$.
\item
	We have the bound
	\begin{equation}\label{Z97}
		|\delta| \le \sqrt{1-\alpha^2}.
	\end{equation}
	Note that $\sqrt{1-\alpha^2}=\sqrt{\rho_1\rho_2}\cdot\frac{2}{\rho_1+\rho_2}$
	if $\max\{\rho_1,\rho_2\}<\infty$.
\item 
	The constant $\omega_{\alpha,\delta}$ and its argument are given by
	\begin{align}
		\omega_{\alpha,\delta}
		&=
		\begin{cases}
			\displaystyle
			\bigl(2\sqrt{1-\alpha^2-\delta^2}\bigr)^{1+\alpha}
			\frac{\Gamma(-\alpha)}{\Gamma(2+\alpha)}
			\cdot\frac{\Gamma\bigl(1+\frac{\alpha}{2}\bigl(1+i\frac{\delta}{\sqrt{1-\alpha^2-\delta^2}\,}\bigr)\bigr)}{%
			\Gamma\bigl(-\frac{\alpha}{2}\bigl(1-i\frac{\delta}{\sqrt{1-\alpha^2-\delta^2}\,}\bigr)\bigr)}
			\hspace*{-20ex}\\[3ex]
			& \text{if} \ \alpha\ne0,\pm1,\, |\delta|<\sqrt{1-\alpha^2},
			\\[1ex]
			\displaystyle
			(i\alpha\delta)^{1+\alpha}\cdot\frac{\Gamma(-\alpha)}{\Gamma(2+\alpha)}
			& \text{if} \ \alpha\ne0,\pm1,\, |\delta|=\sqrt{1-\alpha^2},
			\\[2ex]
			\displaystyle
			\sqrt{1-\delta^2}-i\delta
			& \text{if} \ \alpha=0,
			\\[1ex]
			1
			& \text{if} \ \alpha=1 \ \text{or} \ \alpha=-1,
		\end{cases}
	\label{Z98}
	\\[1ex]
		\arg\omega_{\alpha,\delta}
		&=
		\begin{cases}
			\displaystyle
			-\arctan\biggl[\tan\Bigl(\frac{\pi}{2}\bigl(1-|\alpha|\bigr)\Bigr)
			\tanh\Bigl(\frac{\pi|\alpha|\delta}{2\sqrt{1-\alpha^2-\delta^2}\,}\Bigr)\biggr]
			\hspace*{-15ex}\\[3ex]
			& \text{if} \ \alpha\ne0,\, |\delta|<\sqrt{1-\alpha^2},
			\\[2ex]
			\displaystyle
			-(\sgn\delta)\frac{\pi}{2}\bigl(1-|\alpha|\bigr)
			& \text{if} \ \alpha\ne0,\, |\delta|=\sqrt{1-\alpha^2},
			\\[3ex]
			\displaystyle
			-\arcsin\delta
			& \text{if} \ \alpha=0.
		\end{cases}
	\label{Z6}
	\end{align}
	When $\alpha\in[-1,1]$ is fixed, then the function $\delta\mapsto\arg\omega_{\alpha,\delta}$
	is a strictly decreasing and odd bijection from $[-\sqrt{1-\alpha^2},\sqrt{1-\alpha^2}]$ 
	onto the interval $[-\frac\pi2(1-|\alpha|),\frac\pi2(1-|\alpha|)]$.
	In particular, $\arg\omega_{\alpha,\delta}=0$ if and only if $\delta=0$.
	Moreover, $\omega_{\alpha,\delta}\in i\bb R$ if and only if $|\delta|=1$, which is only 
	possible when $\alpha=0$.
\item
	The modulus of $\omega_{\alpha,\delta}$ and the number $\delta$ in \eqref{Z56} can be expressed 
	in terms of $\alpha$ and $\phi\DE\arg\omega_{\alpha,\delta}$:
	with
	\[
		c \DE \dfrac{2}{\pi|\alpha|}
		\artanh\biggl[\dfrac{\tan\phi}{\tan(\frac{\pi}{2}(1-|\alpha|))}\biggr]
		\qquad\text{if} \ \alpha\ne0, |\phi|<\frac{\pi}{2}(1-|\alpha|)
	\]
	we have
	\begin{align}
		|\omega_{\alpha,\delta}|
		&=
		\begin{cases}
			\displaystyle 2^{1+\alpha}\Bigl(\frac{1-\alpha^2}{1+c^2}\Bigr)^{\frac{1+\alpha}{2}}
			\cdot\frac{|\Gamma(-\alpha)|}{\Gamma(2+\alpha)}\cdot
			\frac{\big|\Gamma\bigl(1+\frac{\alpha}{2}(1-ic)\bigr)\big|}{\big|\Gamma\bigl(-\frac{\alpha}{2}(1+ic)\bigr)\big|}
			& \text{if} \ \alpha\ne0,\; |\phi|<\frac{\pi}{2}(1-|\alpha|),
			\\[3ex]
			\displaystyle |\alpha|^\alpha(1-\alpha^2)^{\frac{1+\alpha}{2}}\frac{\Gamma(1-\alpha)}{\Gamma(2+\alpha)}
			& \text{if} \ \alpha\ne0,\; |\phi|=\frac{\pi}{2}(1-|\alpha|),
			\\[2ex]
			1
			& \text{if} \ \alpha=0,
		\end{cases}
		\label{Z79}
		\\[1ex]
		\delta
		&=
		\begin{cases}
			-\dfrac{c\sqrt{1-\alpha^2}\,}{\sqrt{1+c^2}\,}
			& \text{if} \ \alpha\ne0,\; |\phi|<\frac{\pi}{2}(1-|\alpha|),
			\\[2ex]
			-\sgn\phi\cdot\sqrt{1-\alpha^2} \quad
			& \text{if} \ \alpha\ne0,\; |\phi|=\frac{\pi}{2}(1-|\alpha|),
			\\[2ex]
			-\sin\phi
			& \text{if} \ \alpha=0.
		\end{cases}
		\label{Z11}
	\end{align}
\item
	The function $\ms a_H$ defined in \eqref{Z109} can be rewritten as follows:
	\begin{equation}\label{Z17}
		\ms a_H(r) = rm_1\bigl(\mr t(r)\bigr) = \frac{1}{rm_2\bigl(\mr t(r)\bigr)}\,,
		\qquad r>0.
	\end{equation}
\end{Steps}

\subsection{An inverse theorem}
\label{Z150}

In the following theorem we prove an asymptotic inverse result for canonical systems
in the situation when $\alpha\in(-1,1)$ where $\alpha$ is as in \cref{Z10}.
Moreover, we show that generically every asymptotic behaviour occurs (see also \cref{Z266}), and the relations are constructive.

Before we formulate the theorem let us recall the notion of reparameterisation of Hamiltonians.

\begin{definition}\label{Z154}
	Let $H_1$ and $H_2$ be Hamiltonians on the intervals $[0,L_1)$ and $[0,L_2)$, respectively.
	We say that $H_1$ and $H_2$ are \emph{reparameterisations} of each other
	if there exists a strictly increasing bijection $\gamma:(0,L_1)\to (0,L_2)$
	such that both $\gamma$ and $\gamma^{-1}$ are locally absolutely continuous and
	\begin{equation}\label{Z34}
		H_1(t) = \gamma'(t)H_2(\gamma(t)),  \qquad t\in (0,L_1)\ a.e.
	\end{equation}
\end{definition}

\begin{theorem}\label{Z236}
	Let $H$ be a Hamiltonian on an interval $(0,L)$
	that satisfies the assumptions \textup{(a)--(c)} stated at the beginning of the Introduction,
	let $m_j$ be as in \cref{Z206}, and let $q_H$ be the corresponding Weyl coefficient.
	Further, let $\ms f$ be a function that is regularly varying at $\infty$ with index $\alpha\in(-1,1)$,
	and let $\phi\in\bb R$ with $|\phi|\le\frac{\pi}{2}\bigl(1-|\alpha|\bigr)$.
	Let $\delta$ be given by \eqref{Z11}, let $C_{\alpha,\phi}$ be given by the right-hand side of \eqref{Z79},
	and let $\ms g$ be a strictly decreasing, regularly varying function such that
	\begin{equation}\label{Z246}
		\ms g(r) \sim \frac{1}{C_{\alpha,\phi}}\cdot\frac{\ms f(r)}{r}
		\qquad \text{as} \ r\to\infty
	\end{equation}
	\textup{(}such a function $\ms g$ exists, cf.\ \cref{Z253}\textup{)}.
	Then the following two statements are equivalent:
	\begin{Enumerate}
	\item
		$q_H(ri)\sim ie^{i\phi}\ms f(r)$ \quad as $r\to\infty$;
	\item
		there exists a reparameterisation $\widetilde H$ of $H$ 
		with primitive $\tildeM=\smmatrix{\widetilde m_1}{\widetilde m_3}{\widetilde m_3}{\widetilde m_2}$
		such that
		\begin{equation}\label{Z238}
			\widetilde m_1(t) \sim t,
			\qquad
			\widetilde m_2(t) \sim \frac{1}{t\,[\ms g^{-1}(t)]^2},
			\qquad
			\widetilde m_3(t) 
			\begin{cases}
				\sim \dfrac{\delta}{\ms g^{-1}(t)} & \text{if} \ \delta\ne0,
				\\[3ex]
				\ll \dfrac{1}{\ms g^{-1}(t)} & \text{if} \ \delta=0,
			\end{cases}
		\end{equation}
		as $t\to0$.
	\end{Enumerate}
	For every regularly varying $\ms f$ with $\alpha=\ind\ms f\in(-1,1)$ and
	$\phi\in\bb R$ with $|\phi|\le\frac{\pi}{2}\bigl(1-|\alpha|\bigr)$
	there exists a Hamiltonian $H$ such that \textup{(i)} holds.
\end{theorem}

\noindent
The proof of \cref{Z236} is given in \cref{Z243}.

\subsection{Corollaries and examples}
\label{Z203}

In this subsection we formulate some consequences of \cref{Z10,Z236} and consider
some examples.

To start with, let us note that the implication (i)$\Rightarrow$(ii) in \cref{Z236} shows that under certain assumptions 
the asymptotic behaviour of the Weyl coefficient determines the Hamiltonian uniquely up to asymptotic equivalence 
and reparameterisation.

\begin{corollary}
\label{Z5}
	Let $H_1$ and $H_2$ be Hamiltonians with $q_{H_1}(ir)\sim q_{H_2}(ir)\sim ie^{i\phi}\ms f(r)$ where
	$\ms f$ is regularly varying at $\infty$ with index in $(-1,1)$ and $\phi\in\bb R$.
	Then there are reparameterisations $\widetilde H_1$ and 
	$\widetilde H_2$ of $H_1$ and $H_2$, respectively, such that
	\begin{align*}
		& \widetilde m_{H_1,1}(t)\sim\widetilde m_{H_2,1}(t), \qquad 
		\widetilde m_{H_1,2}(t)\sim\widetilde m_{H_2,2}(t) \qquad \text{as} \ t\to0, \;\; \text{and}
		\\[1ex]
		& \lim_{t\to0}\frac{\widetilde m_{H_1,3}(t)}{\sqrt{\widetilde m_{H_1,1}(t)\widetilde m_{H_1,2}(t)}\,}
		= \lim_{t\to0}\frac{\widetilde m_{H_2,3}(t)}{\sqrt{\widetilde m_{H_2,1}(t)\widetilde m_{H_2,2}(t)}\,},
	\end{align*}
	where $\widetilde m_{H_i,j}$ are the entries of the primitive of $\widetilde H_i$.
\end{corollary}

\noindent
Next, we obtain, again as a consequence of \cref{Z236}, a result on existence of Nevanlinna functions with prescribed asymptotics.

\begin{corollary}\label{Z180}
	Let $\ms f$ be a regularly varying function with index $\alpha\in(-1,1)$
	and let $\phi\in\bb R$ with $|\phi|\le\frac{\pi}{2}\bigl(1-|\alpha|\bigr)$.
	Then there exists a Nevanlinna function $q$ such that
	\begin{equation}\label{Z182}
		q(ri) \sim ie^{i\phi}\ms f(r) \qquad \text{as} \ r\to\infty.
	\end{equation}
\end{corollary}

\noindent
For $\alpha\ne0$ this can also be shown by constructing a measure in 
the integral representation \eqref{Z4} and using Abelian theorems; see, e.g.\ \cite[Section~5]{langer.woracek:kara}.
However, for $\alpha=0$ this is not at all clear as there exist only partial Abelian results.

\begin{remark}\label{Z266}
	The cases $\alpha=\pm1$ are more delicate.
	When $\alpha=-1$ and $r\mapsto r\ms f(r)$ is non-decreasing, we can define
	a measure $\mu$ on $\bb R$ such that $\mu((-r,r))\sim r\ms f(r)$,
	set $q(z)\DE\int_{\bb R}\frac{1}{t-z}\DD\mu(t)$
	and apply \cite[Theorem~5.1]{langer.woracek:kara} to deduce that \eqref{Z182} holds
	with $\phi=0$.
	However, there exist regularly varying functions $\ms f$ with index $\alpha=-1$ 
	that satisfy $\ms f(r)\gg\frac{1}{r}$ such that there is no Nevanlinna function $q$
	so that \eqref{Z182} holds with $\phi=0$.
	An example of such a function is
	\[
		\ms f(r) = \frac{1}{r}\Bigl[\log r+e^{(\log r)^{1/3}\cos((\log r)^{1/3})}\Bigr].
	\]
	If there were a Nevanlinna function $q$ with the desired property, 
	then, by \cite[Theorem~5.1]{langer.woracek:kara} the distribution function of the measure $\mu$
	would satisfy $\mu((-r,r))\sim r\ms f(r)$, which is not possible, as can be seen from the fact
	that, for $r_n=e^{(n\pi)^3}$, one has $\frac{r_{2n+1}\ms f(r_{2n+1})}{r_{2n}\ms f(r_{2n})}\to0$
	as $n\to\infty$.
	The case $\alpha=1$ is similar and can be reduced to the case $\alpha=-1$ by considering
	the functions $\frac{1}{\ms f(r)}$ and $-\frac{1}{q(z)}$.
	For existence of a $q$ in this case one has to assume that $r\mapsto\frac{\ms f(r)}{r}$ is non-increasing.
\end{remark}

\medskip

\noindent
In the following corollary we consider a special situation where the assumption in \cref{Z10} 
that $\mf t$ is regularly varying can be dropped.  
This is used in an application to Sturm--Liouville equations in \cref{Z66}.

\begin{corollary}\label{Z215}
	Assume that $h_1(t)\ne0$ a.e., that $h_2$ does not vanish in any neighbourhood of\, $0$ 
	and that $m_1(t)\gg m_2(t)$ as $t\to0$.  Let \textup{(i)} be as in \cref{Z10}
	and consider the statement
	\begin{Enumerate}
	\item[\textup{(ii)$'$}]
	the function $m_2\circ m_1^{-1}$ is regularly or rapidly varying at $0$, and, 
	if $m_2\circ m_1^{-1}$ is regularly varying, then the limit in \eqref{Z56} exists.
	\end{Enumerate}
	Then \textup{(i)} and \textup{(ii)$'$} are equivalent.
	
	Assume now that \textup{(i)} and \textup{(ii)$'$} hold, set $\rho\DE\ind(m_2\circ m_1^{-1})$ 
	and let $\hat t(r)$ be the solution of
	\begin{equation}\label{Z226}
		\hat t(r)(m_2\circ m_1^{-1})\bigl(\hat t(r)\bigr) = \frac{1}{r^2}
	\end{equation}
	for $r>0$.
	Then
	\begin{equation}\label{Z227}
		q_H(rz) \sim i\omega_{\alpha,\delta}\Bigl(\frac{z}{i}\Bigr)^\alpha r\hat t(r) \qquad \text{as} \ r\to\infty,
	\end{equation}
	locally uniformly for $z\in\bb C^+$, where $\alpha=\frac{\rho-1}{\rho+1}$ and $\omega_{\alpha,\delta}$
	is as in \eqref{Z98}.
	
	Moreover, the asymptotic behaviour of $m_2\circ m_2^{-1}$ can be recovered from the behaviour of $q_H$:
	let $\phi\in\bigl[-\frac{\pi}{2},\frac{\pi}{2}\bigr]$ be such that
	\[
		q_H(ri) \sim ie^{i\phi}|q_H(ri)| \qquad \text{as} \ r\to\infty,
	\]
	let $\alpha$ be the index of $r\mapsto|q_H(ri)|$, let $C_{\alpha,\phi}$
	be given the right-hand side of \eqref{Z11},
	and let $\ms g$ be a strictly decreasing and regularly varying function such that
	\begin{equation}\label{Z247}
		\ms g(r) \sim \frac{1}{C_{\alpha,\phi}}\cdot\frac{|q_H(ri)|}{r} \qquad \text{as} \ r\to\infty,
	\end{equation}
	\textup{(}which exists by \cref{Z253}\textup{)};
	then
	\begin{equation}\label{Z248}
		(m_2\circ m_1^{-1})(t) \sim \frac{1}{t[\ms g^{-1}(t)]^2} \qquad \text{as} \ t\to0.
	\end{equation}
\end{corollary}

\noindent
The proof of \cref{Z215} is given in \cref{Z171}.

\bigskip

\begin{example}\label{Z202}
	Let us consider the situation where not only the $m_i$ but even the $h_i$ are regularly varying.
	Assume that $h_1$ and $h_2$ are regularly varying with indices $\rho_1-1$ and $\rho_2-1$ respectively
	with $\rho_1,\rho_2>0$,
	and let $h_3(t)=\kappa\sqrt{h_1(t)h_2(t)}$ with $|\kappa|\le1$.
	If $\kappa\ne0$, then $h_3$ is regularly varying with index $\frac{\rho_1+\rho_2}{2}-1$, and we obtain
	from \cref{Z122}\,(i)
	that $m_1$, $m_2$, $m_3$ are regularly varying and
	\begin{equation}\label{Z184}
		\frac{m_3(t)}{\sqrt{m_1(t)m_2(t)}\,}
		\sim \frac{\frac{2}{\rho_1+\rho_2}th_3(t)}{\sqrt{\frac{1}{\rho_1}th_1(t)\frac{1}{\rho_2}th_2(t)}\,}
		= \kappa\sqrt{\rho_1\rho_2}\cdot\frac{2}{\rho_1+\rho_2}
		= \kappa\sqrt{1-\alpha^2}
	\end{equation}
	with $\alpha$ as in \eqref{Z78}, which shows that the limit in \eqref{Z56} exists.
	Clearly, when $\kappa=0$, then $\delta=0$.
	Hence (ii) in \cref{Z10} is satisfied.
\end{example}

\medskip

\noindent
Let us explicitly formulate a corollary about Hamiltonians
with power asymptotics at $0$; the deduction from \cref{Z10} is carried out in \cref{Z92}.

\begin{corollary}\label{Z124}
	Let $H$ be a Hamiltonian on an interval $(0,L)$
	that satisfies the assumptions \textup{(a)--(c)} stated at the beginning of the Introduction,
	let $m_j$ be as in \cref{Z206}, and let $q_H$ be the corresponding Weyl coefficient.
	Moreover, assume that $\mf t(t)\sim ct^\sigma$ as $t\to0$ for some $c>0$ and $\sigma>0$.
	Then the following statements are equivalent.
	\begin{Enumerate}
	\item
		There exist $\alpha\in(-1,1)$ and $\omega\in\bb C\setminus\{0\}$
		such that
		\[
			q_H(ri) \sim i\omega r^\alpha \qquad \text{as} \ r\to\infty.
		\]
	\item
		There exist $\rho_1,\rho_2>0$, $c_1,c_2>0$ and $c_3\in\bb R$ such that,
		with $\rho_3=\frac{\rho_1+\rho_2}{2}$,
		\begin{equation}\label{Z135}
			m_1(t) \sim c_1t^{\rho_1}, \qquad
			m_2(t) \sim c_2t^{\rho_2}, \qquad
			m_3(t) 
			\begin{cases}
				\sim c_3t^{\rho_3} & \text{if} \ c_3\ne0,
				\\[1ex]
				\ll t^{\rho_3} & \text{if} \ c_3=0,
			\end{cases}
		\end{equation}
		as $t\to0$.
	\end{Enumerate}
	If \textup{(i)} and \textup{(ii)} are satisfied, then
	\begin{equation}\label{Z136}
		\alpha = \frac{\rho_2-\rho_1}{\rho_2+\rho_1}, \qquad
		\omega = c_1^{\frac{\alpha+1}{2}}c_2^{\frac{\alpha-1}{2}}\omega_{\alpha,\delta},
	\end{equation}
	where $\omega_{\alpha,\delta}$ is as in \eqref{Z98} with $\delta=\frac{c_3}{\sqrt{c_1c_2}\,}$,
	and
	\begin{equation}\label{Z137}
		q_H(z) \sim i\omega\Bigl(\frac zi\Bigr)^\alpha
		\qquad\text{as} \ |z|\to\infty
	\end{equation}
	uniformly in each Stolz angle $\{z\in\bb C\DS \psi\le\arg z\le\pi-\psi\}$
	where $\psi\in(0,\frac\pi2)$.
\end{corollary}

\begin{remark}\label{Z259}
\rule{0ex}{1ex}
\begin{Enumerate}
\item
	Note that, according to item \ding{194} in \cref{Z100}, $\sgn(\arg\omega)=-\sgn c_3$;
	in particular, $\omega>0$ if and only if $c_3=0$.
\item
	If (i) and (ii) of \cref{Z124} are satisfied, then 
	$|c_3|\le\frac{2\sqrt{\rho_1\rho_2}\,}{\rho_1+\rho_2}\cdot\sqrt{c_1c_2}$ by \eqref{Z97}.
\item
	\Cref{Z124} generalises \cite[Corollaries~3.6 and 3.9]{eckhardt.kostenko.teschl:2018}
	where the case $c_3=0$ is considered.
\item
	For $\alpha=0$ \cref{Z124} reduces as follows.
	Under the assumptions of \cref{Z124} the following statements are equivalent:
	\begin{Enumeratealph}
	\item
		$q_H(ri)\to\zeta_0$ as $r\to\infty$ with $\zeta_0\ne0$;
	\item
		the limits $c_i\DE\lim_{t\to0}\frac{m_i(t)}{t^\sigma}=c_i$ exist, $c_1,c_2>0$, and
		\[
			\zeta_0 = \frac{c_3}{c_2}+i\sqrt{\frac{c_1}{c_2}-\Bigl(\frac{c_3}{c_2}\Bigr)^2}.
		\]
	\end{Enumeratealph}
	This reproves \cite[Theorem~3.1]{eckhardt.kostenko.teschl:2018}.
\end{Enumerate}
\end{remark}
\noindent

\begin{example}\label{Z94}
	Assume that we are given a Hamiltonian $H=\smmatrix{h_1}{h_3}{h_3}{h_2}$ with
	\[
		h_1(t) \sim \kappa_1 t^{\rho_1-1}, \qquad
		h_2(t) \sim \kappa_2 t^{\rho_2-1}, \qquad
		h_3(t)
		\begin{cases}
			\sim \kappa_3 t^{\rho_3-1} & \text{if} \ \kappa_3\ne0,
			\\[1ex]
			\ll t^{\rho_3-1} & \text{if} \ \kappa_3=0,
		\end{cases}
	\]
	as $t\to0$, where
	\begin{Itemize}
	\item
		$\rho_1,\rho_2>0$, $\rho_3=\frac{\rho_1+\rho_2}{2}$,
	\item
		$\kappa_1,\kappa_2>0$, $\kappa_3\in\bb R$,
	\item
		$\kappa_3^2\le\kappa_1\kappa_2$.
	\end{Itemize}
	Then \eqref{Z135} is satisfied with $c_i=\frac{\kappa_i}{\rho_i}$,
	$i\in\{1,2,3\}$.  Hence \eqref{Z137} holds with $\alpha$ and $\omega$
	as in \eqref{Z136}.
	Moreover, $\arg\omega$ is, as a function of $\frac{\kappa_3}{\sqrt{\kappa_1\kappa_2}\,}$, 
	a decreasing and odd bijection from $[-1,1]$ onto $[-\frac\pi2(1-|\alpha|),\frac\pi2(1-|\alpha|)]$.
\end{example}

\medskip

\begin{remark}\label{Z191}
	One can find the asymptotic inverse of $m_1m_2$ (which is needed for finding $\ms a_H$)
	in more general situations than those considered in \cref{Z124}.
	Let $\rho>0$, let $\ms g$ be regularly varying at $0$ with index $\gamma\in\bb R$
	and assume that $\ms f$ is a strictly increasing function on $(0,L)$ that satisfies
	\[
		\ms f(t) \sim t^\rho\ms g\bigl(|\log t|\bigr) \qquad \text{as} \ t\to0.
	\]
	Then
	\begin{equation}\label{Z192}
		\ms f^{-1}(x) \sim \rho^{\frac{\gamma}{\rho}}\Bigl(\frac{x}{\ms g(|\log x|)}\Bigr)^{\frac{1}{\rho}}
		\qquad \text{as} \ x\to0.
	\end{equation}
	This covers cases where the $m_i$ are of the form of the functions in \eqref{Z193} with $\log r$
	replaced by $|\log r|$.
	Relation \eqref{Z192} can be proved by showing $(\ms f^{-1}\circ\ms f)(t)\sim t$ and 
	using \cite[Theorem~1.8.6]{bingham.goldie.teugels:1989};
	see also \cite[Remark~A.6]{pruckner.reiffenstein.woracek:sinqB-arXiv}.
\end{remark}

\medskip

\noindent
The following example deals with a situation where $m_2$ is rapidly varying.

\begin{example}\label{Z117}
	Consider a Hamiltonian $H$ such that
	\[
		m_1(t) \sim ct^\gamma, \quad m_2(t) \sim e^{-\frac{1}{t}}
		\qquad\text{as} \ t\to0,
	\]
	for some $\gamma,c>0$.  Clearly, $m_1$ is regularly varying with index $\rho_1=\gamma$
	and $m_2$ is rapidly varying with index $\rho_2=\infty$.
	Hence statement (ii) in \cref{Z10} is true.
	The latter theorem implies that \eqref{Z57} holds with $\alpha=1$ and $\omega_{\alpha,\delta}=1$;
	see \eqref{Z78} and \eqref{Z98}.
	In order to find the asymptotic behaviour of $\mr t$, we write
	\[
		m_1(t)m_2(t) = ct^\gamma e^{-\frac{1}{t}}e^{\eta(t)}, \qquad t>0,
	\]
	with some function $\eta:(0,\infty)\to\bb R$ such that $\lim_{t\to0}\eta(t)=0$.
	The defining relation \eqref{Z38} for $\mr t$ implies that
	\[
		\ln c + \gamma\ln\bigl(\mr t(r)\bigr) - \frac{1}{\mr t(r)} + \eta\bigl(\mr t(r)\bigr)
		= -2\ln r
	\]
	for $r>0$.  If $r$ is large enough so that $\mr t(r)\le1$, then
	\begin{align}
		\frac{1}{\mr t(r)} &= 2\ln r + \ln c + \gamma\ln\bigl(\mr t(r)\bigr)
		+ \eta\bigl(\mr t(r)\bigr)
		\label{Z118}
		\\[1ex]
		&\le 2\ln r + \ln c + \eta\bigl(\mr t(r)\bigr).
		\label{Z119}
	\end{align}
	Using \eqref{Z118} and \eqref{Z119} we obtain
	\begin{equation}\label{Z120}
		\frac{1}{\mr t(r)} \ge 2\ln r + \ln c
		- \gamma\ln\bigl[2\ln r + \ln c + \eta\bigl(\mr t(r)\bigr)\bigr]
		+ \eta\bigl(\mr t(r)\bigr).
	\end{equation}
	Since $\eta(\mr t(r))\to0$ as $r\to\infty$, the relations \eqref{Z119} and \eqref{Z120}
	imply that $\mr t(r)\sim\frac{1}{2\ln r}$.
	Together with \eqref{Z17}, this shows that
	\[
		\ms a_H(r) = rm_1\bigl(\mr t(r)\bigr) \sim \frac{cr}{(2\ln r)^\gamma}
		\qquad\text{as} \ r\to\infty.
	\]
	Now \eqref{Z57} yields
	\[
		q_H(rz) \sim \frac{crz}{(2\ln r)^\gamma} \qquad \text{as} \ r\to\infty
	\]
	locally uniformly for $z\in\bb C^+$.
\end{example}

\medskip

\noindent
Finally, we consider an example to illustrate the inverse result \cref{Z236}.

\begin{example}\label{Z249}
	Let $\phi\in\bigl[-\frac{\pi}{2},\frac{\pi}{2}\bigr]$ and suppose that $H$ is
	a Hamiltonian such that
	\begin{equation}\label{Z257}
		q_H(ri) \sim ie^{i\phi}\log r \qquad \text{as} \ r\to\infty.
	\end{equation}
	Set $\ms f(r)=\log r$ for $r>1$.
	Then $\alpha=\ind\ms f=0$, $\delta=-\sin\phi$ and $C_{\alpha,\delta}=1$ by \eqref{Z11} and \eqref{Z79}.
	Further, set $\ms g(r)=\frac{\log r}{r}$, $r>1$, which satisfies \eqref{Z246}.
	It is easy to check that
	\[
		\ms g^{-1}(t) \sim \frac{|\log t|}{t} \qquad \text{as} \ t\to0;
	\]
	cf.\ \cref{Z191}.  Hence \cref{Z236} implies that there exists a reparameterisation $\widetilde H$
	of $H$ such that
	\begin{equation}\label{Z209}
		\widetilde m_1(t) \sim t, \qquad
		\widetilde m_2(t) \sim \frac{t}{|\log t|^2}, \qquad
		\widetilde m_3(t)
		\begin{cases}
			\sim -\sin\phi\cdot\frac{t}{|\log t|} & \text{if} \ \phi\ne0,
			\\[1ex]
			\ll \frac{t}{|\log t|} & \text{if} \ \phi=0,
		\end{cases}
	\end{equation}
	as $t\to0$.  Conversely, if the relations in \eqref{Z209} hold for 
	a reparameterisation $\widetilde H$ of $H$, then $q_H$ satisfies \eqref{Z257}.
\end{example}

\section{Preliminaries}
\label{Z60}
\subsection{Hamiltonians and their Weyl coefficients}
\label{Z146}

We denote by $W(t;z)=\smmatrix{w_{11}(t;z)}{w_{12}(t;z)}{w_{21}(t;z)}{w_{22}(t;z)}$
the \emph{fundamental solution} of the canonical system \eqref{Z1},
i.e.\ the unique solution of the initial value problem
\begin{equation}\label{Z7}
 	\begin{cases}
		\frac{\partial}{\partial t} W(t;z) J = z W(t;z)H(t), \quad t \in [0,L), \\[1ex]
		W(0,z)=I.
	\end{cases}
\end{equation}
Note that the transposes of the rows of $W$ are solutions of \eqref{Z1}.
The image of the closed upper half-plane under the linear fractional transformation
\[
	\tau\mapsto \frac{w_{11}(t;z)\tau +  w_{12}(t;z)}{w_{21}(t;z)\tau + w_{22}(t;z)}
\]
is a disc in the upper half-plane. For fixed $z\in \bb C^+$, these discs are nested and
converge to a single point as $t\to L$ due to the limit-point assumption $\int_0^L \tr H(t)\RD t=\infty$
made at the beginning of the Introduction.
This limit is denoted by $q_H(z)$, which we state explicitly in the following definition.

\begin{definition}\label{Z153}
	The function
	\begin{equation}\label{Z111}
		q_H(z) \DE \lim_{t\to L}\frac{w_{11}(t;z)\tau +  w_{12}(t;z)}{w_{21}(t;z)\tau + w_{22}(t;z)}\,,
		\qquad z\in\bb C^+,
	\end{equation}
	with arbitrary $\tau\in\bb C^+\cup\bb R\cup\{\infty\}$ is called
	the \emph{Weyl coefficient} corresponding to the Hamiltonian~$H$.
\end{definition}

\noindent
Unless $h_2=0$ a.e.\ (in which case $q_H\equiv\infty$),
the Weyl coefficient is a \emph{Nevanlinna function}, i.e.\ it is analytic in the
open upper half-plane $\bb C^+$ and has non-negative imaginary part.

\begin{remark}\label{Z35}
	In some papers, e.g.\ \cite{eckhardt.kostenko.teschl:2018,romanov:1408.6022v1},
	the equation $Jy'(t)=zH(t)y(t)$ is considered instead of \eqref{Z1}.
	The corresponding Weyl coefficient is $\widetilde q_H(z)=-\qu{q_H(-\qu z)}$.
	The latter function is also the Weyl coefficient (according to our definition)
	of the Hamiltonian $\widehat H\DE\smmatrix{h_1}{-h_3}{-h_3}{h_2}$,
	i.e.\ $q_{\widehat H}(z)=-\qu{q_H(-\qu z)}=\widetilde q_H(z)$.
\end{remark}

\medskip

\noindent
We also need some properties of reparameterisations of Hamiltonians.
\begin{Itemize}
\item If \eqref{Z34} holds and $y$ is a solution of \eqref{Z1} with $H$ replaced by $H_2$,
	then $y\circ\gamma$ is a solution of \eqref{Z1} with $H$ replaced by $H_1$.
	Moreover, if \eqref{Z34} holds and $M_1$, $M_2$ denote the primitives
	of $H_1$, $H_2$, respectively, i.e.\ $M_i(t)=\int_0^t H_i(s)\RD s$, $t\in[0,L_i)$, then
	\begin{equation}\label{Z36}
		M_1(t) = M_2(\gamma(t)), \qquad t\in[0,L_1).
	\end{equation}
\item The Weyl coefficient is unchanged by reparameterisation, i.e.\ $q_{H_1}=q_{H_2}$
	if $H_1$ and $H_2$ are related by \cref{Z34}.
\item For each Hamiltonian $H$ on an interval $[0,L)$ there is a unique
	\emph{trace-normalised} reparameterisation of $H$ defined on $[0,\infty)$ and 
	denoted by $\tn{H}$, i.e.\ $\tr\tn{H}(t)=1$ for a.e.\ $t\in(0,\infty)$.
	In fact, the strictly increasing, bijective and locally absolutely continuous map
	\[
		\mf t: [0,L) \to [0,\infty), \quad x\mapsto \int_0^x \tr H(s)\RD s,
	\]
	can be used to define
	\[
		\tn{H}(t) \DE (\mf t^{-1})'(t)\cdot H\left(\mf t^{-1}(t)\right).
	\]
	Note that $\mf t^{-1}$ is locally absolutely continuous since $\mf t'(x)=\tr H(x)>0$
	almost everywhere (see, e.g.\ \cite[Exercise 13 on p.~271]{natanson:1955}).
\end{Itemize}
Let $\Ham_L$ denote the set of all Hamiltonians on the interval $[0,L)$ with $L\in(0,\infty]$,
set $\Ham\DE\Ham_\infty$,
and let $\NHam\subseteq\Ham$ be the subset of all trace-normalised Hamiltonians.

A deep theorem due to de~Branges, see \cite{debranges:1968}, states that
the map that assigns to each Hamiltonian $H$ its Weyl coefficient $q_H$
is a bijection from the set of all trace-normalised Hamiltonians, $\NHam$,
onto the set $\mc N\cup\{\infty\}$, where $\mc N$ denotes the set of
all Nevanlinna functions.

We equip $\NHam$ with the topology such that this mapping is a homeomorphism
when $\mc N$ carries the topology of locally uniform convergence in $\bb C^+$.
This topology can be defined intrinsically and the space $\NHam$ equipped with this topology 
is compact and metrisable; see, e.g.\ \cite[Lemma~2.9]{pruckner.woracek:limp}.
This topology on $\NHam$ is implicitly used in the work of de~Branges and is studied
more explicitly in, e.g.\ \cite[Chapter~5.2]{remling:2018} and \cite{pruckner.woracek:limp}.
In particular, for $H_n,H\in \NHam$ we have $\lim\limits_{n\to \infty}H_n =H$
if and only if
\begin{equation}\label{Z20}
	\lim_{n\to \infty}\int_0^x H_n(s) \RD s =\int_0^x H(s) \RD s \quad \text{for all} \ x>0.
\end{equation}
A standard compactness argument shows that this convergence is locally uniform
in $x\in(0,\infty)$;
cf.\ \cite[Remark~2.4]{eckhardt.kostenko.teschl:2018}.  For completeness, we provide a proof.

\begin{lemma}\label{Z22}
	Let $H_n,H\in \NHam$ be given and assume that $\lim\limits_{n\to \infty}H_n =H$.
	Then \eqref{Z20} holds locally uniformly for $x\in [0,\infty)$.
\end{lemma}

\begin{proof}
	Fix $T>0$ and consider the family of functions
	\[
		M_H:
		\begin{cases}
			[0,T]\!\!\!&\to \ \bb R^{2\times 2} \\[1ex]
			x & \mapsto \ \displaystyle\int_0^x H(s) \RD s,
		\end{cases}
		\qquad H\in \NHam.
	\]
	Due to $\|M_H(x)-M_H(x')\|_{\ell^1(\bb R^{2\times2})}\le 2 |x-x'|$,
	this family is pointwise bounded and uniformly equicontinuous.
	The Arzel\`a--Ascoli theorem implies pre-compactness in the space of all
	continuous matrix-valued functions on $[0,T]$.

	Now let $H_n,H$ satisfy \eqref{Z20} pointwise for every $x>0$.
	Then, each subsequence of $M_{H_n}$ has a subsubsequence that
	converges uniformly on $[0,T]$.
	Due to \eqref{Z20} the limit is $M_H$.
\end{proof}

\noindent
The notion of convergence on $\NHam$ is pulled back to $\Ham$ in the obvious way.

\begin{definition}\label{Z62}
	We say that a sequence $H_n\in \Ham$ converges to a Hamiltonian $H\in \Ham$
	if the sequence of trace-normalised Hamiltonians $\tn H_n$ converges to $\tn{H}$ in $\NHam$,
	i.e.\
	\[
		\lim_{n\to \infty}\int_0^x\tn H_n(s) \RD s = \int_0^x \tn{H}(s) \RD s, \quad x>0.
	\]
\end{definition}

\begin{lemma}\label{Z74}
	Let $H_n,H\in \Ham$ be given, and set
	\[
		\mf t_n(x) \DE \int_0^x \tr H_n(s)\RD s,\quad \mf t(x) \DE \int_0^x \tr H(s)\RD s,
		\qquad x\ge0,\, n\in\bb N.
	\]
	\begin{Enumerate}
	\item
		Then $\lim\limits_{n\to \infty}H_n =H$ is equivalent to
		\begin{equation}\label{Z23}
			\lim_{n\to \infty}\int_0^{\mf t_n^{-1}(x) } H_n(t) \RD t
			= \int_0^{\mf t^{-1}(x) } H(t) \RD t,
		\end{equation}
		for all $x\in[0,\infty)$.  Further, this is equivalent
		to $\eqref{Z23}$ locally uniformly for $x\in[0,\infty)$,
		and to $q_{H_n}\to q_H$ locally uniformly.
	\item 
		We have
		\begin{equation}\label{Z31}
			\lim_{n\to\infty}\int_0^x H_n(t)\RD t = \int_0^x H(t)\RD t
		\end{equation}
		locally uniformly for $x\in[0,\infty)$ if and only if
		$\lim\limits_{n\to\infty}H_n=H$ and $\lim\limits_{n\to\infty}\mf t_n(x)=\mf t(x)$
		locally uniformly for $x\in[0,\infty)$.
	\end{Enumerate}
\end{lemma}

\begin{proof}
	Let $\widehat H_n$ and $\widehat H$ be the trace-normalised reparameterisations
	of $H_n$ and $H$, respectively.
	Then
	\[
		\int_0^x \widehat H_n(s)\RD s = \int_0^{\mf t_n^{-1}(x)}H_n(t)\RD t, \qquad
		\int_0^x \widehat H(s)\RD s=\int_0^{\mf t^{-1}(x)}H(t)\RD t.
	\]
	The first equivalence in (i) follows. For the second one note \cref{Z22}.
	The last equivalence follows from the fact that the mapping $\NHam\to\mc N$,
	$H\mapsto q_H$ is a homeomorphism.
	
	For the proof of (ii), assume first that \eqref{Z31} holds. Then $\mf t_n\to\mf t$
	locally uniformly, in particular, pointwise.
	Since the functions $\mf t_n$ and $\mf t$ are increasing and continuous, it follows that
	$\lim_{n\to \infty} \mf t_n^{-1}(s) =\mf t^{-1}(s)$ for $s\in [0,\infty)$.
	The fact that the limit in \eqref{Z31} is assumed to be locally uniform in $x$
	implies that we can use it with $\mf t_n^{-1}(x)$ instead of $x$, which yields \eqref{Z23}
	pointwise and hence $H_n\to H$.
	Conversely, assume that $H_n\to H$ and $\mf t_n\to\mf t$ locally uniformly.
	Then \eqref{Z23} holds locally uniformly,
	and it follows that \eqref{Z31} holds locally uniformly.
\end{proof}

\noindent
The next lemma shows that the non-negativity of the limit Hamiltonian is automatic.

\begin{lemma}\label{Z112}
	Let $H_n\in \Ham$, $n\in\bb N$, and $H\in L^1_{\textup{loc}}([0,\infty),\bb R^{2\times2})$
	be given and assume that $\tr H(t)>0$ for a.e.\ $t\in(0,\infty)$ and that \eqref{Z20} holds.
	Then $H\in\Ham$.
\end{lemma}

\begin{proof}
	It follows from \eqref{Z20} and the relations $H_n(t)\ge0$ for a.e.\ $t\in(0,\infty)$
	that, for all $x_1,x_2\in[0,\infty)$ with $x_1<x_2$ and every $\xi\in\bb R^2$, we have
	\[
		\int_{x_1}^{x_2} \xi^T H(s)\xi\RD s \ge0.
	\]
	This implies that $H(t)\ge0$ for a.e.\ $t\in(0,\infty)$;
	hence $H\in\Ham$.
\end{proof}

\noindent
In later sections we shall also use the following inequality, which
follows from the non-negativity of $H$ and the Cauchy--Schwarz inequality:
\begin{equation}\label{Z103}
\begin{aligned}
	|m_3(t)| &\le \int_0^t |h_3(s)|\RD s
	\le \int_0^t \sqrt{h_1(s)h_2(s)}\,\RD s
	\\[1ex]
	&\le \biggl[\int_0^t h_1(s)\RD s\biggr]^{1/2}\biggl[\int_0^t h_2(s)\RD s\biggr]^{1/2}
	= \sqrt{m_1(t)m_2(t)}.
\end{aligned}
\end{equation}

\subsection{Regular Variation}

When dealing with regularly varying functions, we use the usual convention for algebra
in $[0,\infty]$: $0\cdot\infty\DE 0$; $\infty+x\DE\infty$, $\infty-x\DE\infty$ when $x\in\bb R$;
$\frac{x}{0}=\infty$ when $x>0$; $\frac{x}{\infty}=0$ when $x<\infty$.

Note the following simple consequence of the Potter bounds \cite[Theorem~1.5.6\,(iii)]{bingham.goldie.teugels:1989}.
If $\ms f$ is regularly varying at $\infty$ with index $\alpha$, then
\begin{equation}\label{Z27}
	\ms f(r) \gg r^{\alpha-\varepsilon}, \qquad
	\ms f(r) \ll r^{\alpha+\varepsilon}
	\qquad\text{as} \ r\to\infty,
\end{equation}
for every $\varepsilon>0$.
Analogously, if $\ms f$ is regularly varying at $0$ with index $\alpha$, then
\begin{equation}\label{Z18}
	\ms f(r) \ll r^{\alpha-\varepsilon}, \qquad
	\ms f(r) \gg r^{\alpha+\varepsilon}
	\qquad\text{as} \ r\to\infty,
\end{equation}
for every $\varepsilon>0$.

We need the following variant of \cite[Theorem~1.5.11]{bingham.goldie.teugels:1989}.

\begin{lemma}\label{Z122}
	Let $\ms f\in L_{\textup{\textsf{loc}}}^1([0,\infty))$ such that $\ms f(t)>0$ a.e.\
	and define
	\[
		F(t) \DE \int_0^t \ms f(s)\RD s, \qquad t>0.
	\]
	\begin{Enumerate}
	\item
		Assume that $\ms f\in R_\rho(0)$ with $\rho\in[-1,\infty]$.
		Then $F\in R_{\rho+1}(0)$.
		\\
		If, in addition, $\rho\in(-1,\infty)$, then
		\[
			F(t) \sim \frac{1}{\rho+1}t\ms f(t) \qquad\text{as} \ t\to0.
		\]
	\item
		Assume that $\ms f$ is non-decreasing on $[0,t_0]$ for some $t_0>0$
		and that $F\in R_\rho(0)$ with $\rho\in[0,\infty]$.
		Then $\ms f\in R_{\rho-1}(0)$.
		\\
		If, in addition, $\rho\in(0,\infty)$, then
		\[
			\ms f(t) \sim \rho\frac{F(t)}{t} \qquad\text{as} \ t\to0.
		\]
	\end{Enumerate}
\end{lemma}
\begin{proof}
(i)
	The statement for $\rho\in[-1,\infty)$ follows
	from \cite[Theorem~1.5.11 and Proposition~1.5.9a]{bingham.goldie.teugels:1989} by transforming
	the limit towards infinity to $0$.
	Now assume that $\rho=\infty$.  Let $\lambda\in(0,1)$ and $\varepsilon>0$.
	There exists $\delta>0$ such that $\frac{\ms f(\lambda s)}{\ms f(s)}<\varepsilon$
	for all $s\in(0,\delta)$.
	Then, for all $t\in(0,\delta)$, we have
	\[
		\frac{F(\lambda t)}{F(t)}
		= \lambda\frac{\int_0^t \ms f(\lambda s)\RD s}{\int_0^t \ms f(s)\RD s}
		\le \lambda\varepsilon \le \varepsilon,
	\]
	which shows that $F\in R_\infty(0)$.

(ii)
	The case $\rho\in[0,\infty)$ follows directly
	from \cite[Theorem~1.7.2b]{bingham.goldie.teugels:1989}.
	Now consider the case $\rho=\infty$.  Let $\lambda\in(0,1)$.  For $t\in(0,t_0]$ we 
	obtain from the monotonicity of $\ms f$ that
	\begin{align*}
		F(t)-F(\lambda t) &= \int_{\lambda t}^t \ms f(s)\RD s
		\le (1-\lambda)t\ms f(t),
		\\[1ex]
		F(\lambda t)-F(\lambda^2 t) &= \int_{\lambda^2 t}^{\lambda t} \ms f(s)\RD s
		\ge \lambda(1-\lambda)t\ms f(\lambda^2 t)
	\end{align*}
	and hence
	\[
		\frac{\ms f(\lambda^2 t)}{\ms f(t)}
		\le \frac{1}{\lambda}\cdot\frac{F(\lambda t)-F(\lambda^2 t)}{F(t)-F(\lambda t)}
		= \frac{1}{\lambda}\cdot
		\frac{\,\frac{F(\lambda t)}{F(t)}-\frac{F(\lambda^2t)}{F(t)}\,}{1-\frac{F(\lambda t)}{F(t)}}
		\to 0
	\]
	as $t\to0$, which implies that $\ms f\in R_\infty(0)$.
\end{proof}

\medskip

\noindent
A regularly varying function $\ms f$ with non-zero index is asymptotically invertible as the following
theorem shows.

\begin{theorem}\label{Z254}
	Let $\ms f$ be regularly varying at $\infty$ with index $\alpha$.
	\begin{Enumerate}
	\item
		Let $\alpha>0$.  Then there exists $\ms g$ regularly varying at $\infty$ with index $\frac{1}{\alpha}$
		such that
		\begin{equation}\label{Z258}
			\ms f(\ms g(x)) \sim \ms g(\ms f(x)) \sim x, \qquad x\to\infty.
		\end{equation}
		Let $\widetilde{\ms f}$ be a positive measurable function such that $\widetilde{\ms f}(x)\sim \ms f(x)$ as $x\to\infty$
		and let $\widetilde{\ms g}$ be such that \eqref{Z258} holds with $\ms f$ and $\ms g$ replaced
		by $\widetilde{\ms f}$ and $\widetilde{\ms g}$ respectively.
		Then $\widetilde{\ms g}(x)\sim \ms g(x)$ as $x\to\infty$.
	\item
		Let $\alpha<0$.  Then there exists $\ms g$ regularly varying at $0$ with index $\frac{1}{\alpha}$ such that
		\begin{alignat}{2}
			\ms f(\ms g(x)) &\sim x, \qquad && x\to0,
			\label{Z255}
			\\[0.5ex]
			\ms g(\ms f(x)) &\sim x, \qquad && x\to\infty.
			\label{Z256}
		\end{alignat}
		Let $\widetilde{\ms f}$ be a positive measurable function such that $\widetilde{\ms f}(x)\sim \ms f(x)$ as $x\to\infty$
		and let $\widetilde{\ms g}$ be such that \eqref{Z255} and \eqref{Z256} hold with $\ms f$ and $\ms g$ replaced
		by $\widetilde{\ms f}$ and $\widetilde{\ms g}$ respectively.
		Then $\widetilde{\ms g}(x)\sim \ms g(x)$ as $x\to0$.
	\end{Enumerate}
\end{theorem}

\noindent
A function $\ms g$ as in \cref{Z254} is called \emph{asymptotic inverse} of $\ms f$.

\begin{proof}
	The statement in (i) follows from \cite[Theorems~1.5.12 and 1.8.7]{bingham.goldie.teugels:1989}.
	Item (ii) follows by transformation from $\infty$ to $0$: let $\ms h(x)\DE\frac{1}{\ms f(x)}$,
	which has index $-\alpha>0$, and let $\check{\ms h}$ be an asymptotic inverse of $\ms h$ as in (i).  
	Then $\ms g(x)\DE\check{\ms h}\bigl(\frac{1}{x}\bigr)$ satisfies \eqref{Z255} and \eqref{Z256}.
	The second statement is shown in a similar way.
\end{proof}

\medskip

\noindent
We also need the notion of smooth variation; see, e.g.\ \cite[Section~1.8]{bingham.goldie.teugels:1989}.

\begin{definition}\label{Z250}
	A positive function $\ms f$ is called \emph{smoothly varying} at $\infty$ \textup{(}respectively at $0$\textup{)} 
	with index $\alpha$ if it is in $C^\infty$ and $\ms h(x)\DE\log\ms f(e^x)$ satisfies
	\begin{equation}\label{Z251}
		\ms h'(x) \to \alpha, \qquad \ms h^{(n)}(x) \to 0, \quad n\in\{2,3,\ldots\},
	\end{equation}
	as $x\to\infty$ \textup{(}respectively as $x\to-\infty$\textup{)}.
\end{definition}

\noindent
One can show that \eqref{Z251} is equivalent to 
\begin{equation}\label{Z252}
	\lim_{\substack{x\to\infty \\[0.2ex] (x\to0)}}\frac{x^n\ms f^{(n)}(x)}{\ms f(x)} = \alpha(\alpha-1)\cdots(\alpha-n+1),
	\qquad n\in\bb N;
\end{equation}
see \cite[(1.8.1')]{bingham.goldie.teugels:1989}.

The next theorem shows that we can often 
assume, without loss of generality, that a regularly varying function is smoothly varying.

\begin{theorem}[Smooth Variation Theorem]\label{Z253}
	Let $\ms f:\mc I\to(0,\infty)$ be regularly varying at $\infty$ or at $0$ with index $\alpha$, 
	where $\mc I\subseteq(0,\infty)$ is an interval.
	Then there exist smoothly varying functions $\ms g_1,\ms g_2$ such that $\ms g_1(x)\le \ms f(x) \le \ms g_2(x)$
	for all $x\in\mc I$ and $\ms g_1(x)\sim \ms g_2(x)$ as $x\to\infty$ \textup{(}as $x\to0$ respectively\textup{)}.
	
	If $\alpha\ne0$, then there exists a strictly monotone and smoothly varying function $\ms g$ 
	\textup{(}strictly increasing if $\alpha>0$ and strictly decreasing if $\alpha<0$\textup{)}
	such that $\ms f(x)\sim \ms g(x)$ as $x\to\infty$ \textup{(}as $x\to0$ respectively\textup{)}.
\end{theorem}

\begin{proof}
	For regular variation at $\infty$ the existence of smoothly varying functions $\ms g_1,\ms g_2,\ms g$ 
	follows from \cite[Theorem~1.8.2]{bingham.goldie.teugels:1989}.
	Now let $\alpha\ne0$.  The relation in \eqref{Z252} for $n=1$ implies that $\ms g$ is strictly increasing or decreasing
	in a neighbourhood of $\infty$.  On the remaining interval one can change $\ms g$ such
	that it is strictly increasing or decreasing everywhere.
	For regular variation at $0$ one uses a transformation.
\end{proof}

\subsection{Rescaling Hamiltonians}
\label{Z128}

We use a symmetrised variant of the transformations of Hamiltonians
that are used in \cite[Lemma~2.7 and (3.9)]{eckhardt.kostenko.teschl:2018};
for a similar transformation see \cite[Definition~2.10]{langer.pruckner.woracek:gapsatz}.

\begin{definition}\label{Z2}
	For $r,b_1,b_2>0$ we define a map $\mc A_r^{b_1,b_2}\DF\Ham\to\Ham$ by
	\[
		\bigl(\mc A_r^{b_1,b_2}H\bigr)(t)
		\DE
		\begin{pmatrix}
			b_1^2\cdot h_1\bigl(\frac{b_1b_2}rt\bigr) & b_1b_2\cdot h_3\bigl(\frac{b_1b_2}rt\bigr)
			\\[2ex]
			b_1b_2\cdot h_3\bigl(\frac{b_1b_2}rt\bigr) & b_2^2\cdot h_2\bigl(\frac{b_1b_2}rt\bigr)
		\end{pmatrix},
		\qquad t\in(0,\infty),
	\]
	where, as usual, $H=\smmatrix{h_1}{h_3}{h_3}{h_2}$.

	For $\ms b_1,\ms b_2\DF(0,\infty)\to(0,\infty)$, we set (slightly overloading notation)
	\[
		\mc A_r^{\ms b_1,\ms b_2}H\DE \mc A_r^{\ms b_1(r),\ms b_2(r)}H,
		\qquad r>0,\, H\in\Ham.
	\]
\end{definition}

\begin{remark}
Two transforms $\mc A_r^{b_1,b_2}H$ and $\mc A_r^{b_1',b_2'}H$ are reparameterisations
of each other whenever $\frac{b_1}{b_2}=\frac{b_1'}{b_2'}$.
In fact, if this equality holds, we have
\[
	\mc A_r^{b_1',b_2'}H(t)=\big[\mc A_r^{b_1,b_2}H\circ\gamma(t)\big]\cdot\gamma'(t)
\]
with $\gamma(t)\DE\frac{b_1'}{b_1}t$. Still, it turns out to be practical
to keep the two independent parameters $b_1,b_2$.
\end{remark}

An elementary calculation similar to the one in \cite[Lemma~2.7]{eckhardt.kostenko.teschl:2018} 
shows that the fundamental solution $W_{\mc A_r^{b_1,b_2}H}$ corresponding to $\mc A_r^{b_1,b_2}H$
is given by
\[
	W_{\mc A_r^{b_1,b_2}H}(t;z) =
	\begin{pmatrix}
		w_{11}\bigl(\tfrac{b_1b_2}{r}t;rz\bigr) & \frac{b_1}{b_2}w_{12}\bigl(\tfrac{b_1b_2}{r}t;rz\bigr)
		\\[2ex]
		\frac{b_2}{b_1}w_{21}\bigl(\tfrac{b_1b_2}{r}t;rz\bigr) & w_{22}\bigl(\tfrac{b_1b_2}{r}t;rz\bigr)
	\end{pmatrix},
\]
and hence
\begin{equation}\label{Z39}
	q_{\mc A_r^{b_1,b_2}H}(z) = \frac{b_1}{b_2}q_H(rz),
	\qquad z\in\bb C^+,\, r,b_1,b_2>0.
\end{equation}
The asymptotics of the Weyl coefficient is related to the convergence of rescaled Hamiltonians.
Indeed, the following lemma shows that the rescaling transformation $\mc A_r^{\ms b_1,\ms b_2}$
is an appropriate tool to study the behaviour of $q_H$ towards $i\infty$.

\begin{lemma}\label{Z14}
	Let $H,\mr H\in\Ham$ and let $\ms b_1,\ms b_2\DF(0,\infty)\to(0,\infty)$.
	Then the following statements are equivalent:
	\begin{Enumerate}
	\item
		the asymptotic relation
		\[
			q_H(rz) \sim \frac{\ms b_2(r)}{\ms b_1(r)}q_{\mr H}(z), \qquad r\to\infty,
		\]
		holds locally uniformly for $z\in\bb C^+$;
	\item
		$\displaystyle\lim_{r\to\infty}\mc A_r^{\ms b_1,\ms b_2}H=\mr H$.
	\end{Enumerate}
\end{lemma}

\begin{proof}
	It follows from \eqref{Z39} that (i) is equivalent to
	\[
		q_{\mr H}(z) 
		= \lim_{r\to\infty}\frac{\ms b_1(r)}{\ms b_2(r)}q_H(rz)
		= \lim_{r\to\infty}q_{\mc A_r^{\ms b_1,\ms b_2}H}(z)
	\]
	locally uniformly for $z\in\bb C^+$.
	By \cref{Z74}\,(i) the latter is equivalent to (ii).
\end{proof}

\noindent
For later reference note the following lemma.

\begin{lemma} \label{Z19}
	Let $H\in \Ham$ and assume that $\ms b_1,\ms b_2\DF(0,\infty)\to(0,\infty)$ are continuous.
	For $r>0$ define

	\[
		\mf t_r: [0,\infty) \to [0,\infty), \quad
		t\mapsto \ \int_0^t \tr\bigl(\mc A_r^{\ms b_1,\ms b_2}H\bigr)(s) \RD s,
	\]
	and fix $T>0$. Then $r\mapsto \mf t_r^{-1}(T)$ is continuous on $(0,\infty)$.
\end{lemma}

\begin{proof}
	Obviously $t=\mf t_r^{-1}(T)$ is the unique solution of $F(t,r) \DE T-\mf t_r(t)=0$.
	Note that $F$ is continuous in $t$ and $r$, and $t\mapsto F(t,r)$ is one-to-one
	for all $r>0$.
	Now we obtain the result from an application of a variant of the
	implicit function theorem; see \cite{kumagai:1980}.
\end{proof}

\section{The model situation: $\bm{q_H}$ is a power}
\label{Z61}

In this section we prove a direct and inverse spectral theorem for power functions.

Recall that we use the branch of the complex power which is analytic
in $\bb C\setminus(-\infty,0]$ and takes the value $1$ at $1$.
Correspondingly, we understand the argument $\arg z\in(-\pi,\pi)$
for $z\in\bb C\setminus(-\infty,0]$.

\begin{definition}\label{Z131}
	\phantom{}
	\begin{Enumerate}
	\item
		For $\omega\in\bb C$ and $\alpha\in\bb R$ set
		\begin{equation}\label{Z16}
			Q_{\alpha,\omega}(z) \DE i\omega\Big(\frac zi\Big)^\alpha,
			\qquad z\in\bb C^+.
		\end{equation}
	\item
		For $\uprho=(\rho_1,\rho_2,\rho_3),\upkappa=(\kappa_1,\kappa_2,\kappa_3)\in\bb R^3$ set
		\begin{equation}\label{Z45}
			H_{\uprho,\upkappa}(t) \DE
			\begin{pmatrix}
				\kappa_1 t^{\rho_1-1} & \kappa_3 t^{\rho_3-1}
				\\[1ex]
				\kappa_3 t^{\rho_3-1} & \kappa_2 t^{\rho_2-1}
			\end{pmatrix}
			,\qquad t\in(0,\infty).
		\end{equation}
	\end{Enumerate}
\end{definition}

\begin{remark}\label{Z15}
	\phantom{}
	\begin{Enumerate}
	\item
		We have $Q_{\alpha,\omega}\in\mc N$ if and only if
		\begin{Itemize}
		\item $\omega=0$ or
		\item $\alpha\in[-1,1]$ and $|\arg\omega|\le\frac\pi2(1-|\alpha|)$.
		\end{Itemize}
		This is seen by checking onto which sector the upper half-plane is mapped.
		In fact, when $\omega\ne0$, the function $Q_{\alpha,\omega}$ has a continuous extension
		to $(\bb C^+\cup\bb R)\setminus\{0\}$ with
		\begin{equation}\label{Z58}
			\arg Q_{\alpha,\omega}(x)
			=
			\begin{cases}
				\arg\omega+\frac\pi2(1-\alpha) &\text{if}\ x>0,
				\\[0.5ex]
				\arg\omega+\frac\pi2(1+\alpha) &\text{if}\ x<0.
			\end{cases}
		\end{equation}
		Since, for $r>0$, the function $\phi\mapsto\arg Q_{\alpha,\omega}(re^{i\phi})$ is increasing if $\alpha>0$
		and decreasing if $\alpha<0$, we have $Q_{\alpha,\omega}\in\mc N$ 
		if and only if $\arg\omega+\frac{\pi}{2}(1\pm\alpha)\in[0,\pi]$.
	\item
		We have $H_{\uprho,\upkappa}\in\bb H$ if and only if
		\begin{Itemize}
		\item
			$\kappa_1,\rho_1>0$ and $\kappa_2=\kappa_3=0$, or
		\item
			$\kappa_2,\rho_2>0$ and $\kappa_1=\kappa_3=0$, or
		\item
			$\kappa_1,\rho_1,\kappa_2,\rho_2>0$ and $\kappa_3^2\le\kappa_1\kappa_2$,
			such that $\rho_3=\frac 12(\rho_1+\rho_2)$ if $\kappa_3\ne 0$.
		\end{Itemize}
		This is seen by checking $H_{\uprho,\upkappa}(t)\ge 0$ for $t=1$, $t\to 0$,
		and $t\to\infty$, and checking integrability at $0$.
	\end{Enumerate}
\end{remark}

\noindent
With exception of some boundary cases, power Nevanlinna functions $Q_{\alpha,\omega}$
correspond to power Hamiltonians $H_{\uprho,\upkappa}$.
This is proved in \cref{Z132} below.  Before we state this result, let us settle the
mentioned boundary cases.

\begin{remark}\label{Z133}
	\phantom{}
	\begin{Itemize}
	\item
		If $\kappa_1,\rho_1>0$ and $\kappa_2=\kappa_3=0$,
		then $q_{H_{\uprho,\upkappa}}=\infty$.
	\item
		If $\kappa_2,\rho_2>0$ and $\kappa_1=\kappa_3=0$,
		then $q_{H_{\uprho,\upkappa}}=0$.
	\item
		If $\alpha=1$ and $\omega>0$, then $Q_{\alpha,\omega}$
		is the Weyl coefficient of the Hamiltonian
		\[
			H(t)\DE
			\begin{pmatrix}
				\mathds{1}_{(0,\omega]}(t) & 0
				\\[1ex]
				0 & \mathds{1}_{(\omega,\infty)}(t)
			\end{pmatrix}
			,\qquad t\in(0,\infty),
		\]
	\item
		If $\alpha=-1$ and $\omega>0$,
		then $Q_{\alpha,\omega}$ is the Weyl coefficient of the Hamiltonian
		\[
			H(t)\DE
			\begin{pmatrix}
				\mathds{1}_{(\frac 1\omega,\infty)}(t) & 0
				\\[1ex]
				0 & \mathds{1}_{(0,\frac 1\omega]}(t)
			\end{pmatrix}
			,\qquad t\in(0,\infty).
		\]
	\end{Itemize}
\end{remark}

\noindent
Now we formulate the main result of this section.

\begin{theorem}\label{Z132}
	The Weyl coefficient map $H\mapsto q_H$ induces a surjection
	\[
		\bigl\{H_{\uprho,\upkappa}\in\bb H\DS \kappa_1,\kappa_2>0\bigr\}
		\to
		\bigl\{Q_{\alpha,\omega}\in\mc N\DS \omega\ne 0,|\alpha|<1\bigr\}.
	\]
	The data $\uprho,\upkappa$ and $\alpha,\omega$ are related via
	\begin{equation}\label{Z102}
		\alpha = \frac{\rho_2-\rho_1}{\rho_2+\rho_1}
	\end{equation}
	and \textup{(}with $\kappa\DE\sqrt{\kappa_1\kappa_2-\kappa_3^2}$\textup{)}
	\begin{equation}\label{Z134}
		\omega
		= \begin{cases}
		\displaystyle
			\frac{(2\kappa)^{1+\alpha}}{\kappa_2\rho_3^\alpha}
			\cdot\frac{\Gamma(-\alpha)}{\Gamma(1+\alpha)}
			\cdot\frac{\Gamma\bigl(1+\frac{\alpha}{2}\bigl(1+i\frac{\kappa_3}{\kappa}\bigr)\bigr)}{%
			\Gamma\bigl(-\frac{\alpha}{2}\bigl(1-i\frac{\kappa_3}{\kappa}\bigr)\bigr)}\,,
			& \rho_1\ne\rho_2, \ \kappa>0,
			\\[3ex]
			\displaystyle
			\frac{(i\alpha\kappa_3)^{1+\alpha}}{\kappa_2\rho_3^\alpha}
			\cdot\frac{\Gamma(-\alpha)}{\Gamma(1+\alpha)}\,,
			& \rho_1\ne\rho_2, \ \kappa=0,
			\\[3ex]
			\displaystyle
			\frac{\kappa-i\kappa_3}{\kappa_2}\,,
			& \rho_1=\rho_2.
		\end{cases}
	\end{equation}
\end{theorem}

\noindent
We prove this result in several steps formulated as separate lemmas.

\begin{lemma}\label{Z40}
	Assume that $h_2(t)\ne0$, $t\in(0,\infty)$ a.e., and
	that $\frac{h_3}{h_2}$ is differentiable on $(0,\infty)$.
	Further, let $y=\binom{y_1}{y_2}$ be a solution of \eqref{Z1}.
	Then $y_1$ satisfies the differential equation
	\begin{equation}\label{Z41}
		\Bigl(\frac{1}{h_2}y_1'\Bigr)' + \biggl[z\Bigl(\frac{h_3}{h_2}\Bigr)'
		+ z^2\Bigl(h_1-\frac{h_3^2}{h_2}\Bigr)\biggr]y_1 = 0
	\end{equation}
	and, for $z\ne0$, the functions $y_1$ and $y_2$ are related by
	\begin{equation}\label{Z42}
		y_2 = -\frac{1}{h_2}\Bigl(\frac{1}{z}y_1'+h_3y_1\Bigr).
	\end{equation}
\end{lemma}

\begin{proof}
	Equation \eqref{Z1} is equivalent to the two scalar equations
	\begin{align*}
		y_1' &= z(-h_3y_1-h_2y_2),
		\\[0.5ex]
		y_2' &= z(h_1y_1+h_3y_2).
	\end{align*}
	We can solve the first equation for $y_2$, which yields \eqref{Z42}.
	Plugging this expression into the second equation we obtain \eqref{Z41}.
\end{proof}

\noindent
The next lemma cam be deduced from \cite[2.273\,(12)]{kamke:1961}.
In the latter reference the solution is written in terms of Whittaker functions,
which can be expressed in terms of Kummer functions.
For the convenience of the reader we give a direct proof.

\begin{lemma}\label{Z43}
	Let $c,d\in\bb C$ with $d\ne0$ and let $\gamma>0$ with $\gamma\ne1$.
	Two linearly independent solutions of the differential equation
	\begin{equation}\label{Z113}
		u''(x) + \bigl(cx^{\gamma-2}-d^2x^{2\gamma-2}\bigr)u(x) = 0, \qquad x>0,
	\end{equation}
	are given by
	\begin{align*}
		u_+(x) &= x\exp\Bigl(-\frac{d}{\gamma}x^\gamma\Bigr)
		M\biggl(\frac{\gamma+1-\frac{c}{d}}{2\gamma},\frac{\gamma+1}{\gamma},
		\frac{2d}{\gamma}x^\gamma\biggr),
		\\[1ex]
		u_-(x) &= \exp\Bigl(-\frac{d}{\gamma}x^\gamma\Bigr)
		M\biggl(\frac{\gamma-1-\frac{c}{d}}{2\gamma},\frac{\gamma-1}{\gamma},
		\frac{2d}{\gamma}x^\gamma\biggr),
	\end{align*}
	where $M$ is Kummer's confluent hypergeometric function,
	\[
		M(a,b,x) = \Fhyperg{1}{1}(a;b;x) = \sum_{n=0}^\infty \frac{(a)_n}{(b)_n}\cdot\frac{x^n}{n!}
	\]
	with the Pochhammer symbol $(a)_0=1$ and $(a)_n=a(a+1)\cdots(a+n-1)$ for $n\ge1$.
	Further,
	\begin{alignat*}{2}
		u_+(x) &= x - \frac{c}{\gamma(\gamma+1)}x^{\gamma+1} + \BigO\bigl(x^{2\gamma+1}\bigr),
		\quad &
		u_+'(x) &= 1 - \frac{c}{\gamma}x^\gamma + \BigO\bigl(x^{2\gamma}\bigr),
		\\[1ex]
		u_-(x) &= 1 - \frac{c}{(\gamma-1)\gamma}x^\gamma + \BigO\bigl(x^{2\gamma}\bigr),
		\quad &
		u_-'(x) &= -\frac{c}{\gamma-1}x^{\gamma-1} + \BigO\bigl(x^{2\gamma-1}\bigr)
	\end{alignat*}
	as $x\to 0$.
\end{lemma}

\begin{proof}
	We write
	\[
		u(x) = \exp\biggl(-\frac{d}{\gamma}x^\gamma\biggr)y(x).
	\]
	with some new unknown function $y$.
	Since
	\begin{align*}
		u'(x) &= \exp\biggl(-\frac{d}{\gamma}x^\gamma\biggr)y'(x)
		-dx^{\gamma-1}\exp\biggl(-\frac{d}{\gamma}x^\gamma\biggr)y(x)
		\\[1ex]
		u''(x) &= \exp\biggl(-\frac{d}{\gamma}x^\gamma\biggr)y''(x)
		- 2dx^{\gamma-1}\exp\biggl(-\frac{d}{\gamma}x^\gamma\biggr)y'(x)
		\\[0.5ex]
		&\quad + \Bigl(d^2x^{2\gamma-2}-d(\gamma-1)x^{\gamma-2}\Bigr)
		\exp\biggl(-\frac{d}{\gamma}x^\gamma\biggr)y(x),
	\end{align*}
	the differential equation \eqref{Z113} is equivalent to
	\begin{equation}\label{Z114}
		y''(x)-2dx^{\gamma-1}y'(x)+\bigl(c-d(\gamma-1)\bigr)x^{\gamma-2}y(x) = 0.
	\end{equation}
	With $\alpha=0$ or $\alpha=1$ we use
	\[
		y(x) = x^\alpha v\Bigl(\frac{2d}{\gamma}x^\gamma\Bigr).
	\]
	with a new unknown function $v$.
	The first two derivatives of $y$ are
	\begin{align*}
		y'(x) &= 2dx^{\alpha+\gamma-1}v'\Bigl(\frac{2d}{\gamma}x^\gamma\Bigr)
		+ \alpha x^{\alpha-1}v\Bigl(\frac{2d}{\gamma}x^\gamma\Bigr),
		\\[1ex]
		y''(x) &= 4d^2x^{\alpha+2\gamma-2}v''\Bigl(\frac{2d}{\gamma}x^\gamma\Bigr)
		+2d(\alpha+\gamma-1)x^{\alpha+\gamma-2}v'\Bigl(\frac{2d}{\gamma}x^\gamma\Bigr)
		\\[0.5ex]
		&\quad + 2d\alpha x^{\alpha+\gamma-2}v'\Bigl(\frac{2d}{\gamma}x^\gamma\Bigr)
		+ \underbrace{\frac{\RD}{\RD x}\bigl(\alpha x^{\alpha-1}\bigr)}_{=0}
		v\Bigl(\frac{2d}{\gamma}x^\gamma\Bigr)
		\\[1ex]
		&= 4d^2x^{\alpha+2\gamma-2}v''\Bigl(\frac{2d}{\gamma}x^\gamma\Bigr)
		+ 2d(2\alpha+\gamma-1)x^{\alpha+\gamma-2}v'\Bigl(\frac{2d}{\gamma}x^\gamma\Bigr).
	\end{align*}
	Hence \eqref{Z114} is equivalent to
	\begin{align*}
		& 4d^2x^{\alpha+2\gamma-2}v''\Bigl(\frac{2d}{\gamma}x^\gamma\Bigr)
		+ 2d(2\alpha+\gamma-1)x^{\alpha+\gamma-2}v'\Bigl(\frac{2d}{\gamma}x^\gamma\Bigr)
		\\[0.5ex]
		& - 4d^2x^{\alpha+2\gamma-2}v'\Bigl(\frac{2d}{\gamma}x^\gamma\Bigr)
		- 2d\alpha x^{\alpha+\gamma-2}v\Bigl(\frac{2d}{\gamma}x^\gamma\Bigr)
		\\[0.5ex]
		& + \bigl(c-d(\gamma-1)\bigr)x^{\alpha+\gamma-2}v\Bigl(\frac{2d}{\gamma}x^\gamma\Bigr)
		= 0;
	\end{align*}
	dividing by $x^{\alpha+\gamma-2}$ we obtain
	\[
		4d^2x^\gamma v''\Bigl(\frac{2d}{\gamma}x^\gamma\Bigr)
		+ 2d\bigl(2\alpha+\gamma-1-2dx^\gamma\bigr)v'\Bigl(\frac{2d}{\gamma}x^\gamma\Bigr)
		- \bigl(2d\alpha+d(\gamma-1)-c\bigr)v\Bigl(\frac{2d}{\gamma}x^\gamma\Bigr) = 0.
	\]
	Setting $t=\frac{2d}{\gamma}x^\gamma$ we see that this is equivalent to
	\[
		2d\gamma tv''(t) + 2d\bigl(2\alpha+\gamma-1-\gamma t\bigr)v'(t)
		- \bigl(2d\alpha+d(\gamma-1)-c\bigr)v(t) = 0,
	\]
	which, in turn, is equivalent to
	\[
		tv''(t) + \biggl(\frac{2\alpha+\gamma-1}{\gamma}-t\biggr)v'(t)
		- \biggl(\frac{2\alpha+\gamma-1}{2\gamma}-\frac{c}{2d\gamma}\biggr)v(t) = 0.
	\]
	This is Kummer's equation; a solution is given by
	\[
		v(t) = M\biggl(\frac{2\alpha+\gamma-1}{2\gamma}-\frac{c}{2d\gamma},\,
		\frac{2\alpha+\gamma-1}{\gamma},\, t\biggr);
	\]
	see, e.g.\ \cite[(13.2.1) and (13.2.2)]{nist:2010}.
	Substituting back we obtain the solutions $u_+$ and $u_-$
	for the choices $\alpha=1$ and $\alpha=0$ respectively.

	The asymptotic behaviour of the solutions can be obtained from
	\begin{align*}
		u(x) &= x^\alpha\exp\Bigl(-\frac{d}{\gamma}x^\gamma\Bigr)
		M\biggl(\frac{2\alpha+\gamma-1}{2\gamma}-\frac{c}{2d\gamma},\,
		\frac{2\alpha+\gamma-1}{\gamma},\, \frac{2d}{\gamma}x^\gamma\biggr)
		\\[1ex]
		&= x^\alpha\biggl(1-\frac{d}{\gamma}x^\gamma+\ldots\biggr)
		\biggl[1+\biggl(\frac{1}{2}-\frac{c}{2d(2\alpha+\gamma-1)}\biggr)\frac{2d}{\gamma}x^\gamma
		+ \ldots\biggr]
		\\[1ex]
		&= x^\alpha - \frac{c}{\gamma(2\alpha+\gamma-1)}x^{\alpha+\gamma} + \ldots
	\end{align*}
	and differentiation.
\end{proof}

\begin{lemma}\label{Z44}
	Let $\rho_1,\rho_2>0$ with $\rho_1\ne\rho_2$, and
	let $\kappa_1,\kappa_2>0$ and $\kappa_3\in\bb R$ be such that
	$\kappa_3^2<\kappa_1\kappa_2$.
	Set $\rho_3\DE\frac 12(\rho_1+\rho_2)$ and $\kappa\DE\sqrt{\kappa_1\kappa_2-\kappa_3^2}$, and consider the Hamiltonian
	$H_{\uprho,\upkappa}(t)$ where $\uprho=(\rho_1,\rho_2,\rho_3)$, $\upkappa=(\kappa_1,\kappa_2,\kappa_3)$.
	The entries $w_{i1}$, $i=1,2$, of the corresponding fundamental solution $W$,
	defined in \eqref{Z7}, are given by
	\begin{equation}\label{Z49}
	\begin{aligned}
		w_{11}(x,z) &= \exp\Bigl(\frac{\kappa iz}{\rho_3}x^{\rho_3}\Bigr)
		M\Bigl(a_-,\,b_-,\,-\frac{2\kappa iz}{\rho_3}x^{\rho_3}\Bigr),
		\\[1ex]
		w_{21}(x,z) &= -\frac{\kappa_2}{\rho_2}zx^{\rho_2}
		\exp\Bigl(\frac{\kappa iz}{\rho_3}x^{\rho_3}\Bigr)
		M\Bigl(a_+,\,b_+,\,-\frac{2\kappa iz}{\rho_3}x^{\rho_3}\Bigr),
	\end{aligned}
	\end{equation}
	where
	\begin{equation}\label{Z47}
		a_\pm = \frac{1}{2}\biggl(1\pm\frac{\rho_2}{\rho_3}
		-\frac{i\kappa_3}{\kappa}\cdot\frac{\rho_1-\rho_2}{\rho_1+\rho_2}\biggr),
		\qquad
		b_\pm = 1 \pm \frac{\rho_2}{\rho_3}\,.
	\end{equation}
\end{lemma}

\begin{proof}
	To shorten notation, we skip indices and write $H\equiv H_{\uprho,\upkappa}$.
	We multiply both sides of \eqref{Z41} by $\kappa_2$, which gives
	\begin{equation}\label{Z46}
		\bigl(x^{-\rho_2+1}y_1'(x)\bigr)'
		+ \Bigl[z\kappa_3\Bigl(x^{\frac{\rho_1-\rho_2}{2}}\Bigr)'
		+ z^2\bigl(\kappa_1\kappa_2-\kappa_3^2\bigr)x^{\rho_1-1}\Bigr]y_1(x) = 0.
	\end{equation}
	Define the function $u$ by $u(x^{\rho_2}) \DE y_1(x)$.  Then \eqref{Z46} is
	equivalent to
	\[
		\rho_2^2 x^{\rho_2-1}u''(x^{\rho_2})
		+ \biggl[z\kappa_3\frac{\rho_1-\rho_2}{2}x^{\frac{\rho_1-\rho_2}{2}-1}
		+ z^2\bigl(\kappa_1\kappa_2-\kappa_3^2\bigr)x^{\rho_1-1}\biggr]u(x^{\rho_2})
		= 0.
	\]
	Divide both sides by $\rho_2^2x^{\rho_2-1}$, which yields
	\[
		u''(x^{\rho_2})
		+ \biggl[z\kappa_3\frac{\rho_1-\rho_2}{2\rho_2^2}x^{\frac{\rho_1+\rho_2}{2}-2\rho_2}
		+ z^2\frac{\kappa_1\kappa_2-\kappa_3^2}{\rho_2^2}x^{\rho_1+\rho_2-2\rho_2}\biggr]u(x^{\rho_2})
		= 0.
	\]
	Setting $t=x^{\rho_2}$ we obtain
	\[
		u''(t) + \biggl[z\kappa_3\frac{\rho_1-\rho_2}{2\rho_2^2}t^{\frac{\rho_1+\rho_2}{2\rho_2}-2}
		+ z^2\frac{\kappa^2}{\rho_2^2}t^{\frac{\rho_1+\rho_2}{\rho_2}-2}\biggr]
		u(t) = 0.
	\]
	For $z\ne0$ we can use \cref{Z43} with
	\[
		c \DE z\kappa_3\frac{\rho_1-\rho_2}{2\rho_2^2}\,, \qquad
		d \DE -\frac{\kappa iz}{\rho_2}\,, \qquad
		\gamma \DE \frac{\rho_1+\rho_2}{2\rho_2} = \frac{\rho_3}{\rho_2}
	\]
	to obtain two linearly independent solutions,
	\begin{align*}
		u_\pm(x) &= x^{\frac{1\pm1}{2}}\exp\Bigl(-\frac{d}{\gamma}x^\gamma\Bigr)
		M\biggl(\frac{1}{2}\Bigl(1\pm\frac{1}{\gamma}-\frac{c}{\gamma d}\Bigr),\,
		1\pm\frac{1}{\gamma},\, \frac{2d}{\gamma}x^\gamma\biggr)
		\\[1ex]
		&= x^{\frac{1\pm1}{2}}\exp\biggl(\frac{\kappa iz}{\rho_3}
		x^{\frac{\rho_3}{\rho_2}}\biggr)
		M\biggl(a_\pm,\,b_\pm,\,-\frac{2\kappa iz}{\rho_3}
		x^{\frac{\rho_3}{\rho_2}}\biggr)
	\end{align*}
	with $a_\pm,b_\pm$ as in \eqref{Z47}.
	Set
	\[
		y_{\pm,1}(x) \DE u_\pm(x^{\rho_2}), \qquad
		y_{\pm,2}(x) \DE -\frac{1}{\kappa_2}\biggl(\frac{1}{z}x^{-\rho_2+1}y_{\pm,1}'(x)
		+ \kappa_3x^{\frac{\rho_1-\rho_2}{2}}y_{\pm,1}(x)\biggr)
	\]
	according to \eqref{Z42}.  Then $\binom{y_{\pm,1}}{y_{\pm,2}}$ are linearly
	independent solutions of \eqref{Z1} and hence
	\[
		w_{ij}(x,z) = A_i y_{+,j}(x) + B_i y_{-,j}(x),
		\qquad i,j=1,2,
	\]
	with some constants $A_i$ and $B_i$, $i=1,2$.
	To determine the constants $A_i$, $B_i$, we have to study the behaviour
	of $y_{\pm,j}$ at $0$.  \Cref{Z43} implies that
	\begin{alignat*}{2}
		y_{+,1}(x) &= x^{\rho_2} + \BigO\bigl(x^{\frac{\rho_1+3\rho_2}{2}}\bigr),
		\qquad &
		y_{+,1}'(x) &= \rho_2x^{\rho_2-1} + \BigO\bigl(x^{\frac{\rho_1+3\rho_2}{2}-1}\bigr),
		\\[1ex]
		y_{-,1}(x) &= 1 - \frac{\kappa_3z}{\rho_3}x^{\rho_3}
		+ \BigO\bigl(x^{\rho_1+\rho_2}\bigr),
		\qquad &
		y_{-,1}'(x) &= -\kappa_3zx^{\rho_3-1} + \BigO\bigl(x^{\rho_1+\rho_2-1}\bigr)
	\end{alignat*}
	as $x\to 0$, and hence
	\begin{align*}
		y_{+,2}(x) &= -\frac{1}{\kappa_2}\biggl(\frac{1}{z}x^{-\rho_2+1}\rho_2x^{\rho_2-1}
		+ \BigO\bigl(x^{\frac{\rho_1+\rho_2}{2}}\bigr)\biggr)
		\to -\frac{\rho_2}{\kappa_2z}\,,
		\\[1ex]
		y_{-,2}(x) &= -\frac{1}{\kappa_2}
		\biggl(\frac{1}{z}x^{-\rho_2+1}(-\kappa_3z)x^{\frac{\rho_1+\rho_2}{2}-1}
		+ \kappa_3x^{\frac{\rho_1-\rho_2}{2}} + \BigO\bigl(x^{\rho_1}\bigr)\biggr)
		\to 0
	\end{align*}
	as $x\to 0$.  The initial condition $W(0,z)=I$ implies that
	\[
		A_1 = 0, \qquad B_1 = 1, \qquad A_2 = -\frac{\kappa_2z}{\rho_2}\,, \qquad B_2 = 0,
	\]
	which proves \eqref{Z49} when $z\ne0$.
	These representations remain true for $z=0$.
\end{proof}

\begin{remark}\label{Z52}
	In the case of a diagonal Hamiltonian, i.e.\ when $\kappa_3=0$,
	the functions $w_{11}$ and $w_{21}$ can be written in terms of Bessel functions:
	\begin{align*}
		w_{11}(x,z) &= \mf J_{-\frac{\rho_2}{\rho_1+\rho_2}}
		\biggl(\frac{\sqrt{\kappa_1\kappa_2}\,}{\rho_3}zx^{\rho_3}\biggr),
		\\[1ex]
		w_{21}(x,z) &= -\frac{\kappa_2}{\rho_2}z\mf J_{\frac{\rho_2}{\rho_1+\rho_2}}
		\biggl(\frac{\sqrt{\kappa_1\kappa_2}\,}{\rho_3}zx^{\rho_3}\biggr),
	\end{align*}
	where
	\begin{align*}
		\mf J_\nu(x) \DEalign \Gamma(\nu+1)\Bigl(\frac{x}{2}\Bigr)^{-\nu}J_\nu(x)
		\\
		&= \sum_{n=0}^\infty \frac{(-1)^n}{n!(\nu+1)_n}\cdot\Bigl(\frac{x}{2}\Bigr)^{2n}
		= \Fhyperg{0}{1}\Bigl(\,;\nu+1;-\frac{x^2}{4}\Bigr)
	\end{align*}
	is the entire function associated with the Bessel function $J_\nu$, which is used in,
	e.g.\ \cite{askey:1973};
	cf.\ also \cite[Example~2]{eckhardt.kostenko.teschl:2018} for the case
	when $\kappa_1=\rho_1>1$ and $\kappa_2=\rho_2=1$.
\end{remark}

\begin{proof}[Proof of \cref{Z132}; direct part]
	We show that the Weyl coefficient of $H_{\uprho,\upkappa}$ is given
	by the asserted formula.
	Again we shorten notation and write $H\equiv H_{\uprho,\upkappa}$.

	Assume first that $\rho_1\ne\rho_2$ and $\kappa>0$.
	Then we can apply \cref{Z44}, which yields
	\[
		\frac{w_{11}(x,z)}{w_{21}(x,z)}
		= -\frac{\rho_2}{\kappa_2z}x^{-\rho_2}
		\frac{M\bigr(a_-,b_-,-\frac{2\kappa iz}{\rho_3}x^{\rho_3}\bigr)}{M\bigl(a_+,b_+,-\frac{2\kappa iz}{\rho_3}x^{\rho_3}\bigr)}.
	\]
	Let $z\in\bb C^+$.  Then $\Re(-iz)>0$, and hence \cite[\nopp 13.2.23]{nist:2010}
	implies that, as $x\to\infty$,
	\begin{align*}
		\frac{w_{11}(x,z)}{w_{21}(x,z)}
		&\sim -\frac{\rho_2}{\kappa_2z}x^{-\rho_2}
		\frac{\frac{\Gamma(b_-)}{\Gamma(a_-)}\exp\bigl(-\frac{2\kappa iz}{\rho_3}x^{\rho_3}\bigr)
		\bigl(-\frac{2\kappa iz}{\rho_3}x^{\rho_3}\bigr)^{a_--b_-}}{\frac{\Gamma(b_+)}{\Gamma(a_+)}
		\exp\bigl(-\frac{2\kappa iz}{\rho_3}x^{\rho_3}\bigr)
		\bigl(-\frac{2\kappa iz}{\rho_3}x^{\rho_3}\bigr)^{a_+-b_+}}
		\\[1ex]
		&= -\frac{\rho_2}{\kappa_2z}x^{-\rho_2}
		\frac{\Gamma(a_+)\Gamma(b_-)}{\Gamma(a_-)\Gamma(b_+)}
		\Bigl(-\frac{2\kappa iz}{\rho_3}x^{\rho_3}\Bigr)^{-a_-+a_-+b_+-b_-}
		\\[1ex]
		&= -\frac{\rho_2}{\kappa_2z}x^{-\rho_2}
		\frac{\Gamma(a_+)\Gamma(b_-)}{\Gamma(a_-)\Gamma(b_+)}
		\Bigl(-\frac{2\kappa iz}{\rho_3}x^{\rho_3}\Bigr)^{\frac{\rho_2}{\rho_3}}
		\displaybreak[0]\\[1ex]
		&= -\frac{\rho_2}{\kappa_2z}x^{-\rho_2}
		\frac{\Gamma(a_+)\Gamma(b_-)}{\Gamma(a_-)\Gamma(b_+)}
		\Bigl(\frac{2\kappa}{\rho_3}\Bigr)^{\frac{\rho_2}{\rho_3}}
		\Bigl(\frac{z}{i}\Bigr)^{\frac{\rho_2}{\rho_3}}x^{\rho_2}
		\\[1ex]
		&= i\frac{\rho_2}{\kappa_2}\cdot\frac{\Gamma(a_+)\Gamma(b_-)}{\Gamma(a_-)\Gamma(b_+)}
		\Bigl(\frac{2\kappa}{\rho_3}\Bigr)^{\frac{\rho_2}{\rho_3}}
		\Bigl(\frac{z}{i}\Bigr)^{\frac{\rho_2}{\rho_3}-1}.
	\end{align*}
	With $\alpha=\frac{\rho_2}{\rho_3}-1=\frac{\rho_2-\rho_1}{\rho_1+\rho_2}$
	and using \eqref{Z47} we obtain
	\begin{align}
		q_H(z) &= \lim_{x\to\infty}\frac{w_{11}(x,z)}{w_{21}(x,z)}
		\nonumber\\[1ex]
		&= i\frac{\rho_2}{\kappa_2}\Bigl(\frac{2\kappa}{\rho_3}\Bigr)^{\alpha+1}\,
		\frac{\Gamma\bigl(\frac12\bigl[\alpha+2+i\frac{\kappa_3}{\kappa}\alpha\bigr]\bigr)
		\Gamma(-\alpha)}{\Gamma\bigl(\frac12\bigl[-\alpha+i\frac{\kappa_3}{\kappa}\alpha\bigr]\bigr)
		\frac{\rho_2}{\rho_3}\Gamma\bigl(\frac{\rho_2}{\rho_3}\bigr)}
		\Bigl(\frac{z}{i}\Bigr)^\alpha
		\nonumber\\[1ex]
		&= i\frac{(2\kappa)^{\alpha+1}}{\kappa_2\rho_3^\alpha}
		\cdot\frac{\Gamma\bigl(1+\frac{\alpha}{2}\bigl(1+i\frac{\kappa_3}{\kappa}\bigr)\bigr)}{%
		\Gamma\bigl(-\frac{\alpha}{2}\bigl(1-i\frac{\kappa_3}{\kappa}\bigr)\bigr)}
		\cdot\frac{\Gamma(-\alpha)}{\Gamma(1+\alpha)}
		\Bigl(\frac{z}{i}\Bigr)^\alpha,
		\label{Z53}
	\end{align}
	which is the asserted formula for the Weyl coefficient in the case when $\rho_1\neq\rho_2$ and $\kappa>0$.

	Using the reflection formula for the gamma function we can rewrite $\omega$,
	\begin{align}
		\omega &= \frac{(2\kappa)^{\alpha+1}}{\kappa_2\rho_3^\alpha}
		\cdot\frac{\Gamma\bigl(1+\frac{\alpha}{2}\bigl(1+i\frac{\kappa_3}{\kappa}\bigr)\bigr)}{%
		\Gamma\bigl(-\frac{\alpha}{2}\bigl(1-i\frac{\kappa_3}{\kappa}\bigr)\bigr)}
		\cdot\frac{\Gamma(-\alpha)}{\Gamma(1+\alpha)}
		\nonumber\\[1ex]
		&= \frac{(2\kappa)^{\alpha+1}}{\kappa_2\rho_3^\alpha}
		\cdot\frac{\Gamma\bigl(1+\frac{\alpha}{2}\bigl(1+i\frac{\kappa_3}{\kappa}\bigr)\bigr)
		\Gamma\bigl(1+\frac{\alpha}{2}\bigl(1-i\frac{\kappa_3}{\kappa}\bigr)\bigr)
		\sin\bigl(-\pi\frac{\alpha}{2}\bigl(1-i\frac{\kappa_3}{\kappa}\bigr)\bigr)}{\pi}\cdot
		\nonumber\\
		&\quad \cdot\frac{\pi}{\bigl(\Gamma(1+\alpha)\bigr)^2\sin(-\pi\alpha)}
		\nonumber\\[1ex]
		&= \frac{(2\kappa)^{\alpha+1}}{\kappa_2\rho_3^\alpha}
		\cdot\frac{\big|\Gamma\bigl(1+\frac{\alpha}{2}\bigl(1+i\frac{\kappa_3}{\kappa}\bigr)\bigr)\big|^2}{%
		\bigl(\Gamma(1+\alpha)\bigr)^2}
		\cdot\frac{\sin\bigl(\frac{\pi\alpha}{2}\bigl(1-i\frac{\kappa_3}{\kappa}\bigr)\bigr)}{\sin(\pi\alpha)}\,.
		\label{Z21}
	\end{align}
	Next let us consider the case when $\rho_1=\rho_2=\rho_3$ and $\kappa>0$.
	Let $\hat\rho_1\DE\rho_3-\varepsilon$, $\hat\rho_2\DE\rho_3+\varepsilon$, $\hat\rho_3=\rho_3$
	with $\varepsilon\in(0,\rho_3)$ and define $H_\varepsilon\DE H_{\hat\uprho,\upkappa}$.
	Then $\hat\alpha=(\hat\rho_2-\hat\rho_1)/(\hat\rho_1+\hat\rho_2)\to0$ as $\varepsilon\to0$.
	Moreover, $\lim_{\varepsilon\to0}H_\varepsilon=H$, and
	hence the Weyl coefficients converge locally uniformly.
	From \eqref{Z21} we obtain
	\begin{equation}\label{Z54}
		q_H(z) = \lim_{\varepsilon\to0}q_{H_\varepsilon}(z)
		= i\frac{2\kappa}{\kappa_2}\cdot\frac{|\Gamma(1)|^2}{(\Gamma(1))^2}
		\cdot\frac{1}{2}\Bigl(1-i\frac{\kappa_3}{\kappa}\Bigr)\Bigl(\frac{z}{i}\Bigr)^0
		= i\frac{\kappa-i\kappa_3}{\kappa_2}\,.
	\end{equation}
	Consider now the case when $\rho_1\ne\rho_2$ and $\kappa=0$, i.e.\
	$|\kappa_3|=\sqrt{\kappa_1\kappa_2}$.
	Let $\hat\kappa_3\in(-\sqrt{\kappa_1\kappa_2},\sqrt{\kappa_1\kappa_2})$
	and let $H_{\hat\kappa_3}$ be the power Hamiltonian with data $\uprho,\kappa_1,\kappa_2$ and $\hat\kappa_3$.
	Further, set $\hat\kappa\DE\sqrt{\kappa_1\kappa_2-\hat\kappa_3^2}$.
	Since $\lim_{\hat\kappa_3\to\kappa_3}H_{\hat\kappa_3}=H$,
	the relation $\lim_{\hat\kappa_3\to\kappa_3}q_{H_{\hat\kappa_3}}(z)=q_H(z)$
	holds locally uniformly.
	According to \cite[\nopp 5.11.12]{nist:2010}
	we have $\frac{\Gamma(a+w)}{\Gamma(b+w)}\sim w^{a-b}$ as $w\to\infty$
	in a sector $|\arg w|\le \phi<\pi$ for arbitrary $a,b\in\bb C$.
	Using this and \eqref{Z53} we can deduce that, as $\hat\kappa_3\to\kappa_3$,
	\begin{align*}
		q_H(z) &= \lim_{\hat\kappa_3\to\kappa_3}q_{H_{\hat\kappa_3}}(z)
		\\
		&= i\lim_{\hat\kappa_3\to\kappa_3}
		\frac{(2\hat\kappa)^{\alpha+1}}{\kappa_2\rho_3^\alpha}
		\cdot\frac{\Gamma\bigl(1+\frac{\alpha}{2}+\frac{i\alpha\hat\kappa_3}{2\hat\kappa}\bigr)}{%
		\Gamma\bigl(-\frac{\alpha}{2}+\frac{i\alpha\hat\kappa_3}{2\hat\kappa}\bigr)}
		\cdot\frac{\Gamma(-\alpha)}{\Gamma(1+\alpha)}\Bigl(\frac{z}{i}\Bigr)^\alpha
		\\[1ex]
		&= i\lim_{\hat\kappa_3\to\kappa_3}
		\frac{(2\hat\kappa)^{\alpha+1}}{\kappa_2\rho_3^\alpha}
		\Bigl(\frac{i\alpha\hat\kappa_3}{2\hat\kappa}\Bigr)^{1+\frac{\alpha}{2}-(-\frac{\alpha}{2})}
		\frac{\Gamma(-\alpha)}{\Gamma(1+\alpha)}\Bigl(\frac{z}{i}\Bigr)^\alpha
		\\[1ex]
		&= i\frac{(i\alpha\kappa_3)^{1+\alpha}}{\kappa_2\rho_3^\alpha}
		\cdot\frac{\Gamma(-\alpha)}{\Gamma(1+\alpha)}
		\Bigl(\frac{z}{i}\Bigr)^\alpha\,,
	\end{align*}
	which proves the formula for $\omega$ when $\rho_1\ne\rho_2$ and $\kappa=0$.
	In a similar way we can establish the case when $\rho_1=\rho_2$ and $\kappa=0$
	by taking the limit as $\hat\kappa_3\to\kappa_3$ in the formula \eqref{Z54}
	with $\kappa_3$ and $\kappa$ replaced by $\hat\kappa_3$ and $\hat\kappa$ respectively.
	It is easy to check that $\omega\ne0$.
\end{proof}

\noindent
For later reference, let us formulate a key observation about the argument of $\omega$
from \eqref{Z134} as a separate lemma.

\begin{lemma}\label{Z88}
	The argument of the number $\omega$ in \eqref{Z134} satisfies
	\[
		\arg\omega
		=
		\begin{cases}
			\displaystyle
			-\arctan\biggl[\tan\Bigl(\frac{\pi}{2}\bigl(1-|\alpha|\bigr)\Bigr)
			\tanh\Bigl(\frac{\pi|\alpha|\kappa_3}{2\kappa}\Bigr)\biggr],
			& \rho_1\ne\rho_2, \ \kappa>0,
			\\[3ex]
			\displaystyle
			-\arctan\Bigl(\frac{\kappa_3}{\kappa}\Bigr),
			& \rho_1=\rho_2, \ \kappa>0,
			\\[3ex]
			\displaystyle
			-(\sgn\kappa_3)\frac{\pi}{2}\bigl(1-|\alpha|\bigr),
			& \rho_1=\rho_2, \ \kappa=0.
		\end{cases}
	\]
	For each fixed $\rho_1,\rho_2,\kappa_1,\kappa_2$, the function $\kappa_3\mapsto\arg\omega$
	is a strictly decreasing and odd bijection from $[-\sqrt{\kappa_1\kappa_2},\sqrt{\kappa_1\kappa_2}]$
	onto $\bigl[-\frac{\pi}{2}(1-|\alpha|),\frac{\pi}{2}(1-|\alpha|)\bigr]$.
\end{lemma}

\begin{proof}
	Consider first the case when $\rho_1\ne\rho_2$ and $\kappa>0$.
	From \eqref{Z21} we obtain
	\begin{align*}
		\arg\omega
		&= \arg\biggl(\frac{\sin\bigl(\frac{\pi\alpha}{2}\bigl(1-i\frac{\kappa_3}{\kappa}\bigr)\bigr)}{\sin(\pi\alpha)}\biggr)
		\\[1ex]
		&= \arg\biggl(\frac{\sin\bigl(\frac{\pi\alpha}{2}\bigr)\cosh\bigl(\frac{\pi\alpha\kappa_3}{2\kappa}\bigr)
		-i\cos\bigl(\frac{\pi\alpha}{2}\bigr)\sinh\bigl(\frac{\pi\alpha\kappa_3}{2\kappa}\bigr)}{%
		2\sin\bigl(\frac{\pi\alpha}{2}\bigr)\cos\bigl(\frac{\pi\alpha}{2}\bigr)}\biggr)
		\displaybreak[0]\\[1ex]
		&= \arg\biggl(\frac{\cosh\bigl(\frac{\pi|\alpha|\kappa_3}{2\kappa}\bigr)}{\cos\bigl(\frac{\pi|\alpha|}{2}\bigr)}
		-i\frac{\sinh\bigl(\frac{\pi|\alpha|\kappa_3}{2\kappa}\bigr)}{\sin\bigl(\frac{\pi|\alpha|}{2}\bigr)}
		\biggr)
		= -\arctan\biggl[\cot\Bigl(\frac{\pi|\alpha|}{2}\Bigr)
		\tanh\Bigl(\frac{\pi|\alpha|\kappa_3}{2\kappa}\Bigr)\biggr]
		\\[1ex]
		&= -\arctan\biggl[\tan\Bigl(\frac{\pi}{2}\bigl(1-|\alpha|\bigr)\Bigr)
		\tanh\Bigl(\frac{\pi|\alpha|\kappa_3}{2\kappa}\Bigr)\biggr].
	\end{align*}
	The other cases follow easily by taking limits.

	The last statement of the lemma follows by inspecting the given explicit formulae.
\end{proof}

\noindent
To show the inverse part of \cref{Z132}, i.e.\ that the Weyl coefficient map
is surjective, we use the following observation.

\begin{lemma}\label{Z87}
	Let $\alpha\in(-1,1)$, $\kappa_1,\kappa_2>0$, $\sigma>0$
	and $\hat\omega\in\bb C\setminus\{0\}$
	such that $|\arg\hat\omega|\le\frac{\pi}{2}(1-|\alpha|)$.  Set
	\begin{equation}\label{Z37}
	\begin{alignedat}{3}
		\rho_1 &= \sigma, \quad & \rho_2 &= \frac{1+\alpha}{1-\alpha}\sigma
		\qquad & &\text{if} \ \alpha \ge 0,
		\\
		\rho_1 &= \frac{1-\alpha}{1+\alpha}\sigma, \quad & \rho_2 &= \sigma
		\qquad & &\text{if} \ \alpha < 0,
	\end{alignedat}
	\end{equation}
	and $\rho_3=\frac12(\rho_1+\rho_2)$.
	Then there exist $\gamma>0$ and $\kappa_3\in\bb R$ such
	that $\kappa_3^2\le\kappa_1\kappa_2$, $\min\{\rho_1,\rho_2\}=\sigma$,
	and the Weyl coefficient of $H_{\uprho,\upkappa}$ with $\uprho=(\rho_1,\rho_2,\rho_3)$
	and $\upkappa=(\kappa_1,\kappa_2,\kappa_3)$ satisfies
	\[
		q_{H_{\uprho,\upkappa}}(z)
		= \frac{1}{\gamma}i\hat\omega\Bigl(\frac{z}{i}\Bigr)^\alpha,
		\qquad z\in\bb C^+.
	\]
\end{lemma}

\begin{proof}
	It is easy to check that \eqref{Z102} and $\min\{\rho_1,\rho_2\}=\sigma$ hold.
	By \cref{Z88} there exists a unique
	$\kappa_3\in[-\sqrt{\kappa_1\kappa_2},\sqrt{\kappa_1\kappa_2}]$ such that
	$\arg\omega=\arg\hat\omega$ where $\omega$ is from \eqref{Z134}.
	Now the representation for $q_{H_{\uprho,\upkappa}}$ follows from
	the direct part of \cref{Z132}.
\end{proof}

\begin{proof}[Proof of \cref{Z132}; inverse part]
	This is now easily obtained by using rescalings.

	Choose $\kappa_1=\kappa_2=\sigma=1$ and $\hat\omega=\omega$ in \cref{Z87},
	which yields $\rho_1,\rho_2\ge1$, $\tilde\kappa_3\in[-1,1]$ and $\gamma>0$ such that
	\[
		q_{H_{\uprho,\tilde\upkappa}}(z) = \frac{1}{\gamma}i\omega\Bigl(\frac{z}{i}\Bigr)^\alpha,
		\qquad z\in\bb C^+,
	\]
	where $\tilde\upkappa=(1,1,\tilde\kappa_3)$.
	We apply the transformation $\mc A_r^{b_1,b_2}$ from \cref{Z2} with $b_1=\gamma$
	and $b_2=r=1$:
	\[
		\bigl(\mc A_1^{\gamma,1}H_{\uprho,\tilde\upkappa}\bigr)(t)
		=
		\begin{pmatrix}
			\gamma^{\rho_1+1}t^{\rho_1-1} & \gamma^{\rho_3}\tilde\kappa_3t^{\rho_3-1}
			\\[1ex]
			\gamma^{\rho_3}\tilde\kappa_3t^{\rho_3-1} & \gamma^{\rho_2-1}t^{\rho_2-1}
		\end{pmatrix}
		,
	\]
	which is again a power Hamiltonian.  The asserted form of the Weyl coefficient
	follows from \eqref{Z39}.
\end{proof}

\noindent
We can also determine explicitly which power Hamiltonians have the
same Weyl coefficient.

\begin{proposition}\label{Z106}
	Let $\rho_1,\rho_2,\kappa_1,\kappa_2,\tilde\rho_1,\tilde\rho_2,\tilde\kappa_1,
	\tilde\kappa_2>0$ and $\kappa_3,\tilde\kappa_3\in\bb R$ with $\kappa_3^2\le\kappa_1\kappa_2$
	and $\tilde\kappa_3^2\le\tilde\kappa_1\tilde\kappa_2$
	and set $\rho_3=\frac12(\rho_1+\rho_2)$ and $\tilde\rho_3=\frac12(\tilde\rho_1+\tilde\rho_2)$
	and $\uprho=(\rho_1,\rho_2,\rho_3)$, $\upkappa=(\kappa_1,\kappa_2,\kappa_3)$,
	$\tilde\uprho=(\tilde\rho_1,\tilde\rho_2,\tilde\rho_3)$
	and $\tilde\upkappa=(\tilde\kappa_1,\tilde\kappa_2,\tilde\kappa_3)$.
	Then $q_{H_{\uprho,\upkappa}}=q_{H_{\tilde\uprho,\tilde\upkappa}}$
	if and only if there exist $\beta,c>0$ such that
	\begin{equation}\label{Z108}
	\begin{alignedat}{2}
		\tilde\rho_i &= \beta\rho_i, \qquad & & i=1,2,
		\\[1ex]
		\tilde\kappa_i &= \beta c^{\rho_i}\kappa_i, \qquad & & i=1,2,3.
	\end{alignedat}
	\end{equation}
\end{proposition}

\begin{proof}
	Assume that $q_{H_{\uprho,\upkappa}}=q_{H_{\tilde\uprho,\tilde\upkappa}}$.
	Then $H_{\uprho,\upkappa}$ and $H_{\tilde\uprho,\tilde\upkappa}$
	are reparameterisations of each other.
	Hence there exists a strictly increasing bijection $\gamma:[0,\infty)\to[0,\infty)$
	such that $\gamma$ and $\gamma^{-1}$ are locally absolutely continuous
	and $\tildeM(t)=M(\gamma(t))$, $t\in[0,\infty)$; see \eqref{Z36}.
	Explicitly, we have
	\[
		\frac{\tilde\kappa_j}{\tilde\rho_j}t^{\tilde\rho_j}
		= \frac{\kappa_j}{\rho_j}\bigl(\gamma(t)\bigr)^{\rho_j},
		\qquad j=1,2,3; \, t\in[0,\infty),
	\]
	and therefore
	\[
		\gamma(t) = \Bigl(\frac{\tilde\kappa_j\rho_j}{\kappa_j\tilde\rho_j}\Bigr)^{\frac{1}{\rho_j}}
		t^{\frac{\tilde\rho_j}{\rho_j}}
		\ED ct^\beta,
		\qquad j=1,2,3; \, t\in[0,\infty),
	\]
	which implies \eqref{Z108}.
	The converse is clear.
\end{proof}

\noindent
Note that one can choose a unique representative modulo the kernel of the
Weyl coefficient map, e.g.\ by fixing $\rho_1=\kappa_1=1$, i.e.\ $h_1(t)\equiv1$.

\section{Proof of the main theorems}
\label{Z63}

In Sections~\ref{Z244}--\ref{Z179} we prove \cref{Z10}, in \cref{Z243} we prove \cref{Z236},
and the proof of \cref{Z124} is done in \cref{Z92}.
Throughout Sections~\ref{Z244}--\ref{Z179} we fix a Hamiltonian $H$ such that neither $h_1$ nor $h_2$
vanishes on a neighbourhood of $0$.  In Sections~\ref{Z244} and \ref{Z245} we prove the implications
(i)$\Rightarrow$(ii) and (ii)$\Rightarrow$(i) respectively and further statements of \cref{Z10}
for trace-normalised $H$ on $(0,\infty)$.  The general case is then treated in \cref{Z179}.

\subsection[{Proof of (i)$\Rightarrow$(ii) in \cref{Z10} ($\tr H=1$ a.e.)}]{Proof of (i)$\bm\Rightarrow$(ii) in \cref{Z10} ($\bm{\tr H=1$} a.e.)}
\label{Z244}

Assume that we have a regularly varying function $\ms a$ and a non-zero constant $\omega$
such that
\[
	q_H(ri) \sim i\omega\ms a(r)
	\qquad\text{as} \ r\to\infty.
\]
By \cite[Theorem~1.3.3]{bingham.goldie.teugels:1989} we may assume w.l.o.g.\ that $\ms a$ is continuous.

Denote the index of $\ms a$ by $\alpha$.  By the assumption that neither $h_1$ nor $h_2$
vanishes on any neighbourhood of $0$, we have
\begin{equation}\label{Z51}
	\frac{\ms a(r)}{r} \to 0 \quad\text{and}\quad r\ms a(r) \to\infty,
	\qquad r\to\infty;
\end{equation}
see \cref{Z90}\,(ii).  This, together with \eqref{Z27}, implies that $\alpha\in[-1,1]$.

\begin{Steps}
\item 
\textit{We show convergence of rescalings of $H$ to a comparison Hamiltonian.}

	\noindent
	It follows from \cite[Theorem~3.1]{langer.woracek:kara} that $|\arg\omega|\le\frac\pi 2(1-|\alpha|)$ and
	\begin{equation}\label{Z115}
		q_H(rz) \sim i\omega\Bigl(\frac zi\Bigr)^\alpha \ms a(r),
		\qquad r\to\infty,
	\end{equation}
	locally uniformly for $z\in\bb C^+$.

	An appropriate comparison Hamiltonian is obtained from \cref{Z87,Z133}.
	If $\alpha\in(-1,1)$, then there exist
	\begin{Itemize}
	\item
		$\rho_1,\rho_2\geq 1$ with $\min\{\rho_1,\rho_2\}=1$
		and $\alpha=\frac{\rho_2-\rho_1}{\rho_2+\rho_1}$
		(namely, choose $\rho_1,\rho_2$ as in \eqref{Z37} with $\sigma=1$),
	\item
		$\kappa_3\in[-1,1]$ and $\gamma>0$,
	\end{Itemize}
	such that, with $\rho_3\DE\frac 12(\rho_1+\rho_2)$, the Hamiltonian
	\[
		\widetilde H(t)\DE
		\begin{pmatrix}
			t^{\rho_1-1} & \kappa_3t^{\rho_3-1}
			\\[1ex]
			\kappa_3t^{\rho_3-1} & t^{\rho_2-1}
		\end{pmatrix},
		\qquad t\in(0,\infty),
	\]
	(i.e.\ $\widetilde H=H_{\uprho,\upkappa}$ with $\uprho=(\rho_1,\rho_2,\rho_3)$, $\upkappa=(1,1,\kappa_3)$)
	satisfies
	\begin{equation}\label{Z93}
		i\omega\Bigl(\frac zi\Bigr)^\alpha = \gamma q_{\widetilde H}(z),
		\qquad z\in\bb C^+.
	\end{equation}
	If $\alpha=1$, we have $\gamma\DE\omega>0$, and \eqref{Z93} holds for the Hamiltonian
	\[
		\widetilde H(t)\DE
		\begin{pmatrix}
			\mathds{1}_{[0,1]}(t) & 0
			\\[1ex]
			0 & \mathds{1}_{(1,\infty)}(t)
		\end{pmatrix}.
	\]
	If $\alpha=-1$, we again have $\gamma\DE\omega>0$, and \eqref{Z93} holds for
	\[
		\widetilde H(t)\DE
		\begin{pmatrix}
			\mathds{1}_{(1,\infty)}(t) & 0
			\\[1ex]
			0 & \mathds{1}_{[0,1]}(t)
		\end{pmatrix}.
	\]
	Let $\ms b_1,\ms b_2\DF(0,\infty)\to(0,\infty)$ be any two functions
	(a particular choice will be made later) such that
	\begin{equation}\label{Z155}
		\frac{\ms b_2(r)}{\ms b_1(r)} = \gamma\ms a(r),
		\qquad r>0.
	\end{equation}
	Then we obtain from \eqref{Z115}, \eqref{Z155} and \eqref{Z93} that
	\[
		\lim_{r\to\infty}\frac{\ms b_1(r)}{\ms b_2(r)}q_H(rz)
		= \frac{1}{\gamma}\lim_{r\to\infty}\frac{q_H(rz)}{\ms a(r)}
		= q_{\tilde H}(z)
	\]
	locally uniformly for $z\in\bb C^+$, and hence \cref{Z14} yields
	\begin{equation}\label{Z24}
		\lim_{r\to\infty}\mc A_r^{\ms b_1,\ms b_2}H = \widetilde H.
	\end{equation}
	Recall that, with the notation
	\begin{equation}\label{Z75}
		\mf t_r(x)\DE\int_0^x\tr(\mc A_r^{\ms b_1,\ms b_2}H)(s)\RD s,
		\qquad
		\tilde{\mf t}(x)\DE\int_0^x\tr\widetilde H(s)\RD s,
	\end{equation}
	the limit relation \eqref{Z24} means that
	\begin{equation}\label{Z25}
		\lim_{r\to\infty}\int_0^{\mf t_r^{-1}(T)}(\mc A_r^{\ms b_1,\ms b_2}H)(s)\RD s
		= \int_0^{\tilde{\mf t}^{-1}(T)}\widetilde H(s)\RD s
	\end{equation}
	holds locally uniformly for $T\in[0,\infty)$; see \cref{Z74}\,(i).

\item \textit{We relate the functions $\mf t_r^{-1}$ and $\tilde{\mf t}^{-1}$.}

	\noindent
	Let us evaluate integrals over the entries of $\mc A_r^{\ms b_1,\ms b_2}H$.
	For $x\ge0$ we have
	\begin{equation}\label{Z70}
	\begin{aligned}
		& \binom 10^*\int_0^x (\mc A_r^{\ms b_1,\ms b_2}H)(s)\RD s\,\binom 10
		= \int_0^x \ms b_1(r)^2 h_1\Bigl(\frac{\ms b_1(r)\ms b_2(r)}{r}s\Bigr)\RD s
		\\[1ex]
		&= \ms b_1(r)^2 \int_0^{\frac{\ms b_1(r)\ms b_2(r)}{r}x}
		\frac{r}{\ms b_1(r)\ms b_2(r)}h_1(t)\RD t
		= r\frac{\ms b_1(r)}{\ms b_2(r)}m_1\Bigl(\frac{\ms b_1(r)\ms b_2(r)}{r}x\Bigr),
	\end{aligned}
	\end{equation}
	and similarly,
	\begin{align}
		\binom 01^*\int_0^x (\mc A_r^{\ms b_1,\ms b_2}H)(s)\RD s\,\binom 01
		&= r\frac{\ms b_2(r)}{\ms b_1(r)}m_2\Bigl(\frac{\ms b_1(r)\ms b_2(r)}{r}x\Bigr),
		\label{Z71}
		\\[1ex]
		\binom 10^*\int_0^x (\mc A_r^{\ms b_1,\ms b_2}H)(s)\RD s\,\binom 01
		&= r m_3\Bigl(\frac{\ms b_1(r)\ms b_2(r)}{r}x\Bigr).
		\label{Z72}
	\end{align}
	Hence, with the notation
	$\tildeM(x)=\smmatrix{\widetilde m_1(x)}{\widetilde m_3(x)}{\widetilde m_3(x)}{\widetilde m_2(x)}
	=\int_0^x \widetilde H(t)\RD t$, the limit relation \eqref{Z25} can be written as
	\begin{align}
		r\frac{\ms b_1(r)}{\ms b_2(r)}m_1\Bigl(\frac{\ms b_1(r)\ms b_2(r)}r\cdot\mf t_r^{-1}(T)\Bigr)
		&= \widetilde m_1\bigl(\tilde{\mf t}^{-1}(T)\bigr)+R_1(T,r),
		\label{Z29}
		\\[1ex]
		r\frac{\ms b_2(r)}{\ms b_1(r)}m_2\Bigl(\frac{\ms b_1(r)\ms b_2(r)}r\cdot\mf t_r^{-1}(T)\Bigr)
		&= \widetilde m_2\bigl(\tilde{\mf t}^{-1}(T)\bigr)+R_2(T,r),
		\label{Z30}
		\\[1ex]
		rm_3\Bigl(\frac{\ms b_1(r)\ms b_2(r)}r\cdot\mf t_r^{-1}(T)\Bigr)
		&= \widetilde m_3\bigl(\tilde{\mf t}^{-1}(T)\bigr)+R_3(T,r),
		\label{Z95}
	\end{align}
	for $T\in[0,\infty)$, $r>0$, where
	\[
		\lim_{r\to\infty}R_1(T,r) =\lim_{r\to\infty}R_2(T,r)=\lim_{r\to\infty}R_3(T,r) = 0
	\]
	locally uniformly for $T\in[0,\infty)$.

	Since $\tr H=1$ a.e., we have $m_1(x)+m_2(x)=x$, $x\ge0$.
	Dividing \eqref{Z29} by $\ms b_1(r)^2$ and \eqref{Z30} by $\ms b_2(r)^2$,
	and taking the sum, we obtain
	\begin{align}
		\mf t_r^{-1}(T) &=
		\frac{1}{\ms b_1(r)^2}\widetilde m_1\bigl(\tilde{\mf t}^{-1}(T)\bigr)
		+\frac{1}{\ms b_2(r)^2}\widetilde m_2\bigl(\tilde{\mf t}^{-1}(T)\bigr)
		\nonumber
		\\[1ex]
		&\quad +\Bigl(\frac 1{\ms b_1(r)^2}+\frac 1{\ms b_2(r)^2}\Bigr)
		\biggl[\underbrace{\frac{\ms b_2(r)^2}{\ms b_1(r)^2+\ms b_2(r)^2}R_1(T,r)+
		\frac{\ms b_1(r)^2}{\ms b_1(r)^2+\ms b_2(r)^2}R_2(T,r)}_{\ED\tilde R(T,r)}\biggr]
	\label{Z26}
	\end{align}
	for all $T\in[0,\infty)$ and $r>0$.  Clearly, $\lim_{r\to\infty}\tilde R(T,r)=0$
	locally uniformly for $T\in[0,\infty)$.

\item \textit{We specify a choice of $\ms b_1$ and $\ms b_2$.}

	\noindent
	It turns out to be appropriate to let $\ms b_1$ and $\ms b_2$ balance the contributions
	of upper and lower entries in the following way: set
	\[
		\ms b_2(r) \DE \sqrt{1+(\gamma\ms a(r))^2},
		\qquad
		\ms b_1(r) \DE \frac{\ms b_2(r)}{\gamma\ms a(r)}.
	\]
	A couple of properties of these functions are obvious:
	\begin{Itemize}
	\item 
		\eqref{Z155} holds;
	\item 
		$\ms b_2(r)\ge\max\{1,\gamma\ms a(r)\}$ and $\ms b_1(r)\ge\max\{1,(\gamma\ms a(r))^{-1}\}$;
	\item 
		$\frac 1{\ms b_1(r)^2}+\frac 1{\ms b_2(r)^2}=1$;
	\item 
		$\ms b_1(r)\ms b_2(r)=[\gamma\ms a(r)]+[\gamma\ms a(r)]^{-1}$.
	\end{Itemize}
	For later reference, observe that $\ms b_1\ms b_2$ is continuous since $\ms a$ is, and that
	\[
		\lim_{r\to\infty}\frac{\ms b_1(r)\ms b_2(r)}r=0
	\]
	by \eqref{Z51}.
	From now on this choice of $\ms b_1$ and $\ms b_2$ is kept fixed.

\item \textit{We show that the limit $\delta$ in \eqref{Z56} involving the
	off-diagonal entries exists.}

	\noindent
	Assume that $\alpha\in(-1,1)$, fix $x_1>1$ and set $T\DE\tilde{\mf t}(x_1)$. The function $r\mapsto\mf t_r^{-1}(T)$ is continuous by \cref{Z19}
	and bounded by \eqref{Z26}. Moreover, $\widetilde m_1(x_1)\widetilde m_2(x_1)>0$.
	Using \eqref{Z29}--\eqref{Z95} we conclude that
	\begin{align*}
		& \lim_{t\to 0}\frac{m_3(t)}{\sqrt{m_1(t)m_2(t)}\,}
		= \lim_{r\to\infty}\frac{m_3\bigl(\frac{\ms b_1(r)\ms b_2(r)}{r}\mf t_r^{-1}(T)\bigr)}{%
		\sqrt{m_1\bigl(\frac{\ms b_1(r)\ms b_2(r)}{r}\mf t_r^{-1}(T)\bigr)
		m_2\bigl(\frac{\ms b_1(r)\ms b_2(r)}{r}\mf t_r^{-1}(T)\bigr)}\,}
		\\[1ex]
		&= \lim_{r\to\infty}\frac{rm_3\bigl(\frac{\ms b_1(r)\ms b_2(r)}{r}\mf t_r^{-1}(T)\bigr)}{%
		\sqrt{r\frac{\ms b_1(r)}{\ms b_2(r)}m_1\bigl(\frac{\ms b_1(r)\ms b_2(r)}{r}\mf t_r^{-1}(T)\bigr)\cdot
		r\frac{\ms b_2(r)}{\ms b_1(r)}m_2\bigl(\frac{\ms b_1(r)\ms b_2(r)}{r}\mf t_r^{-1}(T)\bigr)}\,}
		\\[1ex]
		&= \frac{\widetilde m_3(x_1)}{\sqrt{\widetilde m_1(x_1)\widetilde m_2(x_1)}\,}.
	\end{align*}

\item \textit{We show that rescalings converge on fixed intervals.}

	\noindent
	First, we establish
	\begin{equation}\label{Z69}
		\lim_{r\to\infty}\mf t_r^{-1}(T) = \tilde{\mf t}^{-1}(T),\qquad T\in[0,\tilde{\mf t}(1)],
	\end{equation}
	by distinguishing cases for $\alpha$.
	If $\alpha>0$, then $\lim_{r\to\infty}\ms a(r)=\infty$ and $\widetilde m_1(x)=x$
	for all $x\in[0,1]$.
	Thus $\lim_{r\to\infty}\ms b_2(r)=\infty$, which implies that $\lim_{r\to\infty}\ms b_1(r)=1$.
	Together with \eqref{Z26}, this shows that
	\begin{equation}\label{Z55}
		\lim_{r\to\infty}\mf t_r^{-1}(T) = \widetilde m_1\bigl(\tilde{\mf t}^{-1}(T)\bigr),
		\qquad T\in[0,\infty),
	\end{equation}
	which, in turn, establishes \eqref{Z69}.
	If $\alpha<0$, then $\lim_{r\to\infty}\ms a(r)=0$ and $\widetilde m_2(x)=x$
	for all $x\in[0,1]$.
	Thus $\lim_{r\to\infty}\ms b_1(r)=\infty$ and $\lim_{r\to\infty}\ms b_2(r)=1$.
	Also in this case \eqref{Z69} follows from \eqref{Z26}.
	If $\alpha=0$, then $\widetilde m_1(x)=\widetilde m_2(x)=x$ for all $x\ge 0$,
	and again \eqref{Z26} implies \eqref{Z69}.

	Since $\tilde{\mf t}$ is an increasing bijection from $[0,1]$ onto $[0,\tilde{\mf t}(1)]$,
	and all functions $\mf t_r$ are non-decreasing, the limit relation \eqref{Z69} implies that
	\[
		\lim_{r\to\infty}\mf t_r(x)=\tilde{\mf t}(x), \qquad x\in(0,1).
	\]
	As \eqref{Z25} holds locally uniformly in $T$, we obtain that
	\begin{equation}\label{Z28}
		\lim_{r\to\infty}\int_0^x(\mc A_r^{\ms b_1,\ms b_2}H)(s)\RD s
		= \int_0^x\widetilde H(s)\RD s,
		\qquad x\in(0,1).
	\end{equation}

\item \textit{We prove regular variation of $m_1$ and $m_2$.}
	
	\noindent
	Assume first that $\alpha\in(-1,1]$.  Then $\widetilde m_1(x)=\frac 1{\rho_1}x^{\rho_1}$
	for all $x\in[0,1]$ (here $\rho_1\DE 1$ if $\alpha=1$).
	We can use \eqref{Z70} to rewrite the left upper entry of the left-hand side of \eqref{Z28},
	which then yields
	\[
		\lim_{r\to\infty}r\frac{\ms b_1(r)}{\ms b_2(r)}m_1\Bigl(\frac{\ms b_1(r)\ms b_2(r)}rx\Bigr)
		= \frac{x^{\rho_1}}{\rho_1},
		\qquad x\in(0,1),
	\]
	and therefore
	\[
		\lim_{t\to 0}\frac{m_1(xt)}{m_1(t)}
		= \lim_{r\to\infty}\frac{r\frac{\ms b_1(r)}{\ms b_2(r)}\cdot m_1\bigl(x\cdot\frac{\ms b_1(r)\ms b_2(r)}rx\bigr)}{%
		r\frac{\ms b_1(r)}{\ms b_2(r)}\cdot m_1\bigl(\frac{\ms b_1(r)\ms b_2(r)}rx\bigr)}
		= x^{\rho_1},
		\qquad x\in(0,1).
	\]
	It follows from \cite[Theorem~1.4.1]{bingham.goldie.teugels:1989} that $m_1$ is regularly varying 
	at $0$ with index $\rho_1$.

	Assume now that $\alpha\in[-1,1)$. This case is completely dual. We have
	$\widetilde m_2(x)=\frac 1{\rho_2}x^{\rho_2}$ for all $x\in[0,1]$
	(where $\rho_2\DE 1$ if $\alpha=-1$).
	Based on \eqref{Z71} we obtain
	\[
		\lim_{t\to 0}\frac{m_2(xt)}{m_2(t)}
		= \lim_{r\to\infty}\frac{r\frac{\ms b_2(r)}{\ms b_1(r)}\cdot m_2\bigl(x\cdot\frac{\ms b_1(r)\ms b_2(r)}rx\bigr)}{%
		r\frac{\ms b_2(r)}{\ms b_1(r)}\cdot m_2\bigl(\frac{\ms b_1(r)\ms b_2(r)}rx\bigr)}
		= x^{\rho_2},
		\qquad x\in(0,1),
	\]
	and conclude that $m_2$ is regularly varying at $0$ with index $\rho_2$.

\item \textit{The case of rapid variation of $m_1$ or $m_2$.}

	\noindent
	Assume that $\alpha=1$.  Then
	\[
		\widetilde m_1(t)
		= \begin{cases}
			t, & t\in[0,1],
			\\[0.5ex]
			1, & t\in(1,\infty),
		\end{cases}
		\qquad
		\widetilde m_2(t)
		= \begin{cases}
			0, & t\in[0,1],
			\\[0.5ex]
			t-1, & t\in(1,\infty),
		\end{cases}
	\]
	and $\tilde{\mf t}(t)=t$ for all $t\geq 0$. Fix $x_1>1$ and $x\in(0,1)$.
	It follows from \eqref{Z55} that
	\[
		\lim_{r\to\infty}\mf t_r^{-1}(x_1)=1
		\qquad\text{and}\qquad
		\lim_{r\to\infty}\frac{\mf t_r^{-1}\bigl(\sqrt{x}\bigr)}{x}
		= \frac{\sqrt{x}\,}{x} >1.
	\]
	Hence, for large enough $r$ we have $\mf t_r^{-1}(x_1)\le \frac{1}{x}\mf t_r^{-1}(\sqrt{x})$
	and, by the monotonicity of $m_2$,
	\[
		m_2\Big(\frac{\ms b_1(r)\ms b_2(r)}r\cdot\mf t_r^{-1}(x_1)\Big)
		\le m_2\Big(\frac{\ms b_1(r)\ms b_2(r)}r\cdot\frac{1}{x}\mf t_r^{-1}(\sqrt{x})\Big).
	\]
	We can use \eqref{Z30} to obtain
	\begin{align*}
		\lim_{t\to 0}\frac{m_2(xt)}{m_2(t)}
		&= \lim_{r\to\infty}\frac{m_2\bigl(x\cdot \frac{\ms b_1(r)\ms b_2(r)}r\frac 1x\mf t_r^{-1}(\sqrt{x})\bigr)}{%
		m_2\bigl(\frac{\ms b_1(r)\ms b_2(r)}r\frac{1}{x}\mf t_r^{-1}(\sqrt{x})\bigr)}
		\\[1ex]
		&\le \lim_{r\to\infty}
		\frac{r\frac{\ms b_2(r)}{\ms b_1(r)}\cdot m_2\bigl(\frac{\ms b_1(r)\ms b_2(r)}r\mf t_r^{-1}(\sqrt{x})\bigr)}{%
		r\frac{\ms b_2(r)}{\ms b_1(r)}\cdot m_2\bigl(\frac{\ms b_1(r)\ms b_2(r)}r\mf t_r^{-1}(x_1)\bigr)}
		= \frac{\widetilde m_2(\sqrt{x})}{\widetilde m_2(x_1)}=0,
	\end{align*}
	which shows that $m_2$ is rapidly varying at $0$
	because the case when $x\in(1,\infty)$ can be reduced to the case $x\in(0,1)$
	by considering reciprocals.

	The proof that $m_1$ is rapidly varying when $\alpha=-1$ is completely dual; we skip details.
\end{Steps}

\noindent
All properties stated in \cref{Z10}\,\Enumref{2} are established;
see Steps~\ding{175}, \ding{177}, \ding{178}.

\subsection[Proof of (ii)$\Rightarrow$(i) in \cref{Z10} ($\tr H=1$ a.e.)]{Proof of (ii)$\bm\Rightarrow$(i) in \cref{Z10} ($\bm{\tr H=1}$ a.e.)}
\label{Z245}

Assume that $m_1$ and $m_2$ are regularly or rapidly varying at $0$
with indices $\rho_1$ and $\rho_2$ respectively, and that
\[
	\delta\DE\lim_{t\to0}\frac{m_3(t)}{\sqrt{m_1(t)m_2(t)}\,}
\]
exists if $\max\{\rho_1,\rho_2\}<\infty$.
Further, set $\rho_3\DE\frac 12(\rho_1+\rho_2)$ if $\max\{\rho_1,\rho_2\}<\infty$.

It follows from \eqref{Z18} that
\[
	\lim_{t\to 0}\frac{m_j(t)}{t}=
	\begin{cases}
		0, & \rho_j>1,
		\\[0.5ex]
		\infty, & \rho_j<1.
	\end{cases}
\]
The fact that $\tr H=1$ a.e.\ implies that $m_1(t)+m_2(t)=t$, and thus
\begin{equation}\label{Z101}
	\min\{\rho_1,\rho_2\}=1.
\end{equation}
We define a matrix-valued function $\widetilde H\DF(0,\infty)\to\bb R^{2\times 2}$
such that $\widetilde H(t)$ equals
\[
	\begin{pmatrix}
		\rho_1 t^{\rho_1-1} & \delta\rho_3 t^{\rho_3-1}
		\\[1ex]
		\delta\rho_3 t^{\rho_3-1} & \rho_2 t^{\rho_2-1}
	\end{pmatrix}
	\quad\text{or}\quad
	\begin{pmatrix}
		\mathds{1}_{[0,1]}(t) & 0
		\\[1ex]
		0 & \mathds{1}_{(1,\infty)}(t)
	\end{pmatrix}
	\quad\text{or}\quad
	\begin{pmatrix}
		\mathds{1}_{(1,\infty)}(t) & 0
		\\[1ex]
		0 & \mathds{1}_{[0,1]}(t)
	\end{pmatrix}
\]
a.e.\ according to the cases `$\max\{\rho_1,\rho_2\}<\infty$' or `$\rho_2=\infty$'
or `$\rho_1=\infty$', respectively. Further, we set
\[
	\ms b_1(r) \DE r\sqrt{\mr t(r)m_2\bigl(\mr t(r)\bigr)},
	\qquad
	\ms b_2(r)\DE r\sqrt{\mr t(r)m_1\bigl(\mr t(r)\bigr)}.
\]
The essence of the proof is to show that corresponding rescalings of $H$ converge
to $\widetilde H$.  Recall that $\mr t(r)$ is the unique number
with $(m_1m_2)\bigl(\mr t(r)\bigr) = \frac{1}{r^2}$; see \eqref{Z38}.

We start with the generic case when $m_1$ and $m_2$ are both regularly varying.
Relations \eqref{Z70}--\eqref{Z72} applied with the present functions
$\ms b_1,\ms b_2$ show that
\begin{align}
	\binom10^*\int_0^x \bigl(\mc A_r^{\ms b_1\ms b_2}H\bigr)(s)\RD s\,\binom10
	&= r\sqrt{\frac{m_2\bigl(\mr t(r)\bigr)}{m_1\bigl(\mr t(r)\bigr)}\,}
	m_1\bigl(\mr t(r)x\bigr)
	= \frac{m_1\bigl(\mr t(r)x\bigr)}{m_1\bigl(\mr t(r)\bigr)},
	\label{Z80}
	\\[1ex]
	\binom01^*\int_0^x \bigl(\mc A_r^{\ms b_1\ms b_2}H\bigr)(s)\RD s\,\binom01
	&= r\sqrt{\frac{m_1\bigl(\mr t(r)\bigr)}{m_2\bigl(\mr t(r)\bigr)}\,}
	m_2\bigl(\mr t(r)x\bigr)
	= \frac{m_2\bigl(\mr t(r)x\bigr)}{m_2\bigl(\mr t(r)\bigr)},
	\label{Z82}
	\\[1ex]
	\binom10^*\int_0^x \bigl(\mc A_r^{\ms b_1\ms b_2}H\bigr)(s)\RD s\,\binom01
	&= r m_3\bigl(\mr t(r)x\bigr)
	= \frac{m_3\bigl(\mr t(r)x\bigr)}{\sqrt{m_1\bigl(\mr t(r)\bigr)m_2\bigl(\mr t(r)\bigr)}\,}
	\nonumber
	\\[1ex]
	&\hspace*{-10ex} 
	= \frac{m_3\bigl(\mr t(r)x\bigr)}{\sqrt{m_1\bigl(\mr t(r)x\bigr)m_2\bigl(\mr t(r)x\bigr)}\,}
	\cdot\sqrt{\frac{m_1\bigl(\mr t(r)x\bigr)}{m_1\bigl(\mr t(r)\bigr)}\,}
	\cdot\sqrt{\frac{m_2\bigl(\mr t(r)x\bigr)}{m_2\bigl(\mr t(r)\bigr)}\,}.
	\label{Z84}
\end{align}
for all $x>0$.

\begin{lemma}\label{Z96}
	Assume that \textup{(ii)} holds with $\max\{\rho_1,\rho_2\}<\infty$.
	Then $\widetilde H\in\Ham$ and
	\begin{equation}\label{Z77}
		\lim_{r\to\infty}\mc A_r^{\ms b_1,\ms b_2}H = \widetilde H.
	\end{equation}
\end{lemma}

\begin{proof}
	To start with, note that
	\[
		\frac{\ms b_1(r)\ms b_2(r)}{r}
		= r\mr t(r)\sqrt{m_1\bigl(\mr t(r)\bigr)m_2\bigl(\mr t(r)\bigr)}
		= \mr t(r),
		\qquad r>0.
	\]
	It follows from \eqref{Z80}--\eqref{Z84}, \eqref{Z103}, the assumption in (ii),
	and the Uniform Convergence Theorems
	\cite[Theorems~1.2.1 and 2.4.1]{bingham.goldie.teugels:1989} that
	\begin{equation}\label{Z73}
		\lim_{r\to\infty} \int_0^x \bigl(\mc A_r^{\ms b_1\ms b_2}H\bigr)(s)\RD s
		= \int_0^x \widetilde H(s)\RD s
	\end{equation}
	holds locally uniformly for $x\in(0,\infty)$. Now \cref{Z112} implies
	that $\widetilde H\in\Ham$, and \cref{Z74}\,(ii) yields \eqref{Z77}.
\end{proof}

\noindent
The case when one of $m_1$ and $m_2$ is rapidly varying requires a
slightly different argument.
We elaborate on the case when $m_2$ is rapidly varying;
the case when $m_1$ is rapidly varying is completely analogous.

\begin{lemma}\label{Z129}
	Let $H$ be a Hamiltonian defined on $(0,\infty)$ such that neither $h_1$ nor $h_2$ 
	vanishes a.e.\ in a neighbourhood of $0$.
	Assume that $m_1$ is regularly varying at $0$ with index $\rho_1\in(0,\infty)$
	and that $m_2$ is rapidly varying at $0$.  Then
	\begin{equation}\label{Z130}
		\lim_{r\to\infty}\mc A_r^{\ms b_1,\ms b_2}H
		=
		\begin{pmatrix}
			\mathds{1}_{[0,1]} & 0
			\\[1ex]
			0 & \mathds{1}_{(1,\infty)}
		\end{pmatrix}
	\end{equation}
	and
	\begin{equation}\label{Z141}
		\lim_{t\to 0}\frac{m_3(t)}{\sqrt{m_1(t)m_2(t)}\,}=0.
	\end{equation}
\end{lemma}

\begin{proof}
	We start again with the formulae \eqref{Z80}--\eqref{Z84}
	and abbreviate the left-hand sides by $m_{r,1}(x)$, $m_{r,2}(x)$ and $m_{r,3}(x)$ respectively.
	Then
	\begin{equation}\label{Z138}
		\lim_{r\to\infty}m_{r,1}(x) = x^{\rho_1}
	\end{equation}
	for $x\in(0,\infty)$ locally uniformly (by the Uniform Convergence Theorem), and
	\[
		\lim_{r\to\infty}m_{r,2}(x)
		=
		\begin{cases}
			0, & x\in(0,1),
			\\[0.5ex]
			\infty, & x\in(1,\infty).
		\end{cases}
	\]
	Let $\mf t_r$ be again as in \eqref{Z75}.  Then we see that
	\[
		\lim_{r\to\infty}\mf t_r(x)
		=
		\begin{cases}
			x^{\rho_1}, & x\in(0,1),
			\\[0.5ex]
			\infty, & x\in(1,\infty),
		\end{cases}
	\]
	and therefore
	\begin{equation}\label{Z139}
		\lim_{r\to\infty}\mf t_r^{-1}(T)
		=
		\begin{cases}
			T^{\frac 1{\rho_1}}, & T\in(0,1),
			\\[0.5ex]
			1, & T\in(1,\infty).
		\end{cases}
	\end{equation}
	
	Let $r_n\to\infty$, and assume that the limit
	$\lim_{n\to\infty}\mc A_{r_n}^{\ms b_1,\ms b_2}H$ exists,
	which, by \cref{Z112}, is in $\Ham$.
	Denote the (unique) trace-normalised reparameterisation of this limit 
	by $\widehat H$ with entries $\hat h_i$.  Then
	\begin{equation}\label{Z140}
		\lim_{n\to\infty}\int_0^{\mf t_{r_n}^{-1}(T)}
		\bigl(\mc A_{r_n}^{\ms b_1,\ms b_2}H\bigr)(t)\RD t
		= \int_0^T \widehat H(t)\RD t
	\end{equation}
	locally uniformly for $T\in[0,\infty)$.
	Let $\hatM$ be the primitive of $\widehat H$ with entries $\hat m_i$.
	From \eqref{Z140}, \eqref{Z138} and \eqref{Z139} we obtain
	\[
		\hat m_1(T)
		= \lim_{n\to\infty}m_{r_n,1}\bigl(\mr t_{r_n}^{-1}(T)\bigr)
		=
		\begin{cases}
			T, & T\in(0,1),
			\\[0.5ex]
			1, & T\in(1,\infty).
		\end{cases}
	\]
	Thus $\hat h_1(t)=1$ for $t\in(0,1)$ a.e., and $\hat h_1(t)=0$ for $t\in(1,\infty)$ a.e.
	Since $\tr\widehat H=1$, it follows that $\hat h_2(t)=0$ for $t\in(0,1)$
	and $\hat h_2(t)=1$ for $t\in(1,\infty)$ a.e.
	Further, the non-negativity of $\widehat H$ implies that $\hat h_3=0$ a.e. 
	The limit relation \eqref{Z130} follows since $\NHam$ is compact;
	see, e.g.\ \cite[Lemma~2.9]{pruckner.woracek:limp}.
	In particular, we have
	\begin{equation}\label{Z147}
		\lim_{r\to\infty}\int_0^{\mf t_r^{-1}(T)}\bigl(\mc A_r^{\ms b_1,\ms b_2}H\bigr)\RD t
		= \hatM(T)
	\end{equation}
	for all $T\ge0$.

	To prove \eqref{Z141}, fix $T>1$.  For $i\in\{1,2\}$ we have
	\[
		\lim_{r\to\infty}\frac{m_i\bigl(\mr t(r)\mf t_r^{-1}(T)\bigr)}{m_i\bigl(\mr t(r)\bigr)}
		= \lim_{r\to\infty}m_{r,i}\bigl(\mf t_r^{-1}(T)\bigr)
		= \hat m_i(T)
	\]
	by \eqref{Z80}, \eqref{Z82} and \eqref{Z147}.
	It follows from \eqref{Z147} and \eqref{Z84} that
	\begin{align*}
		0 &= \hat m_3(T) = \lim_{r\to\infty}m_{r,3}\bigl(\mf t_r^{-1}(T)\bigr)
		\\[1ex]
		&= \lim_{r\to\infty}\Biggl[\frac{m_3\bigl(\mr t(r)\mf t_r^{-1}(T)\bigr)}{\sqrt{m_1\bigl(\mr t(r)\mf t_r^{-1}(T)\bigr)
		m_2\bigl(\mr t(r)\mf t_r^{-1}(T)\bigr)}\,}
		\cdot\underbrace{\sqrt{\frac{m_1\bigl(\mr t(r)\mf t_r^{-1}(T)\bigr)}{m_1\bigl(\mr t(r)\bigr)}\,}}_{
		\to\hat m_1(T)^{\frac 12}=1}
		\cdot\underbrace{\sqrt{\frac{m_2\bigl(\mr t(r)\mf t_r^{-1}(T)\bigr)}{m_2\bigl(\mr t(r)\bigr)}\,}}_{
		\to\hat m_2(T)^{\frac 12}=\sqrt{T-1}}\,\Biggr].
	\end{align*}
	Since the second and third factors on the right-hand side tend to non-zero numbers,
	the first factor must converge to $0$.
	The function $r\mapsto\mr t(r)\mf t_r^{-1}(T)$ is continuous, positive, and
	$\lim_{r\to\infty}\mr t(r)\mf t_r^{-1}(T)=0$.
	Hence the required limit relation \eqref{Z141} follows.
\end{proof}

\begin{proof}[Proof of \cref{Z10} `\textup{(ii)}$\bm\Rightarrow$\textup{(i)}'
and the formulae \eqref{Z98} and \eqref{Z86}--\eqref{Z17} in \cref{Z100}]
	\hfill
	\begin{Itemize}
	\item
		\Cref{Z96,Z14,Z129} show that
		\[
			q_H(rz) \sim \underbrace{\frac{\ms b_2(r)}{\ms b_1(r)}}_{=\ms a_H(r)}\cdot\,
			q_{\widetilde H}(z),
			\qquad r\to\infty,
		\]
		locally uniformly for $z\in\bb C^+$.
	\item
		It follows from \cref{Z132} (when $\max\{\rho_1,\rho_2\}<\infty$)
		or \cref{Z133} (when $\max\{\rho_1,\rho_2\}=\infty$) that
		\[
			q_{\widetilde H}(z) = i\omega\Bigl(\frac{z}{i}\Bigr)^\alpha,
			\qquad z\in\bb C^+,
		\]
		where $\alpha$ is as in \eqref{Z78}, and $\omega$ is as in \cref{Z132}
		if $\max\{\rho_1,\rho_2\}<\infty$ and $\omega=1$ otherwise.
		Note that that $\omega\ne0$ and that $\omega=\omega_{\alpha,\delta}$ where $\omega_{\alpha,\delta}$
		is as in \eqref{Z98}.
	\item
		From \cite[Theorem~3.1]{langer.woracek:kara} we obtain that $\ms a_H$ is regularly varying with index $\alpha$,
		that $\alpha\in[-1,1]$ and that $|\arg\omega_{\alpha,\delta}|\le\frac{\pi}{2}\bigl(1-|\alpha|\bigr)$.
	\item
		Relation \eqref{Z86} is just \eqref{Z101}, and \eqref{Z194} follows directly from \eqref{Z78}.
		
		For the proof of \eqref{Z97}, \eqref{Z98} and \eqref{Z6} assume first that $\max\{\rho_1,\rho_2\}<\infty$.
		The relation $\sqrt{1-\alpha^2}=\frac{2\sqrt{\rho_1\rho_2}\,}{\rho_1+\rho_2}$
		can be easily shown, from which we obtain (with the notation from \cref{Z132})
		\[
			\frac{\kappa}{\rho_3} = \frac{\sqrt{\rho_1\rho_2-(\delta\rho_3)^2}\,}{\rho_3}
			= \sqrt{\frac{\rho_1\rho_2}{\rho_3^2}-\delta^2}
			= \sqrt{1-\alpha^2-\delta^2}.
		\]
		This, together with \cref{Z132,Z88}, implies \eqref{Z98} and \eqref{Z6}.
		Since $\widetilde H\in\Ham$, we have $\rho_1\rho_2-(\delta\rho_3)^2\ge0$,
		which shows the inequality in \eqref{Z97}.
		In the case when $\max\{\rho_1,\rho_2\}=\infty$, relation \eqref{Z98}
		has already been shown, \eqref{Z6} is clear, and \eqref{Z97} follows
		from \cref{Z129}.
		
		Finally, in both cases, \eqref{Z79} and \eqref{Z11} follow from \eqref{Z98},
		and \eqref{Z17} follows from \eqref{Z38}.
	\end{Itemize}
\end{proof}

\subsection{Deducing the non-trace-normalised case}
\label{Z179}

Assume that we are given a Hamiltonian $H$ as in \cref{Z10}.
We pass to the trace-normalised reparameterisation of $H$,
i.e.\ we consider the Hamiltonian $\widehat H$ defined by
\[
	\widehat H(t) \DE H\bigl(\mf t^{-1}(t)\bigr)\cdot(\mf t^{-1})'(t),\qquad t\in(0,\infty),
\]
where $\mf t(x)=\int_0^x \tr H(s)\RD s$.
Let $\widehat m_j$ and similar notation have the corresponding meaning. Then
\[
	q_{\widehat H} = q_H, \qquad \hatM=M\circ\mf t^{-1},
\]
and, with $g(r)=\frac{1}{r^2}$,
\[
	\widehat{\mr t} = (\widehat m_1\widehat m_2)^{-1}\circ g
	= \bigl((m_1m_2)\circ\mf t^{-1}\bigr)^{-1}\circ g
	= \mf t\circ\mr t.
\]
Since $\mf t$ is regularly varying with positive index, also $\mf t^{-1}$
is regularly varying by \cref{Z254}.
Thus $\widehat m_j$ is regularly or rapidly varying if and only if $m_j$ is.

From these facts it is clear that neither (i) nor (ii)
in \cref{Z10} changes its truth value when passing from $H$ to $\widehat H$.
Further, $\hat\delta=\delta$ and
\[
	\ms a_{\widehat H} = \sqrt{\frac{\widehat m_1}{\widehat m_2}}\circ\widehat{\mr t}
	= \sqrt{\frac{m_1}{m_2}}\circ\mf t^{-1}\circ\mf t\circ\mr t
	= \ms a_H,
\]
and hence $\widehat\omega_{\alpha,\delta}=\omega_{\alpha,\delta}$ 
for $\widehat\omega_{\alpha,\delta}$ and $\omega_{\alpha,\delta}$ as in \eqref{Z57}.
Denote the index of $\mf t$ by $\sigma$.  Then $\widehat\rho_j=\frac{\rho_j}{\sigma}$
and hence $\widehat\alpha=\alpha$.
This shows that the formulae \eqref{Z86}--\eqref{Z17} hold true.

\subsection{The inverse theorem (proof of \cref{Z236})}
\label{Z243}

\begin{proof}[Proof of \cref{Z236}]
	Suppose first that (i) holds.  Without loss of generality we can assume that $\tr H(t)=1$ a.e.
	Then we can apply \cref{Z10}, which yields that $m_1$ and $m_2$ are regularly varying with $\ind m_i=\rho_i\in(0,\infty)$
	and the limit in \eqref{Z56} exists.
	By \cref{Z253} there exists a smoothly varying function $\eta:(0,L)\to(0,\infty)$ such that $m_1(t)\sim\eta(t)$ as $t\to0$
	and $\eta'(t)>0$ for all $t\in(0,L)$.
	Consider the reparameterisation $\widetilde H(t)\DE \frac{1}{\eta'(t)}H(\eta^{-1}(t))$ (see \cref{Z154})
	whose primitive is $\tildeM(t)=M(\eta^{-1}(t))$ by \eqref{Z36}.
	Then
	\[
		\widetilde m_1(t) = m_1\bigl(\eta^{-1}(t)\bigr) \sim t.
	\]
	Let $\mr{\tilde t}(r)$ be the solution of 
	\begin{equation}\label{Z237}
		(\widetilde m_1\widetilde m_2)\bigl(\mr{\tilde t}(r)\bigr) = \frac{1}{r^2}
	\end{equation}
	for $r>0$.
	From (i) in \cref{Z236}, \eqref{Z57} and \eqref{Z17} we obtain
	\[
		ie^{i\phi}\ms f(r) \sim q_H(ri) = q_{\widetilde H}(ri) 
		\sim i\omega_{\alpha,\delta}\,r\widetilde m_1\bigl(\mr{\tilde t}(r)\bigr) 
		\sim i\omega_{\alpha,\delta}\,r\mr{\tilde t}(r),
	\]
	which yields
	\[
		\mr{\tilde t}(r) \sim \frac{1}{|\omega_{\alpha,\delta}|}\cdot\frac{\ms f(r)}{r} \sim \ms g(r);
	\]
	note that $|\omega_{\alpha,\delta}|=C_{\alpha,\phi}$.
	By assumption, $\ind\ms g=\alpha-1<0$ and hence $\mr{\tilde t}^{-1}(t)\sim\ms g^{-1}(t)$
	by \cref{Z254}.
	It follows from \eqref{Z237} that
	\[
		t\widetilde m_2(t) \sim \widetilde m_1(t)\widetilde m_2(t) = \frac{1}{\bigl(\mr{\tilde t}^{-1}(t)\bigr)^2}
		\sim \frac{1}{(\ms g^{-1}(t))^2},
	\]
	which implies the second relation in \eqref{Z238}.
	The third relation in \eqref{Z238} follows from the first two and the limit relation \eqref{Z56}.
	Hence (ii) is satisfied.
	
	Now assume that (ii) holds.  Then $\widetilde m_1$ and $\widetilde m_2$ are regularly varying
	and the limit \eqref{Z56} exists.  Moreover,
	\[
		\frac{1}{r^2} = (\widetilde m_1\widetilde m_2)\bigl(\mr{\tilde t}(r)\bigr)
		\sim \frac{1}{\bigl[\ms g^{-1}\bigl(\mr{\tilde t}(r)\bigr)\bigr]^2},
	\]
	which implies that $\mr{\tilde t}(r)\sim\ms g(r)$ by \cref{Z254}.  From \eqref{Z17} we therefore obtain
	\begin{align*}
		q_H(ri) &= q_{\widetilde H}(ri) \sim i\omega_{\alpha,\delta}\,r\widetilde m_1\bigl(\mr{\tilde t}(r)\bigr)
		\sim i\omega_{\alpha,\delta}\,r\mr{\tilde t}(r)
		\\[1ex]
		&\sim i\omega_{\alpha,\delta}\,r\ms g(r) \sim i\omega_{\alpha,\delta}\,r\frac{1}{C_{\alpha,\phi}}\cdot\frac{\ms f(r)}{r} 
		= ie^{i\phi}\ms f(r),
	\end{align*}
	which shows that (i) is satisfied.
	
	Finally, let $\ms f$ be given as in the last part of the theorem.
	By \cref{Z253} we can choose $\ms g$ to be smoothly varying such that \eqref{Z246} holds.
	The function $\eta_2(t)\DE\frac{1}{t[\ms g^{-1}(t)]^2}$ is smoothly varying with index $\frac{1+\alpha}{1-\alpha}>0$
	and hence strictly increasing on $(0,t_0]$ with some $t_0>0$ and $\lim_{t\to0}\eta_2(t)=0$.  
	Define the Hamiltonian $H$ on $(0,\infty)$ by
	\[
		h_1(t) = 1, \qquad 
		h_2(t) =
		\begin{cases}
			\eta_2'(t), & t\in(0,t_0),
			\\[1ex]
			0, & t\in(t_0,\infty),
		\end{cases}
		\qquad h_3(t) = \frac{\delta}{\sqrt{1-\alpha^2}\,}\sqrt{h_2(t)}.
	\]
	The first two relations in \eqref{Z238} with $m_i=\widetilde m_i$ are clearly satisfied.
	When $\delta=0$, then also the third relation in \eqref{Z238} holds.
	Now assume that $\delta\ne0$.
	Since $\ind h_2+1=\frac{1+\alpha}{1-\alpha}$ and $\ind h_3+1=\frac{1}{1-\alpha}$, we have,
	for $t\in(0,t_0)$ and as $t\to0$,
	\[
		\frac{1}{t[\ms g^{-1}(t)]^2} = m_2(t) \sim \frac{1-\alpha}{1+\alpha}th_2(t),
		\qquad
		m_3(t) \sim \frac{\delta}{\sqrt{1-\alpha^2}\,}(1-\alpha)t\sqrt{h_2(t)},
	\]
	which implies that the third relation in \eqref{Z238} with $m_3=\widetilde m_3$ 
	also holds in the case when $\delta\ne0$.
	By what we have already shown in the first part of the proof, it follows that (i) is satisfied.
\end{proof}

\subsection{Proof of \cref{Z215}}
\label{Z171}

\begin{proof}[Proof of \cref{Z215}]
	Since $h_1(t)\ne0$ a.e., the function $m_1$ is strictly increasing and $m_1^{-1}$
	is absolutely continuous by \cite[Exercise 13 on p.~271]{natanson:1955}.
	Define the Hamiltonian $\widetilde H(t)\DE H(m_1^{-1}(t))(m_1^{-1})'(t)$
	(cf.\ \cref{Z154}) so that the corresponding primitive $\tildeM$ satisfies
	$\tildeM(t)=M(m_1^{-1}(t))$ by \eqref{Z36}.  Clearly,
	\[
		\widetilde m_1(t) = t, \qquad \widetilde m_2(t) = m_2\bigl(m_1^{-1}(t)\bigr) \ll t, \qquad 
		\widetilde{\mf t}(t) = \widetilde m_1(t)+\widetilde m_2(t) \sim t,
	\]
	and hence $\widetilde H$ satisfies the assumptions of \cref{Z10}.
	The statement in (i) is unaffected by rescaling.
	Moreover, (ii) in \cref{Z10} for $\widetilde H$ is equivalent to (ii)$'$.
	Hence the equivalence of (i) and (ii)$'$ follows from \cref{Z10}.
	
	Assume now that (i) and (ii)$'$ hold and let $\tilde t(r)$ be as in \eqref{Z226}.
	Then $\hat t(r)$ solves $(\widetilde m_1\widetilde m_2)(\hat t(r))=\frac{1}{r^2}$
	and hence $\hat t(r)=\mr t(r)$ where $\mr t(r)$ is as in \cref{Z10} applied to $\widetilde H$.
	Relation \eqref{Z57} yields
	\[
		q_H(rz) = q_{\widetilde H}(rz) \sim i\omega_{\alpha,\delta}\Bigl(\frac{z}{i}\Bigr)^\alpha\ms a_{\widetilde H}(r)
		\qquad \text{as} \ r\to\infty
	\]
	locally uniformly for $z\in\bb C^+$, and from \eqref{Z17} we obtain
	\[
		\ms a_{\widetilde H}(r) = r\widetilde m_1\bigl(\mr t(r)\bigr) = r\hat t(r).
	\]
	Further, $\ind(\widetilde m_1)=1$, $\ind(\widetilde m_2)=\rho$ and hence $\alpha=\frac{\rho-1}{\rho+1}$
	by \eqref{Z78}.
	
	Let us prove the last statement.  Relation \eqref{Z79} implies that $|\omega_{\alpha,\delta}|=C_{\alpha,\phi}$.
	It follows from \eqref{Z227} and \eqref{Z247} that $\hat t(r)\sim\ms g(r)$.
	Further, \eqref{Z226} can be rewritten as $t(m_2\circ m_1^{-1})(t) = \frac{1}{[\hat t^{-1}(t)]^2}$,
	which implies \eqref{Z248}.
\end{proof}

\subsection{The case of power asymptotics (proof of \cref{Z124})}
\label{Z92}

\begin{proof}[Proof of \cref{Z124}]
	Assume first that (ii) in the corollary is satisfied.
	Then statement (ii) in \cref{Z10} is true with $\delta=\frac{c_3}{\sqrt{c_1c_2}\,}$,
	and hence also (i) and the remaining statements are true.  Let us determine $\ms a_H$.
	By \eqref{Z38} we have $c_1c_2(\mr t(r))^{\rho_1+\rho_2}\sim\frac{1}{r^2}$
	and hence
	\[
		\mr t(r) \sim \bigl(c_1c_2r^2\bigr)^{-\frac{1}{\rho_1+\rho_2}},
		\qquad r\to\infty.
	\]
	This yields
	\[
		\ms a_H(r) \sim \sqrt{\frac{c_1}{c_2}\bigl(\mr t(r)\bigr)^{\rho_1-\rho_2}}
		\sim \sqrt{\frac{c_1}{c_2}\,}
		\bigl(c_1c_2r^2\bigr)^{\frac{\rho_2-\rho_1}{2(\rho_1+\rho_2)}}
		= c_1^{\frac{\alpha+1}{2}}c_2^{\frac{\alpha-1}{2}}r^\alpha
	\]
	where $\alpha=\frac{\rho_2-\rho_1}{\rho_2+\rho_1}$ as in the statement of the corollary.
	Note that $\alpha\in(-1,1)$ since $\rho_1$ and $\rho_2$ are finite.
	Now, with $z$ ranging in the compact set $\{e^{i\phi}\DS \psi\le\phi\le\pi-\psi\}$
	with $\psi\in(0,\frac{\pi}{2})$, we obtain from \eqref{Z57} that
	\[
		q_H(re^{i\phi}) \sim i\omega_{\alpha,\delta}\Bigl(\frac{e^{i\phi}}{i}\Bigr)^\alpha\ms a_H(r)
		\sim ic_1^{\frac{\alpha+1}{2}}c_2^{\frac{\alpha-1}{2}}\omega_{\alpha,\delta}
		\Bigl(\frac{re^{i\phi}}{i}\Bigr)^\alpha
	\]
	uniformly in $\phi\in[\psi,\pi-\psi]$.
	This shows that (i) and \eqref{Z137} in the corollary hold.

	Now assume that (i) is satisfied.
	Then (i) in \cref{Z10} and hence also (ii) holds.
	Since $\rho_1$ and $\rho_2$ are finite by \eqref{Z78} and the fact that $\alpha\in(-1,1)$,
	the functions $m_1$ and $m_2$ are regularly varying with indices $\rho_1$ and $\rho_2$
	respectively.
	By assumption, $\ms a_H(r)\sim dr^\alpha$ with some $d>0$.
	We distinguish two cases.
	\\[1ex]
	Case~1: $\rho_1\ne\rho_2$. \\
	We consider only the case when $\rho_1<\rho_2$; the case $\rho_1>\rho_2$ is
	completely analogous.  It follows from \eqref{Z18} that $\frac{m_2(t)}{m_1(t)}\to0$
	as $t\to0$.  Now the assumption that $m_1(t)+m_2(t)\sim ct^\sigma$ implies
	that $m_1(t)\sim ct^\sigma$.  This, together with \eqref{Z17}, shows that
	\[
		c\bigl(\mr t(r)\bigr)^\sigma \sim m_1\bigl(\mr t(r)\bigr)
		= \frac{\ms a_H(r)}{r} \sim dr^{\alpha-1}
	\]
	and hence $\mr t(r)\sim c'r^{\frac{\alpha-1}{\sigma}}$ with some $c'>0$, which,
	in turn, implies that $\mr t^{-1}(t)\sim c''t^{\frac{\sigma}{\alpha-1}}$
	with $c''>0$.  Hence
	\[
		m_2(t) = \frac{(m_1m_2)(t)}{m_1(t)} = \frac{1}{m_1(t)\bigl(\mr t^{-1}(t)\bigr)^2}
		\sim c_2t^{\rho_2}
	\]
	with $c_2>0$.
	\\[1ex]
	Case 2: $\rho_1=\rho_2$. \\
	In this case we have $\alpha=0$ and hence $\ms a_H(r)\to d$, which implies that
	\[
		\lim_{t\to0}\frac{m_1(t)}{m_2(t)}
		= \lim_{r\to\infty}\frac{m_1(\mr t(r))}{m_2(\mr t(r))} = d^2.
	\]
	Now the relation $m_1(t)+m_2(t)\sim ct^\sigma$ yields
	\[
		m_1(t) \sim \frac{cd^2}{d^2+1}t^\sigma, \qquad
		m_2(t) \sim \frac{c}{d^2+1}t^\sigma.
	\]

	Hence in both cases the first two limit relations in \eqref{Z135} hold.
	Since $\sqrt{m_1(t)m_2(t)} \sim \sqrt{c_1c_2}\,t^{\rho_3}$, the existence of the 
	limit in \eqref{Z56} implies that the third limit relation in \eqref{Z135} with $c_3=\delta\sqrt{c_1c_2}$
	holds.
\end{proof}

\section{The spectral measure}
\label{Z64}

In this section we translate the asymptotic behaviour of the Weyl coefficient to
the asymptotic behaviour of the spectral measure, more precisely, the distribution
functions $r\mapsto\mu_H((0,r))$ and $r\mapsto\mu_H((-r,0))$.
Throughout this section we assume that the following assumption is satisfied.

\begin{assumption}\label{Z165}
	Let $H$ be a Hamiltonian satisfying (a)--(c) at the beginning of the Introduction,
	let $m_j$ and $\mf t$ be as in \cref{Z206},
	and assume that neither $h_1$ nor $h_2$ vanishes a.e.\ on any neighbourhood of $0$
	and that $\mf t$ is regularly varying at $0$ with positive index.
	Further, let $\mu_H$ be the measure in the representation \eqref{Z4} of the 
	Weyl coefficient $q_H$.
	Recall also the notation $\ms a_H(r)$ from \eqref{Z109}.
\end{assumption}

\noindent
The first theorem provides information about the asymptotic behaviour of the spectral measure in many cases.

\begin{theorem}\label{Z170}
	Let \cref{Z165} be satisfied, assume that the equivalent 
	conditions \textup{(i)} and \textup{(ii)} in \cref{Z10} hold with $\alpha\in(-1,1)$,
	and set $\beta\DE\alpha+1$.  Then
	\begin{align}
		\lim_{r\to\infty}\frac{\mu_H\bigl((0,r)\bigr)}{r\ms a_H(r)}
		&= \frac{|\omega_{\alpha,\delta}|}{\pi\beta}\sin\Bigl(\frac{\pi\beta}{2}-\arg\omega_{\alpha,\delta}\Bigr),
		\label{Z207}
		\\[1ex]
		\lim_{r\to\infty}\frac{\mu_H\bigl((-r,0)\bigr)}{r\ms a_H(r)}
		&= \frac{|\omega_{\alpha,\delta}|}{\pi\beta}\sin\Bigl(\frac{\pi\beta}{2}+\arg\omega_{\alpha,\delta}\Bigr),
		\label{Z208}
	\end{align}
	where $\omega_{\alpha,\delta}$ is as in \eqref{Z98}.
\end{theorem}

\noindent
This theorem follows directly from \cref{Z10} and the Tauberian theorem \cite[Theorem~3.2]{langer.woracek:kara}.

\begin{remark}\label{Z91}
\rule{0ex}{1ex}
\begin{Enumerate}
\item
	If the right-hand side of \eqref{Z207} or \eqref{Z208} is non-zero, then
	the corresponding distribution function is regularly varying with index $\beta$, 
	and  the function $r\mapsto\mu_H((-r,r))$ has the same property if at least one of the right-hand sides is non-zero.
	Both right-hand sides (of \eqref{Z207} and \eqref{Z208}) vanish if and only if $\beta=1$  (i.e.\ $\alpha=0$)
	and $|\arg\omega_{\alpha,\delta}|=\frac{\pi}{2}$.  The relation $|\arg\omega_{\alpha,\delta}|=\frac{\pi}{2}$ 
	is equivalent to $|\delta|=1$.
\item
	There exist Hamiltonians $H$ such that (i) and (ii) in \cref{Z10} with $|\delta|=1$ are satisfied
	but $\mu_H((-\infty,0))=0$ and $r\mapsto\mu_H((0,r))$ is not regularly varying;
	this follows from \cite[Example~5.3]{langer.woracek:kara} combined with de~Branges' inverse theorem.
\item
	In the two cases when $q_H(ri)\to\zeta_0$, $\zeta_0\setminus\{0\}$ or $\alpha\in(-1,1)\setminus\{0\}$, $\delta=0$
	a similar result is proved in \cite[Corollaries~3.3 and 3.10]{eckhardt.kostenko.teschl:2018}.
\end{Enumerate}
\end{remark}

\noindent
In the remainder of this section we consider theorems that contain also inverse spectral results, i.e.\
where we obtain information about $H$ from the asymptotic behaviour of the spectral measure.
In the following $\beta$ denotes the index of regular variation of the distribution function $r\mapsto\mu_H((-r,r))$
if the latter is regularly varying.
Since $\mu_H$ satisfies $\int_{\bb R}\frac{\D\mu_H(t)}{1+t^2}<\infty$, we have $\beta\in[0,2]$;
see, e.g.\ \cite[Lemma~4.1]{langer.woracek:kara}.
It can happen that $\beta\ne\alpha+1$ where $\alpha$ is as in \cref{Z10};
see \cref{Z195}.
We distinguish several cases depending on $\beta$ and formulate the results in five theorems.
In most cases we can prove equivalences including the boundary cases $\beta=0$ and $\beta=2$.
The proofs are given after \cref{Z196}.

In the first two theorems we consider the cases $\beta=2$ and $\beta\in(1,2)$ respectively.

\begin{theorem}\label{Z166}
	Let \cref{Z165} be satisfied.
	Then the following statements are equivalent:
	\begin{Enumeratealph}
	\item
		$r\mapsto\mu_H((-r,r))$ is regularly varying with index $2$;
	\item
		$m_1$ is regularly varying and $m_2$ is rapidly varying at $0$.
	\end{Enumeratealph}
	If \textup{(a)} and \textup{(b)} are satisfied, then
	\begin{equation}\label{Z167}
		\ms a_H(r) \sim 2r\int_r^\infty \frac{\mu_H((-t,t))}{t^3}\RD t,
		\qquad r\to\infty.
	\end{equation}
\end{theorem}

\begin{theorem}\label{Z168}
	Let \cref{Z165} be satisfied.
	Then the following statements are equivalent:
	\begin{Enumeratealph}
	\item
		$r\mapsto\mu_H((-r,r))$ is regularly varying with index $\beta\in(1,2)$ and the limit
		\begin{equation}\label{Z197}
			\zeta \DE \lim_{r\to\infty}\frac{\mu_H((-r,0))}{\mu_H((0,r))}
		\end{equation}
		exists in $[0,\infty]$ \textup{(}the limit is interpreted as $\infty$ when $\mu_H((0,\infty))=0$\textup{)};
	\item
		$m_1$ and $m_2$ are regularly varying at $0$ with indices $\rho_1,\rho_2$ respectively 
		with $\rho_1<\rho_2$, and the limit
		\[
			\delta \DE \lim_{t\to0}\frac{m_3(t)}{\sqrt{m_1(t)m_2(t)}\,}
		\]
		exists.
	\end{Enumeratealph}
	If \textup{(a)} and \textup{(b)} are satisfied, then $\beta=\frac{2\rho_2}{\rho_1+\rho_2}$,
	\begin{align}
		\ms a_H(r) &\sim \frac{c_{\beta,\zeta}}{|\omega_{\beta-1,\delta}|}\cdot\frac{\mu_H((-r,r))}{r},
		\qquad r\to\infty,
		\label{Z186}
		\\[1ex]
		\log\zeta &= \frac{\pi(\beta-1)\delta}{\sqrt{1-(\beta-1)^2-\delta^2}\,},
		\label{Z187}
	\end{align}
	where $\omega_{\beta-1,\delta}$ is as in \eqref{Z98} and
	\[
		c_{\beta,\zeta} = \frac{\pi\beta}{2}\biggl[
		\Big|\cot\frac{\pi\beta}{2}\Big|+\Bigl(\frac{\zeta-1}{\zeta+1}\Bigr)^2\Big|\tan\frac{\pi\beta}{2}\Big|\biggr].
	\]
	The right-hand side of \eqref{Z187} is interpreted as $\pm\infty$ when $\delta=\pm\sqrt{1-(\beta-1)^2}$
	and $\log\zeta$ is interpreted as $\infty$ when $\zeta=\infty$ and as $-\infty$ when $\zeta=0$.
\end{theorem}

\medskip

\noindent
The case when $\beta<1$ is slightly more complicated since the constant $a_H$ in \eqref{Z4}
dominates the integral if the former is non-zero.
This applies also to some cases when $\beta=1$.
Let us introduce the following notation: for $\gamma\in\bb R$ set
\begin{equation}\label{Z175}
	H^{(\gamma)}(t) \DE \begin{pmatrix} 1 & -\gamma \\[0.5ex] 0 & 1 \end{pmatrix}H(t)
	\begin{pmatrix} 1 & 0 \\[0.5ex] -\gamma & 1 \end{pmatrix}.
\end{equation}
Then the primitive of $H^{(\gamma)}$ is given by
\begin{equation}\label{Z176}
	M^{(\gamma)} = \begin{pmatrix} m_1^{(\gamma)} & m_3^{(\gamma)} \\[1ex] m_3^{(\gamma)} & m_2^{(\gamma)} \end{pmatrix}
	= \begin{pmatrix} m_1-2\gamma m_3+\gamma^2m_2 & m_3-\gamma m_2 \\[1ex] m_3-\gamma m_2 & m_2 \end{pmatrix}.
\end{equation}
It is well known that the Weyl coefficient of $H^{(\gamma)}$ satisfies
\begin{equation}\label{Z177}
	q_{H^{(\gamma)}}(z) = q_H(z)-\gamma;
\end{equation}
see, e.g.\ \cite[Lemma~3.2]{winkler:1995}.

The next two theorems deal with the cases $\beta=0$ and $\beta\in(0,1)$ respectively.

\begin{theorem}\label{Z172}
	Let \cref{Z165} be satisfied and let $H^{(\gamma)}$ and $m_i^{(\gamma)}$ be as 
	in \eqref{Z175} and \eqref{Z176} respectively.
	Then the following statements are equivalent:
	\begin{Enumeratealph}
	\item
		$r\mapsto\mu_H((-r,r))$ is slowly varying;
	\item
		there exists $\gamma\in\bb R$ such that
		\begin{Itemize}
		\item
			$\displaystyle\lim_{t\to0}\frac{m_1(t)}{m_2(t)} = \gamma^2$,
		\item
			$m_1^{(\gamma)}$ is rapidly varying at $0$.
		\end{Itemize}
	\end{Enumeratealph}
	If \textup{(a)} and \textup{(b)} are satisfied, then $\lim_{t\to0}\frac{m_3(t)}{m_2(t)}=\gamma$ and
	\begin{equation}\label{Z178}
		\ms a_{H^{(\gamma)}}(r) \sim \frac{\mu_H((-r,r))}{r}, \qquad r\to\infty.
	\end{equation}
\end{theorem}

\begin{theorem}\label{Z173}
	Let \cref{Z165} be satisfied, let $H^{(\gamma)}$ and $m_i^{(\gamma)}$ be as 
	in \eqref{Z175} and \eqref{Z176} respectively,
	and set $\sigma\DE\ind(\tr M)$.
	Then the following statements are equivalent:
	\begin{Enumeratealph}
	\item
		$r\mapsto\mu_H((-r,r))$ is regularly varying with index $\beta\in(0,1)$ and the limit
		in \eqref{Z197} exists in $[0,\infty]$;
	\item
		there exists $\gamma\in\bb R$ such that
		\begin{Itemize}
		\item
			$\displaystyle\lim_{t\to0}\frac{m_1(t)}{m_2(t)} = \gamma^2$,
		\item
			$m_1^{(\gamma)}$ is regularly varying at $0$ with $\tilde\rho_1\DE\ind m_1^{(\gamma)}>\sigma$,
		\item
			the limit 
			\begin{equation}\label{Z199}
				\delta \DE \lim_{t\to0}\frac{m_3^{(\gamma)}(t)}{\sqrt{m_1^{(\gamma)}(t)m_2^{(\gamma)}(t)}\,}
			\end{equation}
			exists.
		\end{Itemize}
	\end{Enumeratealph}
	If \textup{(a)} and \textup{(b)} are satisfied, then $m_2$ is regularly varying with $\ind m_2=\sigma$,
	$\beta=\frac{2\sigma}{\tilde\rho_1+\sigma}$ and
	the relations in \eqref{Z186} and \eqref{Z187} with $\ms a_H$ replaced by $\ms a_{H^{(\gamma)}}$ hold.
\end{theorem}

\begin{remark}\label{Z195}
	Assume that (a) and (b) in \cref{Z172} or \cref{Z173} hold with $\gamma\ne0$.
	Then $q_H(ri)\to\gamma$ as $r\to\infty$ (see \cref{Z174}) and hence $\alpha=0$ in \cref{Z10},
	which shows that $\beta\ne\alpha+1$ in this situation.
	Moreover, the right-hand sides of \eqref{Z207} and \eqref{Z208} vanish.
\end{remark}

\noindent
It remains to consider the case when $\beta=1$.  Here we do not have equivalent statements;
see \cref{Z91}.  For the Tauberian direction one can use \cref{Z170} in many cases, namely when $|\delta|<1$.
The following theorem contains an Abelian implication in the situation where the measure on
one half-axis dominates the measure on the other half-axis.

\begin{theorem}\label{Z196}
	Let \cref{Z165} be satisfied and let $H^{(\gamma)}$ and $m_i^{(\gamma)}$ be as 
	in \eqref{Z175} and \eqref{Z176} respectively.
	Assume that $r\mapsto\mu_H((-r,r))$ is regularly varying with index $1$ and 
	the limit in \eqref{Z197} exists in $[0,\infty]$ with $\zeta\ne1$.
	\begin{Enumerate}
	\item
		If $\int_{\bb R}\frac{\RD\mu_H(t)}{1+|t|}=\infty$, then
		\begin{Itemize}
		\item
			$m_1$ and $m_2$ are regularly varying with equal indices;
		\item
			$\displaystyle \ms a_H(r) \sim \frac{|\zeta-1|}{\zeta+1}\int_1^r\frac{\mu_H((-t,t))}{t^2}\,\RD t$,
			\qquad $r\to\infty$;
		\item
			$\displaystyle \sgn(\zeta-1) = \lim_{t\to0}\frac{m_3(t)}{\sqrt{m_1(t)m_2(t)}\,}$.
		\end{Itemize}
	\item
		If $\int_{\bb R}\frac{\RD\mu_H(t)}{1+|t|}<\infty$, then there exists $\gamma\in\bb R$ such that
		\begin{Itemize}
		\item
			$\displaystyle\lim_{t\to0}\frac{m_1(t)}{m_2(t)} = \gamma^2$;
		\item
			$m_1^{(\gamma)}$ is regularly varying at $0$ with $\ind m_1^{(\gamma)}=\ind(\tr M)$;
		\item
			$\displaystyle \ms a_{H^{(\gamma)}}(r) \sim \frac{|\zeta-1|}{\zeta+1}\int_r^\infty\frac{\mu_H((-t,t))}{t^2}\,\RD t$,
			\qquad $r\to\infty$;
		\item
			$\displaystyle \sgn(\zeta-1) = -\lim_{t\to0}\frac{m_3^{(\gamma)}(t)}{\sqrt{m_1^{(\gamma)}(t)m_2^{(\gamma)}(t)}\,}$.
		\end{Itemize}
	\end{Enumerate}
\end{theorem}

\bigskip

\noindent
Let us now come to the proofs of Theorems~\ref{Z166}--\ref{Z196}.

\begin{proof}[Proof of \cref{Z166}]
	It follows from \cite[Theorems~5.1 and 3.1]{langer.woracek:kara}, \cref{Z10}
	and \eqref{Z78}, \eqref{Z98} (with $\alpha=1$) that both (a) and (b) are equivalent to
	the fact that there exists a regularly varying function $\ms a$ such that \eqref{Z110}
	holds with $\omega=1$.  The above mentioned theorems also imply that one can choose either
	\[
		\ms a(r) = \ms a_H(r)
		\qquad\text{or}\qquad
		\ms a(r) = 2r\int_r^\infty\frac{\mu_H((-t,t))}{t^3}\RD t,
	\]
	which implies that \eqref{Z167} holds if (a) and (b) are satisfied.
\end{proof}

\begin{proof}[Proof of \cref{Z168}]
	The equivalence of (a) and (b) is shown in a similar way as in the proof of \cref{Z166},
	again based on \cref{Z10} and \cite[Theorem~5.1]{langer.woracek:kara}.
	Now assume that (a) and (b) are satisfied.  Then $\beta=\alpha+1=\frac{2\rho_2}{\rho_1+\rho_2}$
	by \eqref{Z78},
	\[
		q_H(ri) \sim i\omega_{\alpha,\delta}\ms a_H(r)
	\]
	by \cref{Z10} with $\omega_{\alpha,\delta}$ as in \eqref{Z98}, and
	\[
		q_H(ri) \sim i\hat\omega_{\beta,\zeta}\frac{\pi\beta}{|\sin(\pi\beta)|}\cdot\frac{\mu_H((-r,r))}{r}
	\]
	with
	\begin{equation}\label{Z185}
		\hat\omega_{\beta,\zeta} = -\biggl(\cos\frac{\pi\beta}{2}+i\frac{\zeta-1}{\zeta+1}\cdot\sin\frac{\pi\beta}{2}\biggr),
	\end{equation}
	by \cite[Theorem~5.1]{langer.woracek:kara}.
	These relations imply that $\arg\omega_{\alpha,\delta}=\arg\hat\omega_{\beta,\zeta}$ and
	\begin{equation}\label{Z190}
		\ms a_H(r) \sim \frac{|\hat\omega_{\beta,\zeta}|}{|\omega_{\alpha,\delta}|}\cdot\frac{\pi\beta}{|\sin(\pi\beta)|}
		\cdot\frac{\mu_H((-r,r))}{r}
		= \frac{c_{\beta,\zeta}}{|\omega_{\alpha,\delta}|}\cdot\frac{\mu_H((-r,r))}{r},
	\end{equation}
	which is \eqref{Z186}.
	To show \eqref{Z187}, assume first that $\zeta\in(0,\infty)$.  Then 
	\begin{equation}\label{Z188}
		\tan(\arg\hat\omega_{\beta,\zeta}) = \tan\frac{\pi\beta}{2}\cdot\frac{\zeta-1}{\zeta+1}
		= \tan\frac{\pi\beta}{2}\cdot\tanh\frac{\log\zeta}{2}
	\end{equation}
	by \eqref{Z185}.  When $|\delta|<\sqrt{1-\alpha^2}$ we obtain from \eqref{Z6} that
	\begin{equation}\label{Z189}
		\tan(\arg\omega_{\alpha,\delta}) = -\tan\Bigl(\frac{\pi}{2}(1-\alpha)\Bigr)
		\cdot\tanh\frac{\pi\alpha\delta}{2\sqrt{1-\alpha^2-\delta^2}\,}
		= \tan\frac{\pi\beta}{2}\cdot\tanh\frac{\pi\alpha\delta}{2\sqrt{1-\alpha^2-\delta^2}\,}.
	\end{equation}
	From \eqref{Z188} and \eqref{Z189} we can deduce \eqref{Z187} when $|\delta|<\sqrt{1-\alpha^2}$.
	The cases $\delta=\sqrt{1-\alpha^2}$ and $\delta=-\sqrt{1-\alpha^2}$ correspond to $\zeta=\infty$
	and $\zeta=-\infty$ respectively.
\end{proof}

\medskip

\noindent
For the next proofs we need the following lemma.

\begin{lemma}\label{Z174}
	Let $q_H$ be the Weyl coefficient of a Hamiltonian $H$.
	Then the following equivalences and implications hold.
	\begin{Enumerate}
	\item
		$\lim\limits_{r\to\infty}q_H(ri)=0$ \quad $\Leftrightarrow$ \quad $\lim\limits_{t\to0}\dfrac{m_1(t)}{m_2(t)}=0$
		\quad $\Rightarrow$ \quad $\lim\limits_{t\to0}\dfrac{m_3(t)}{m_2(t)}=0$.
	\item
		For $\gamma\in\bb R\setminus\{0\}$ we have
		\begin{equation}\label{Z198}
			\lim_{r\to\infty}q_H(ri)=\gamma \quad\Leftrightarrow\quad
			\biggl(\lim_{t\to0}\frac{m_3(t)}{m_2(t)}=\gamma \quad\text{and}\quad
			\lim_{t\to0}\frac{m_1(t)}{m_2(t)}=\gamma^2\biggr).
		\end{equation}
	\end{Enumerate}
\end{lemma}

\begin{proof}
	(i)
	The equivalence follows directly from \cite[Theorem~1.1]{langer.pruckner.woracek:heniest};
	see also \cite[Remark~1.2\,(vi)]{langer.pruckner.woracek:heniest}.
	For the second implication one uses \eqref{Z103}.

	(ii) The limit relations are unaffected by reparameterisations.
	Hence we can assume that $\tr H(t)=1$ a.e.
	It follows from \cref{Z259} 
	that the limit relation on the left-hand side of $\Leftrightarrow$ in \eqref{Z198} is equivalent to
	\[
		\lim_{t\to0}\frac{m_i(t)}{t}=c_i, \quad i\in\{1,2,3\}, \qquad
		c_1,c_2>0, \qquad \gamma=\frac{c_3}{c_2}, \qquad \frac{c_1}{c_2}=\Bigl(\frac{c_3}{c_2}\Bigr)^2,
	\]
	which, in turn, is equivalent to the limit relations on the right-hand side of \eqref{Z198} 
	since $m_1(t)+m_2(t)=t$.
\end{proof}

\begin{proof}[Proof of \cref{Z172,Z173}]
	Assume first that (a) in \cref{Z172} or \cref{Z173} holds.
	For the case of \cref{Z172} we set $\beta=0$.
	Then $\int_{\bb R}\frac{\RD\mu_H(t)}{1+|t|}<\infty$ 
	(see, e.g.\ \cite[Lemma~4.1]{langer.woracek:kara} and \eqref{Z27}) and hence
	\begin{equation}\label{Z204}
		q_H(z) = \gamma+\int_{\bb R}\frac{1}{t-z}\,\RD\mu_H(t)
	\end{equation}
	with some $\gamma\in\bb R$, which, in turn, 
	implies that $q_H(ri)\to\gamma$ as $r\to\infty$.
	Hence we can apply \cref{Z174}, which yields
	\begin{align}
		& m_2(t) \sim \frac{1}{\gamma^2+1}\tr M(t),
		\label{Z200}
		\\[1ex]
		& \frac{\tr M^{(\gamma)}(t)}{m_2(t)} 
		= \frac{m_1(t)}{m_2(t)}-2\gamma\frac{m_3(t)}{m_2(t)}+\gamma^2+1
		\to 1
		\nonumber
	\end{align}
	as $t\to0$, and therefore $m_2$ and $\tr M^{(\gamma)}$ are regularly varying 
	with index $\sigma=\ind(\tr M)>0$.
	Relations \eqref{Z204} and \eqref{Z177} imply that 
	\[
		q_{H^{(\gamma)}}(z) = \int_{\bb R}\frac{1}{t-z}\,\RD\mu_H(t),
	\]
	and hence we can apply \cite[Theorem~5.1]{langer.woracek:kara} to deduce that
	\[
		q_{H^{(\gamma)}}(ri) \sim i\hat\omega_{\beta,\zeta}\frac{\pi\beta}{|\sin(\pi\beta)|}\cdot\frac{\mu_H((-r,r))}{r}
	\]
	with
	\[
		\hat\omega_{\beta,\zeta} = \cos\frac{\pi\beta}{2}+i\frac{\zeta-1}{\zeta+1}\cdot\sin\frac{\pi\beta}{2}
	\]
	and $\frac{\pi\beta}{|\sin(\pi\beta)|}$ is interpreted as $1$ when $\beta=0$.
	Since $r\mapsto\frac{\mu_H((-r,r))}{r}$ is regularly varying with index $\beta-1$,
	we obtain from \cref{Z10} that $m_1^{(\gamma)}$ is rapidly or regularly varying 
	with $\ind m_1^{(\gamma)}>\ind m_2=\ind\tr M$
	and that the limit in \eqref{Z199} exists in the case of \cref{Z173}.
	This shows that (b) is satisfied in \cref{Z172} or \cref{Z173} respectively.
	
	Now assume that (b) holds in \cref{Z172} or \cref{Z173}.
	It follows that \eqref{Z200} is satisfied and hence $\ind m_1^{(\gamma)}>\ind\mf t=\ind m_2$,
	which, in turn, implies that
	\begin{equation}\label{Z201}
		\frac{m_1(t)}{m_2(t)}-2\gamma\frac{m_3(t)}{m_2(t)}+\gamma^2 = \frac{m_1^{(\gamma)}(t)}{m_2(t)} \to 0
	\end{equation}
	as $t\to0$.  This shows, in particular, that $\tr M^{(\gamma)}$ is regularly varying with positive index.
	Hence we can apply \cref{Z10}, which yields that \eqref{Z57} holds.
	We also obtain from \eqref{Z201} that $\frac{m_3(t)}{m_2(t)}\to\gamma$ when $\gamma\ne0$.
	By \cref{Z174} we have $q_H(ri)\to\gamma$ and hence $q_{H^{(\gamma)}}(ri)\to0$ as $r\to\infty$.
	Therefore we can apply \cite[Theorem~5.1]{langer.woracek:kara},
	which implies that (a) holds.
	
	In a similar way as in the proofs of \cref{Z166,Z168} one shows the remaining
	statements under the assumption that (a) and (b) hold.
\end{proof}

\begin{proof}[Proof of \cref{Z196}]
	(i)
	It follows from \cite[Theorem~5.1]{langer.woracek:kara} that
	\begin{equation}\label{Z116}
		q_H(ri) \sim \frac{\zeta-1}{\zeta+1}\int_1^r\frac{\mu_H((-t,t))}{t^2}\,\RD t,
		\qquad r\to\infty,
	\end{equation}
	and that the integral on the right-hand side is slowly varying.
	Hence, by \cref{Z10}, $m_1$ and $m_2$ are regularly varying with equal indices
	and the limit in \eqref{Z56} exists.
	Combining \eqref{Z116} with \eqref{Z57} we obtain that $\omega_{\alpha,\delta}\in i\bb R$ and 
	hence $\omega_{\alpha,\delta}=-i\delta\in\{i,-i\}$ by \eqref{Z98}.
	From this all assertions follow.

	(ii)
	The proof is very similar with the additional transformation from \eqref{Z175}.
\end{proof}

\section{Scalar differential equations}
\label{Z65}

In this section we apply our results to three scalar differential equations:
Sturm--Liouville equations, Krein strings and generalised indefinite strings.
In all three cases we rewrite the differential equation as a canonical system
and characterise when the corresponding Weyl coefficient is asymptotically equal
to a constant times a regularly varying function along the imaginary axis.
We also prove some inverse results.

\subsection{Sturm--Liouville equations}
\label{Z66}

Let $p,q,w$ be real-valued functions on $(0,L)$ such that $p(x)\ne0,w(x)>0$ a.e.\ 
and $\frac{1}{p},q,w\in L^1_{\textsf{loc}}([0,L))$,
and consider the Sturm--Liouville equation
\begin{equation}\label{Z210}
	-(pu')'+qu = zwu
\end{equation}
on $(0,L)$, where $z\in\bb C$.  Assume that \eqref{Z210} is in the limit point case at $L$, i.e.\ for non-real $z$,
equation \eqref{Z210} has only one linearly independent solution in the space $L^2(w)\DE L^2((0,L),w(x)\DD x)$.
One could also consider the limit circle situation or a regular endpoint and impose a boundary condition at $L$;
this would require only slight modifications.
The Titchmarsh--Weyl coefficient (corresponding to a Dirichlet boundary condition at $0$) is defined as
\[
	\qSL(z) \DE \frac{\psi^{[1]}(0;z)}{\psi(0;z)}, \qquad z\in\bb C^+,
\]
where $\psi(\cdot;z)$ is a non-zero solution of \eqref{Z210} in $L^2(w)$ and $u^{[1]}\DE pu'$ 
denotes the first quasi-derivative of $u$.

We can rewrite \eqref{Z210} as a canonical system as follows.  Let $\theta$ and $\varphi$ be solutions
of \eqref{Z210} with $z=0$ that satisfy the initial conditions
\begin{equation}\label{Z217}
\begin{alignedat}{2}
	\theta(0) &= 1, \qquad & \varphi(0) &= 0,
	\\[1ex]
	\theta^{[1]}(0) &= 0, \qquad & \varphi^{[1]}(0) &= -1,
\end{alignedat}
\end{equation}
and define the Hamiltonian
\begin{equation}\label{Z216}
	H = w
	\begin{pmatrix}
		\theta^2 & \theta\varphi \\[1ex] \theta\varphi & \varphi^2
	\end{pmatrix}.
\end{equation}
If \eqref{Z210} is in the limit circle case or regular at $L$, one has to extend $H$ with
a constant $H$ with $\det H=0$ to the right of $L$.
One can easily show that if $y=\binom{y_1}{y_2}$ is a solution of \eqref{Z1},
then $u=\theta y_1+\varphi y_2$ is a solution of \eqref{Z210}
and $\int_0^L |u(t)|^2w(t)\DD t=\int_0^L y(t)^*H(t)y(t)\DD t$.
Moreover, since $u(0)=y_1(0)$ and $u^{[1]}(0)=-y_2(0)$, the Weyl coefficients satisfy
\[
	\qSL(z) = q_H(z), \qquad z\in\bb C^+;
\]
(note that the Weyl coefficient of a canonical system satisfies $q_H(z)=-\frac{y_2(0;z)}{y_1(0;z)}$
if $y(\cdot;z)$ is a solution of \eqref{Z1} that satisfies $\int_0^L y^*Hy<\infty$).

The following theorem gives a characterisation of the situation when the Titchmarsh--Weyl coefficient
behaves like a constant times a regularly varying function along the imaginary axis.
The theorem contains also an inverse result.

\begin{theorem}\label{Z212}
	Let $p,q,w$ be as at the beginning of the section, set
	\begin{equation}\label{Z221}
		P(x) \DE \int_0^x \frac{1}{p(t)}\DD t, \qquad
		W(x) \DE \int_0^x w(t)\DD t, \qquad
		x\in[0,L),
	\end{equation}
	and assume that
	\begin{equation}\label{Z214}
		\max_{s\in[0,t]}|P(s)| \asymp |P(t)|
	\end{equation}
	as $t\to0$.
	Then the following statements are equivalent:
	\begin{Enumeratealph}
	\item
		there exist a regularly varying function $\ms a:(0,\infty)\to(0,\infty)$ and 
		a constant $\omega\in\bb C\setminus\{0\}$ such that
		\[
			\qSL(ri) \sim i\omega\ms a(r) \qquad \text{as} \ r\to\infty;
		\]
	\item
		the mapping
		\[
			x \mapsto G(x) \DE \int_0^x \bigl[P\bigl(W^{-1}(s)\bigr)\bigr]^2\DD s
		\]
		is regularly or rapidly varying at $0$, and, if $G$ is regularly varying, the limit
		\begin{equation}\label{Z220}
			\delta \DE -\lim_{x\to0}\frac{\int_0^x P(W^{-1}(s))\DD s}{\sqrt{xG(x)}}
		\end{equation}
		exists.
	\end{Enumeratealph}
	Assume that \textup{(a)} and \textup{(b)} hold and that $G$ is regularly varying.  
	Further, let $\hat x(r)$ be the solution of
	\begin{equation}\label{Z229}
		\hat x(r)G\bigl(\hat x(r)\bigr) = \frac{1}{r^2}
	\end{equation}
	for $r>0$ and set $\alpha\DE\ind\ms a$ and $\rho\DE\ind G$.  Then $\alpha=\frac{\rho-1}{\rho+1}$ and
	\begin{equation}\label{Z234}
		\qSL(rz) \sim i\omega_{\alpha,\delta}\Bigl(\frac{z}{i}\Bigr)^\alpha r\hat x(r) \qquad \text{as} \ r\to\infty
	\end{equation}
	locally uniformly for $z\in\bb C^+$ where $\omega_{\alpha,\delta}$ is as in \eqref{Z98}.

	Moreover, the asymptotic behaviour of $G$ can be recovered from the behaviour of $\qSL$:
	let $\phi\in\bigl[-\frac{\pi}{2},\frac{\pi}{2}\bigr]$ be such that
	\[
		\qSL(ri) \sim ie^{i\phi}|\qSL(ri)| \qquad \text{as} \ r\to\infty,
	\]
	let $\alpha$ be the index of $r\mapsto|\qSL(ri)|$, let $C_{\alpha,\phi}$
	be given the right-hand side of \eqref{Z11},
	and let $\ms g$ be a strictly decreasing and regularly varying function such that
	\begin{equation}\label{Z261}
		\ms g(r) \sim \frac{1}{C_{\alpha,\phi}}\cdot\frac{|\qSL(ri)|}{r} \qquad \text{as} \ r\to\infty;
	\end{equation}
	then
	\begin{equation}\label{Z263}
		G(x) \sim \frac{1}{x[\ms g^{-1}(x)]^2} \qquad \text{as} \ x\to0.
	\end{equation}
\end{theorem}

\medskip

\noindent
Before we prove the theorem, let us state some comments and a corollary.

\begin{remark}\label{Z223}
\rule{0ex}{1ex}
\begin{Enumerate}
\item
	It follows from \eqref{Z214} that either $P(t)>0$ for all $t\in(0,L)$
	or $P(t)<0$ for all $t\in(0,L)$.
\item
	Clearly, \eqref{Z214} is satisfied if $p(t)>0$ a.e.\ or if $p(t)<0$ a.e.
\end{Enumerate}
\end{remark}

\begin{corollary}\label{Z222}
	Let $p,q,w$ be as at the beginning of the section, let $P$ and $W$ be as in \eqref{Z221} and 
	assume that \eqref{Z214} holds.
	Consider the statement
	\begin{Enumeratealph}
	\item[\textup{(b)$'$}]
		the mapping $|P\circ W^{-1}|$ is regularly or rapidly varying at $0$,
	\end{Enumeratealph}
	and let \textup{(a)} be as in \cref{Z212}.  Then the following statements are true.
	\begin{Enumerate}
	\item
	\textup{(b)$'$}\,$\Rightarrow$\,\textup{(a)}.
	\item
	If, in addition, $p(t)>0$ a.e.\ or $p(t)<0$ a.e., then also \textup{(a)}\,$\Rightarrow$\,\textup{(b)$'$}.
	\end{Enumerate}
	Assume that \textup{(a)} and \textup{(b)$'$} hold and that $|P\circ W^{-1}|$ 
	is regularly varying at $0$ with index $\gamma$.
	Let $\tilde t(r)$ be the solution of
	\begin{equation}\label{Z231}
		W\bigl(\tilde t(r)\bigr)\big|P\bigl(\tilde t(r)\bigr)\big| = \frac{1}{r}
	\end{equation}
	for $r>0$.  Then $\alpha\DE\ind\ms a=\frac{\gamma}{\gamma+1}$ and
	\begin{equation}\label{Z235}
		\qSL(rz) \sim -\varepsilon(-\varepsilon z)^\alpha \alpha^\alpha(1-\alpha)^\alpha
		\frac{\Gamma(1-\alpha)}{\Gamma(1+\alpha)}\cdot\frac{1}{|P(\tilde t(r))|}
		\qquad \text{as} \ r\to\infty
	\end{equation}
	locally uniformly for $z\in\bb C^+$, where $\varepsilon\DE\sgn P(t)$, $t\in(0,L)$, 
	which is well defined by \cref{Z223}\,\textup{(i)}, and $\alpha^\alpha$ is interpreted as $1$ when $\alpha=0$.

	Moreover, the asymptotic behaviour of $|(P\circ W^{-1}|$ can be recovered from the behaviour of $\qSL$:
	let
	\begin{equation}\label{Z262}
		C_{\alpha,\phi} = \alpha^\alpha(1-\alpha^2)^{\frac{1+\alpha}{2}}\frac{\Gamma(1-\alpha)}{\Gamma(2+\alpha)}
	\end{equation}
	and let $\ms g$ be a strictly decreasing and regularly varying function such that \eqref{Z261} holds; then
	\begin{equation}\label{Z264}
		\big|\bigl(P\circ W^{-1}\bigr)(x)\big| \sim \sqrt{\frac{1+\alpha}{1-\alpha}}\cdot\frac{1}{x\ms g^{-1}(x)}.
	\end{equation}
\end{corollary}

\begin{remark}\label{Z225}
\rule{0ex}{1ex}
\begin{Enumeratealph}
\item
	In \cite[Theorem~4.1]{bennewitz:1989} the implication (b)$'$\,$\Rightarrow$\,(a) is proved
	under the stronger assumption $\max_{s\in[0,t]}P(s)\sim P(t)$ instead of \eqref{Z214}.
	There is no converse implication in that paper such as (a)\,$\Rightarrow$\,(b) in \cref{Z212}.
\item
	Section~5 in \cite{bennewitz:1989} contains some partial inverse results under quite restrictive assumptions.
\item
	If $p(x)>0$ a.e., then one can transform the Sturm--Liouville equation with a change of the independent variable
	to obtain $p\equiv1$.  In this situation one has a one-to-one correspondence between the asymptotic behaviours
	of $\qSL$ and $W$.
\end{Enumeratealph}
\end{remark}

\medskip

\noindent
Before we prove \cref{Z212,Z222}, we need a lemma.

\begin{lemma}\label{Z211}
	Assume that \eqref{Z214} holds and let $\varphi$ be as at the beginning of this section.
	Then
	\[
		\varphi(t) \sim -P(t) \qquad \text{as} \ t\to 0.
	\]
\end{lemma}

\begin{proof}
	One can rewrite the differential equation $(p\varphi')'=q\varphi$
	on the interval $(0,t)$ with $t>0$
	in a standard way as $v(x)=v_0+(T_tv)(x)$ where 
	\[
		v(x)=\binom{\varphi(x)}{\varphi^{[1]}(x)}, \qquad v_0 = \binom{0}{-1}, \qquad
		(T_tv)(x) = \int_0^x\begin{pmatrix} 0 & \frac{1}{p(s)} \\ q(s) & 0 \end{pmatrix}v(s) \DD s
	\]
	in the space $X_t=C([0,t],\bb R^2)$ with $\big\|\binom{v_1}{v_2}\big\|=\|v_1\|_\infty+\|v_2\|_\infty$.
	It is easy to see that
	\[
		\|T_t\| \le \max\biggl\{\int_0^t \frac{1}{|p(s)|}\DD s,\int_0^t |q(s)|\DD s\biggr\}.
	\]
	Let $t$ be small enough so that $\|T_t\|<1$.  Then $v=(I-T_t)^{-1}v_0$.  Further, let $P_1$ be the projection
	onto the first component in $\bb R^2$.  Then, with some $C>0$,
	\begin{align*}
		|\varphi(t)+P(t)| &= \big|P_1\bigl[v(t)-\bigl(v_0+(T_tv_0)(t)\bigr)\bigr]\big| \le \|v-(v_0+T_tv_0)\| 
		= \bigg\|\sum_{n=2}^\infty T_t^n v_0\bigg\|
		\\[1ex]
		&\le \frac{\|T_t\|}{1-\|T_t\|}\|T_tv_0\|
		= \frac{\|T_t\|}{1-\|T_t\|}\max_{s\in[0,t]}|P(s)|
		\le \frac{\|T_t\|}{1-\|T_t\|}C|P(t)|
		\ll |P(t)|
	\end{align*}
	since $\|T_t\|\to0$ as $t\to0$.
\end{proof}

\begin{proof}[Proof of \cref{Z212}]
	Let $H$ be as in \eqref{Z216}.
	It follows from \eqref{Z217} and \cref{Z211} that
	\begin{equation}\label{Z218}
		m_1(t) \sim W(t), \qquad m_2(t) \sim \int_0^t w(\xi)P(\xi)^2\DD\xi,
		\qquad m_3(t) \sim -\int_0^t w(\xi)P(\xi)\DD\xi
	\end{equation}
	as $t\to0$.  Since $m_1(t)\gg m_2(t)$, we can apply \cref{Z215}.
	Let us consider the behaviour of $m_2\circ m_1^{-1}$.
	With the substitution $s=W(\xi)$ we obtain
	\begin{equation}\label{Z219}
	\begin{aligned}
		\bigl(m_2\circ m_1^{-1}\bigr)(t) &\sim \int_0^{m_1^{-1}(t)}w(\xi)P(\xi)^2\DD\xi
		= \int_0^{W(m_1^{-1}(t))}P\bigl(W^{-1}(s)\bigr)^2\DD s
		\\[1ex]
		&= G\circ\bigl(W\circ m_1^{-1}\bigr)(t).
	\end{aligned}
	\end{equation}
	Note that 
	\begin{equation}\label{Z228}
		\bigl(W\circ m_1^{-1}\bigr)(t)\sim t
	\end{equation}
	by \eqref{Z218}.

	Assume first that (b) holds.  Then $G$ is regularly or rapidly varying and hence, 
	by \eqref{Z219} and \eqref{Z228}, also $m_2\circ m_1^{-1}$ is regularly or rapidly varying.
	Moreover, the limit
	\begin{align*}
		& \lim_{t\to0}\frac{m_3(t)}{\sqrt{m_1(t)m_2(t)}\,}
		= -\lim_{t\to0}\frac{\int_0^t w(\xi)P(\xi)\DD\xi}{\sqrt{W(t)\int_0^t w(\xi)P(\xi)^2\DD\xi}\,}
		\\[1ex]
		&= -\lim_{t\to0}\frac{\int_0^{W(t)}P(W^{-1}(s))\DD s}{\sqrt{W(t)\int_0^{W(t)}P(W^{-1}(s))^2\DD s}\,}
		= -\lim_{x\to0}\frac{\int_0^x P(W^{-1}(s))\DD s}{\sqrt{x\int_0^x P(W^{-1}(s))^2\DD s}\,}
	\end{align*}
	exists.  Hence \cref{Z215} implies that (a) is satisfied.
	
	Conversely, assume that (a) holds.  By \cref{Z215}, the function $m_2\circ m_1^{-1}$ 
	is regularly or rapidly varying, and hence $G(t)\sim (m_2\circ m_1^{-1})\circ (W\circ m_1^{-1})^{-1}(t)$
	is regularly or rapidly varying.  Furthermore, the limit in \eqref{Z220} exists.
	This shows that (b) is satisfied.
	
	Now assume that (a) and (b) hold and that $G$ is regularly varying.  
	It follows from \eqref{Z219} and \eqref{Z228} that $(m_2\circ m_1^{-1})(t)\sim G(t)$.
	Since $t\mapsto t(m_2\circ m_1^{-1})(t)$ is regularly varying with index $\rho+1\ge1$,
	the solution $\hat t(r)$ of \eqref{Z226} satisfies $\hat t(r)\sim \hat x(r)$.
	The remaining assertions follow from \cref{Z215}.
\end{proof}

\begin{proof}[Proof of \cref{Z222}]
	First assume that (b)$'$ holds and let $\gamma=\ind|P\circ W^{-1}|$.  
	Then $G$ is regularly or rapidly varying by \cref{Z122}\,(i).  If $|P\circ W^{-1}|$ is regularly varying, then
	\begin{equation}\label{Z230}
		G(x) \sim \frac{1}{2\gamma+1}xP\bigl(W^{-1}(x)\bigr)^2
	\end{equation}
	and
	\begin{equation}\label{Z260}
		\delta = -\lim_{x\to0}\frac{\int_0^x P(W^{-1}(s))\DD s}{\sqrt{x\int_0^x P(W^{-1}(s))^2\DD s}\,}
		= -\lim_{x\to0}\frac{\frac{1}{\gamma+1}xP(W^{-1}(x))}{\sqrt{x\frac{1}{2\gamma+1}xP(W^{-1}(x))^2}\,}
		= -\varepsilon\frac{\sqrt{2\gamma+1}\,}{\gamma+1}
	\end{equation}
	exists; hence (b) in \cref{Z212} holds.  By the latter also (a) holds, $\rho=2\gamma+1$
	and $\alpha=\frac{\gamma}{\gamma+1}$, where $\rho=\ind G$.
	
	Next assume that $p(t)>0$ a.e.\ or $p(t)<0$ a.e., and that (a) holds.  Then $G$ is regularly varying by \cref{Z212}.
	Since $|P\circ W^{-1}|$ is monotonic increasing, it follows from \cref{Z122}\,(ii)
	that $|P\circ W^{-1}|$ is regularly or rapidly varying.
	
	Now, assume that (a) and (b)$'$ hold and that $|P\circ W^{-1}|$ is regularly varying.
	Define
	\[
		\tilde x(r) \DE W(\tilde t(r)), \qquad
		F(x) \DE x\,\big|\bigl(P\circ W^{-1}\bigr)(x)\big|, \qquad
		\lambda \DE (2\gamma+1)^{\frac{1}{2(\gamma+1)}}.
	\]
	It follows from \eqref{Z229} and \eqref{Z230} that
	\begin{equation}\label{Z232}
		\frac{1}{r} = \bigl[\hat x(r)G\bigl(\hat x(r)\bigr)\bigr]^{\frac12}
		\sim \frac{1}{\sqrt{2\gamma+1}\,}F\bigl(\hat x(r)\bigr).
	\end{equation}
	Since $F$ is regularly varying at $0$ with index $\gamma+1$, we obtain from \eqref{Z231} that
	\begin{equation}\label{Z233}
	\begin{aligned}
		\frac{1}{r} &= W\bigl(\tilde t(r)\bigr)\big|P\bigl(\tilde t(r)\bigr)\big| = F\bigl(\tilde x(r)\bigr)
		= F\Bigl(\frac{1}{\lambda}\lambda\tilde x(r)\Bigr)
		\\[1ex]
		&\sim \Bigl(\frac{1}{\lambda}\Bigr)^{\gamma+1}F\bigl(\lambda\tilde x(r)\bigr)
		= \frac{1}{\sqrt{2\gamma+1}\,}F\bigl(\lambda\tilde x(r)\bigr).
	\end{aligned}
	\end{equation}
	Combining \eqref{Z232} and \eqref{Z233} we arrive at $\hat x(r)\sim\lambda\tilde x(r)$ by \cref{Z254};
	note that $\ind F=\gamma+1>0$.
	Hence \eqref{Z234} and \eqref{Z231} yield
	\[
		\qSL(rz) \sim i\omega_{\alpha,\delta}\Bigl(\frac{z}{i}\Bigr)^\alpha r\lambda W\bigl(\tilde t(r)\bigr)
		= i\lambda\omega_{\alpha,\delta}\Bigl(\frac{z}{i}\Bigr)^\alpha\frac{1}{|P(\tilde t(r))|}.
	\]
	Let us simplify the coefficient of $\frac{1}{|P(\tilde t(r))|}$.
	First note that $|\delta|=\sqrt{1-\alpha^2}$ by \eqref{Z260}.
	If $\alpha=0$, then $\gamma=0$, $\lambda=1$ and $\omega_{\alpha,\delta}=-i\delta=i\varepsilon$
	and hence $i\lambda\omega_{\alpha,\delta}\bigl(\frac{z}{i}\bigr)^\alpha=-\varepsilon$.
	Now let $\alpha\in(0,1)$.  Then
	\begin{align*}
		i\lambda\omega_{\alpha,\delta}\Bigl(\frac{z}{i}\Bigr)^\alpha 
		&= i(2\gamma+1)^{\frac{1}{2(\gamma+1)}}(i\alpha\delta)^{1+\alpha}\frac{\Gamma(-\alpha)}{\Gamma(2+\alpha)}
		\Bigl(\frac{z}{i}\Bigr)^\alpha
		\\[1ex]
		&= i\Bigl(\frac{1+\alpha}{1-\alpha}\Bigr)^{\frac{1-\alpha}{2}}
		\biggl[-i\alpha\varepsilon(1-\alpha)\Bigl(\frac{1+\alpha}{1-\alpha}\Bigr)^{\frac12}\biggr]^{1+\alpha}
		\frac{\Gamma(1-\alpha)}{-\alpha(1+\alpha)\Gamma(1+\alpha)}\Bigl(\frac{z}{i}\Bigr)^\alpha
		\\[1ex]
		&= -i(-i\varepsilon)^{1+\alpha}\Bigl(\frac{z}{i}\Bigr)^\alpha \alpha^\alpha(1-\alpha)^\alpha
		\frac{\Gamma(1-\alpha)}{\Gamma(1+\alpha)}
		\\[1ex]
		&= -\varepsilon(-\varepsilon z)^\alpha \alpha^\alpha(1-\alpha)^\alpha\frac{\Gamma(1-\alpha)}{\Gamma(1+\alpha)},
	\end{align*}
	which proves \eqref{Z235}.
	
	Let us show the last statement.  Since $|\delta|=\sqrt{1-\alpha^2}$, 
	we have $|\arg\omega_{\alpha,\delta}|=\frac{\pi}{2}(1-\alpha)$,
	and hence $C_{\alpha,\phi}$, as introduced in \cref{Z212}, coincides with \eqref{Z262} by \eqref{Z79}.
	From \eqref{Z263} and \eqref{Z230} we can now derive \eqref{Z264}.
\end{proof}

\begin{remark}\label{Z213}
	When \eqref{Z210} is singular at $0$ but in the limit circle case
	and non-oscillatory for $z=0$ (i.e.\ zeros of solutions of \eqref{Z210} with $z=0$ 
	do not accumulate at $0$), 
	then, for $\varphi$ and $\theta$, one can choose a principal and a non-principal solution 
	of \eqref{Z210} with $z=0$, i.e.\ $\varphi(t)\ll\theta(t)$ as $t\to0$; 
	see, e.g.\ \cite[Theorem~2.1]{niessen.zettl:1992}.
	If one knows the asymptotic behaviour of $\varphi$ and $\theta$,
	then one can obtain the asymptotic behaviour of $\qSL$.
\end{remark}

\medskip

\noindent
Let us consider a simple example to illustrate \cref{Z222}.

\begin{example}\label{Z267}
	Let $p,q,w$ satisfy the conditions at the beginning of the section, assume that
	\[
		P(x) \sim c_1x^{\rho_1}, \qquad W(x) \sim c_2x^{\rho_2} \qquad \text{as} \ x\to0
	\]
	with $c_1,c_2,\rho_1,\rho_2>0$ and suppose that \eqref{Z214} holds.
	Then $(P\circ W^{-1})(x)\sim cx^{\frac{\rho_1}{\rho_1}}$ with some $c>0$ and hence (b)$'$
	in \cref{Z222} is satisfied with $\gamma=\frac{\rho_1}{\rho_2}$.
	It is easy to see that $\tilde t(r)$, the solution of \eqref{Z231}, fulfils
	\[
		\tilde t(r) \sim (c_1c_2)^{-\frac{1}{\rho_1+\rho_2}}r^{-\frac{1}{\rho_1+\rho_2}}
		\qquad \text{as} \ r\to\infty.
	\]
	Plugging this expression into \eqref{Z235} we obtain
	\[
		\qSL(z) \sim -C(-z)^{\frac{\rho_1}{\rho_1+\rho_2}} \qquad\text{as} \ |z|\to\infty
	\]
	uniformly in each sector $\{z\in\bb C:\psi\le\arg z\le\pi-\psi\}$ with $\psi\in(0,\frac{\pi}{2})$,
	where
	\[
		C = \biggl(\frac{\rho_1\rho_2}{(\rho_1+\rho_2)^2}\biggr)^{\frac{\rho_1}{\rho_1+\rho_2}}
		\cdot\frac{\Gamma\bigl(\frac{\rho_2}{\rho_1+\rho_2}\bigr)}{\Gamma\bigl(\frac{2\rho_1+\rho_2}{\rho_1+\rho_2}\bigr)}
		\cdot c_1^{-\frac{\rho_2}{\rho_1+\rho_2}}c_2^{\frac{\rho_1}{\rho_1+\rho_2}}.
	\]
	When $c_1=c_2=1$ and $\rho_1=\rho_2=1$, this reduces to $\qSL(z)\sim-(-z)^{\frac12}$ as $|z|\to\infty$.
\end{example}

\subsection{Krein strings}
\label{Z67}

A Krein string is a pair $\strS[L,\ms m]$ where
\begin{Itemize}
\item
	$L\in[0,\infty]$;
\item
	$\ms m$ is a non-decreasing, left-continuous, $[0,\infty)$-valued function
	on $[0,L)$ with $\ms m(0)=0$.
\end{Itemize}
The number $L$ is called the length of the string, and the function $\ms m$
its mass distribution function.
The string equation can be written as an integro-differential equation,
\[
	y_+'(x) + z\int_{[0,x]}y(t)\DD\ms m(t) = 0, \qquad x\in[0,L),
\]
where $z\in\bb C$ is the spectral parameter and $y_+'$ denotes the right-hand
derivative of $y$; see, e.g.\ \cite{kac.krein:1968}.

Given a string, one can construct a function $\qstr$, the
\emph{principal Titchmarsh--Weyl coefficient} of the string;
see \cite{kac.krein:1968}.
This function belongs to the Stieltjes class $\mc S$, i.e.\ it is analytic
on $\bb C\setminus[0,\infty)$, has non-negative imaginary part in
the upper half-plane, and takes positive values on $(-\infty,0)$.
The class $\mc S$ can also be characterised as follows: a function $q$ belongs to $\mc S$
if and only if its restriction to $\bb C^+$ is a Nevanlinna function
and $zq(z)$ is also a Nevanlinna function.
It is a fundamental theorem proved by M.\,G.~Krein that the assignment
\[
	\strS[L,\ms m]\mapsto \qstr
\]
sets up a bijection from the set of all strings onto the Stieltjes class.

With a string $\strS[L,\ms m]$ we can associate the following Hamiltonian
\begin{equation}\label{Z123}
	H(t)
	=
	\begin{cases}
		\begin{pmatrix} 1 & - \ms m(t) \\[0.5ex] - \ms m(t) & \ms m(t)^2 \end{pmatrix}
		& \text{if} \ t\in(0,L),
		\\[3ex]
		\begin{pmatrix} 0 & 0 \\ 0 & 1 \end{pmatrix}
		& \text{if} \ L+\int_0^L \ms m(s)^2\RD s<\infty \ \text{and} \ t>L.
	\end{cases}
\end{equation}
By \cite[Theorem~4.2]{kaltenbaeck.winkler.woracek:bimmel} the Weyl coefficient, $q_\strS$,
of the string is connected with the Weyl coefficient of the Hamiltonian $H$ via
\[
	\qstr(z) = \frac{q_H(z)}{z}\,, \qquad z\in\bb C^+.
\]
In the following theorem we characterise those Krein strings whose principal Titchmarsh--Weyl coefficients
are regularly varying on the negative real axis.

\begin{theorem}\label{Z121}
	Let $\strS[L,\ms m]$ be a Krein string, assume that $L\in(0,\infty]$ and $\min(\supp \ms m)=0$,
	and let $q_{\strS}$ be the principal Titchmarsh--Weyl coefficient.
	Then the following statements are equivalent:
	\begin{Enumeratealph}
	\item
		$\ms m$ is regularly or rapidly varying at $0$;
	\item
		the function $r\mapsto\qstr(-r)$ is regularly varying at $\infty$.
	\end{Enumeratealph}
	Assume that \textup{(a)} and \textup{(b)} hold and set $\gamma\DE\ind\ms m$, $\beta\DE\ind(r\mapsto\qstr(-r))$. 
	Then $\gamma\in[0,\infty]$, $\beta\in[-1,0]$,
	\begin{equation}\label{Z142}
		\beta
		=
		\begin{cases}
			-\frac{1}{\gamma+1} & \text{when} \ \gamma<\infty,
			\\[0.5ex]
			0 & \text{when} \ \gamma=\infty,
		\end{cases}
	\end{equation}
	and
	\[
		\qstr(rz) \sim C_\beta(-z)^\beta\mr t(r), \qquad r\to\infty,
	\]
	locally uniformly for $z\in\bb C\setminus[0,\infty)$,
	where $\mr t(r)$ is the unique solution of
	\begin{equation}\label{Z143}
		\mr t(r)\int_0^{\mr t(r)} \ms m(s)^2\RD s = \frac{1}{r^2}\,,
		\qquad r>0,
	\end{equation}
	and
	\begin{equation}\label{Z145}
		C_\beta =
		\begin{cases}
			|\beta|^{\frac{\beta}{2}+1}(\beta+1)^\beta(\beta+2)^{\frac{\beta}{2}}
			\frac{\Gamma(-\beta)}{\Gamma(\beta+1)}
			& \text{when} \ \beta\in(-1,0),
			\\[1ex]
			1 & \text{when} \ \beta=0 \ \text{or} \ \beta=-1.
		\end{cases}
	\end{equation}
\end{theorem}

\begin{remark}\label{Z224}
\rule{0ex}{1ex}
\begin{Enumeratealph}
\item
	Most statements were proved in \cite[Theorem~2]{kasahara:1975}
	apart from the relation of the asymptotic behaviour of $\qstr$ and $\ms m$
	in the case when $\ms m$ is rapidly varying.
	Note that in the formula below equation (12) in \cite{kasahara:1975}
	it should be $\Gamma(1-\alpha)$ instead of $(1-\alpha)$ in the denominator.
	
	In \cite[Theorem~5.1]{eckhardt.kostenko.teschl:2018} a similar theorem is stated for $\gamma\in(0,\infty)$
	but only some special cases are proved.
\item
	When $\gamma<\infty$, one knows the asymptotic behaviour of the left-hand side of \eqref{Z143},
	and hence, instead of $\mr t(r)$, one can use the solution $\hat t(r)$ of 
	\[
		\hat t(r)\ms m\bigl(\hat t(r)\bigr) = \frac{\sqrt{2\gamma+1}\,}{r}.
	\]
	One can use this relation also to recover the asymptotic behaviour of $\ms m$ from the
	asymptotic behaviour of $\qstr$.
\item
	It is well known that $\lim_{r\to\infty}\qstr(-r) = \min(\supp\ms m)$.
	Hence, under the assumption of the theorem we have $\lim_{r\to\infty}\qstr(-r)=0$.
\end{Enumeratealph}
\end{remark}

\medskip

\noindent
For the proof of \cref{Z121} we need a lemma.

\begin{lemma}\label{Z239}
	Let $q$ be a function from the Stieltjes class $\mc S$ and let $\ms g$ be a regularly varying function
	with index $\beta$.
	Assume that there exists $z_0\in\bb C\setminus[0,\infty)$ such that
	\begin{equation}\label{Z240}
		q(rz_0) \sim (-z_0)^\beta\ms g(r), \qquad r\to\infty.
	\end{equation}
	Then
	\begin{equation}\label{Z241}
		q(rz) \sim (-z)^\beta\ms g(r), \qquad r\to\infty,
	\end{equation}
	locally uniformly for $z\in\bb C\setminus[0,\infty)$.
\end{lemma}

\begin{proof}
	It follows from \cite[Lemma~S1.5.1]{kac.krein:1968a} that
	\[
		q(z^2) = \frac{\hat q(z)}{z}, \qquad z\in\bb C^+,
	\]
	where $\hat q$ is a symmetric Nevanlinna function, i.e.\ a Nevanlinna function 
	that satisfies $\hat q(-\overline z)=-\overline{\hat q(z)}$ for $z\in\bb C^+$.
	Set $\ms f(r)\DE r\ms g(r^2)$, which is regularly varying with index $\alpha\DE 2\beta+1$.
	Further, let $w_0\in\bb C^+$ such that $w_0^2=z_0$.  Then
	\[
		\frac{\hat q(rw_0)}{\ms f(r)} = \frac{rw_0q(r^2w_0^2)}{r\ms g(r^2)}
		= w_0\frac{q(r^2z_0)}{\ms g(r^2)} \to w_0(-z_0)^\beta \ne 0,
		\qquad r\to\infty.
	\]
	Hence we can apply \cite[Theorem~3.1]{langer.woracek:kara} to deduce that
	\[
		\hat q(rw) \sim i\omega\Bigl(\frac{w}{i}\Bigr)^\alpha\ms f(r), \qquad r\to\infty,
	\]
	locally uniformly for $w\in\bb C^+$ with some $\omega\ne0$.
	The following relations hold locally uniformly for $w\in\bb C^+$:
	\[
		q(rw^2) = \frac{\hat q(r^{\frac12}w)}{r^{\frac12}w}
		\sim i\omega\Bigl(\frac{w}{i}\Bigr)^\alpha\frac{\ms f(r^{\frac12})}{r^{\frac12}w}
		= \omega\Bigl(\frac{w}{i}\Bigr)^{\alpha-1}\frac{\ms f(r^{\frac12})}{r^{\frac12}}
		= \omega(-w^2)^\beta\ms g(r), \qquad r\to\infty.
	\]
	Replacing $w^2$ by $z\in\bb C\setminus[0,\infty)$ and comparing with \eqref{Z240}
	we obtain $\omega=1$, and hence \eqref{Z241} holds.
\end{proof}

\begin{proof}[Proof of \cref{Z121}]
\begin{Steps}
\item
	Define $H$ as in \eqref{Z123}.
	The assumption $\min(\supp\ms m)=0$ implies that $h_2$ does not vanish
	in any neighbourhood of $0$.
	Since $\lim_{t\to0}h_2(t)$ exists, the function $\mf t$ as defined in \cref{Z10}
	is regularly varying at $0$ with index $1$.
	Hence the assumptions of \cref{Z10} are satisfied.

\item
	First assume that $\ms m\in R_\gamma(0)$.
	\Cref{Z122}\,(i) implies that $m_2\in R_{2\gamma+1}(0)$.
	If $\gamma\in[0,\infty)$, then, again by \cref{Z122}\,(i), we have
	\begin{align*}
		m_2(t) &= \int_0^t \ms m(s)^2\RD s \sim \frac{1}{2\gamma+1}t\ms m(t)^2,
		\\[1ex]
		m_3(t) &= -\int_0^t \ms m(s)\RD s \sim -\frac{1}{\gamma+1}t\ms m(t)
	\end{align*}
	as $t\to0$, and hence the limit
	\[
		\delta = \lim_{t\to0}\frac{m_3(t)}{\sqrt{m_1(t)m_2(t)}\,}
		= \lim_{t\to0}\frac{-\frac{1}{\gamma+1}t\ms m(t)}{\sqrt{t\frac{1}{2\gamma+1}t\ms m(t)^2}\,}
		= -\frac{\sqrt{2\gamma+1}\,}{\gamma+1}
	\]
	exists.
	\Cref{Z10} implies that \eqref{Z57} and \eqref{Z98} hold and that $\delta=0$ when $\gamma=\infty$.
	Let $\alpha=\ind\ms a_H$.
	Then $\alpha=\frac{\gamma}{\gamma+1}$ when $\gamma\in[0,\infty)$
	and $\alpha=1$ when $\gamma=\infty$ by \eqref{Z78}.
	Since $|\delta|=\sqrt{1-\alpha^2}$, we obtain from \eqref{Z6} that $\arg\omega_{\alpha,\delta}=-\frac{\pi}{2}(\alpha-1)$.
	Relations \eqref{Z57} and \eqref{Z17} imply that, locally uniformly for $z\in\bb C^+$ and as $r\to\infty$,
	\begin{align}
		\qstr(rz) &= \frac{q_H(rz)}{rz}
		\sim i\omega_{\alpha,\delta}\frac{1}{rz}\Bigl(\frac{z}{i}\Bigr)^\alpha\ms a_H(r)
		= \omega_{\alpha,\delta}\Bigl(\frac{z}{i}\Bigr)^{\alpha-1}\frac{\ms a_H(r)}{r}
		\nonumber\\[1ex]
		&= |\omega_{\alpha,\delta}|e^{-i\frac{\pi}{2}(\alpha-1)}\Bigl(\frac{z}{i}\Bigr)^{\alpha-1}
		m_1\bigl(\mr t(r)\bigr)
		= |\omega_{\alpha,\delta}|(-z)^{\alpha-1}\,\mr t(r).
		\label{Z242}
	\end{align}
	Since $\mr t(r)=\frac{\ms a_H(r)}{r}$ is regularly varying with index $\beta=\alpha-1$,
	we can use \cref{Z239} to extend the asymptotic relation \eqref{Z242} to $z\in\bb C\setminus[0,\infty)$, 
	again locally uniformly.
	Clearly, the defining relation \eqref{Z38} for $\mr t$ turns into \eqref{Z143}.
	Let us determine the value of $C_\beta=|\omega_{\alpha,\delta}|$, where $\omega_{\alpha,\delta}$ is given by \eqref{Z98}.
	When $\gamma=0$, then $\alpha=0$ and $\delta=-1$ and hence $\omega_{\alpha,\delta}=i$;
	when $\gamma=\infty$, then $\alpha=1$ and hence $\omega_{\alpha,\delta}=1$.
	Now let $\gamma\in(0,\infty)$.  Then
	\[
		\beta = -\frac{1}{\gamma+1}, \qquad
		\delta = -\frac{\sqrt{2\gamma+1}\,}{\gamma+1} = \beta\sqrt{-\frac{2}{\beta}-1}
		= \sqrt{|\beta|(\beta+2)}
	\]
	and hence
	\begin{align*}
		|\omega_{\alpha,\delta}| &= \bigg|(i\alpha\delta)^{1+\alpha}
		\frac{\Gamma(-\alpha)}{\Gamma(2+\alpha)}\bigg|
		= (\beta+1)^{\beta+2}\bigl[|\beta|(\beta+2)\bigr]^{\frac{\beta}{2}+1}
		\frac{|\Gamma(-\beta-1)|}{\Gamma(\beta+3)}
		\\[1ex]
		&= |\beta|^{\frac{\beta}{2}+1}(\beta+1)^\beta(\beta+2)^{\frac{\beta}{2}}
		\frac{\Gamma(-\beta)}{\Gamma(\beta+1)},
	\end{align*}
	which proves \eqref{Z145}.

\item
	Conversely, assume that (b) holds and set $\ms g(r)\DE\qstr(-r)$.
	It follows from \cref{Z239} with $z_0=-1$ that
	\[
		q_H(ri) = ri\qstr(ri) \sim i(-i)^\beta r\ms g(r),
	\]
	where the function $r\mapsto r\ms g(r)$ is regularly varying at $\infty$
	with index $\alpha=\beta+1$.
	This shows that (i) in \cref{Z10} is satisfied and hence also (ii) in that theorem,
	where, by \eqref{Z78}, the index of $m_2$
	is $\rho_2=-\frac{\beta+2}{\beta}$ if $\beta\in[-1,0)$ and $\rho_2=\infty$ if $\beta=0$.
	Since $\ms m^2$ is non-decreasing, \cref{Z122}\,(ii) implies
	that $\ms m^2\in R_{\rho_2-1}(0)$ and therefore $\ms m \in R_{-\frac{\beta+1}{\beta}}(0)$
	if $\beta\in[-1,0)$ and $\ms m\in R_{\infty}(0)$ if $\beta=0$.
\end{Steps}
\end{proof}

\subsection{Generalised indefinite strings}
\label{Z68}

Let $L\in(0,\infty]$, let $\upsilon$ be a positive Borel measure on $[0,L)$, 
and let $\chi\in H_{\textsf{loc}}^{-1}([0,L))$ be real, where $H_{\textsf{loc}}^{-1}([0,L))$
is the dual space of $H_{\textsf{c}}^1([0,L))=\{f\in H^1([0,L)):\supp f \ \text{compact in} \ [0,L)\}$.
We consider the equation
\begin{equation}\label{Z125}
	-y'' = zy\chi + z^2y\upsilon
\end{equation}
on $[0,L)$, which is to be understood in a distributional sense 
with $y\in H_{\textsf{loc}}^1([0,L))$ and $H_{\textsf{c}}^1([0,L))$ as test function space.

Equations of the form in \eqref{Z125} were introduced by J.~Eckhardt and A.~Kostenko 
in \cite{eckhardt.kostenko:2016} and generalise Krein strings in three directions.
First, instead of a measure $\chi$ one allows $\chi$ to be more singular,
namely a distribution in $H_{\textsf{loc}}^{-1}([0,L))$;
second the positivity of $\chi$ is dropped; third, a term that depends on $z^2$ is included.
These equations, called \emph{generalised indefinite strings}, can also be seen as generalisations of
Krein strings with a signed measure or equations of the form \eqref{Z125}
with $\chi$ a signed measure and $\upsilon$ a positive measure as studied 
by H.~Langer in \cite{langer:1976}.
One motivation to consider generalised indefinite strings is the use of the 
inverse scattering approach for solving non-linear wave equations such as the Camassa--Holm equation.

The principal Weyl--Titchmarsh function associated with \eqref{Z125} is 
defined by $\qigs(z)\DE\frac{\psi'(0-;z)}{z\psi(0;z)}$,
$z\in\bb C^+$, where $\psi(\,\cdot\,;z)$ is the (up to a scalar multiple) unique non-trivial solution
of \eqref{Z125} that satisfies 
$\psi'(\,\cdot\,;z)\in L^2(0,L)$, $\psi(\,\cdot\,;z)\in L^2((0,L);\upsilon)$,
and also $\lim_{x\to L}\psi(x;z)=0$ if $L<\infty$;
for details see \cite[Section~5]{eckhardt.kostenko:2016}.

Let $\sfw\in L_{\textsf{loc}}^2([0,L))$ be the normalised anti-derivative of $\chi$
defined by
\[
	\chi(h) = -\int_0^L \sfw(x)h'(x)\RD x,
	\qquad h\in H_{\textsf{c}}^1([0,L));
\]
note that if $\chi$ is a Borel measure, then $\sfw(x)=\chi([0,x))$.
Further, we set
\begin{equation}\label{Z164}
	g(x) \DE \int_0^x (\sfw(t))^2\RD t + \upsilon([0,x)),
	\qquad x\in[0,L),
\end{equation}
which measures the total contribution of the coefficients on $[0,x)$.
Note that the function $g$ appears implicitly in \cite[Proposition~2.6]{eckhardt.kostenko:2024}
in connection with convergence of generalised indefinite strings, 
and in \cite[Theorem~5.3]{eckhardt.kostenko:2023} in connection with discreteness of the spectrum.

\begin{theorem}\label{Z148}
	Assume that $\upsilon(\{0\})=0$ and $g(x)>0$ for every $x\in(0,L)$.
	Then the following conditions are equivalent:
	\begin{Enumeratealph}
	\item
		there exist a regularly varying function $\ms a$ with index $\alpha\in(-1,1)$
		and a constant $\omega\in\bb C\setminus\{0\}$ such that
		\begin{equation}\label{Z169}
			\qigs(ri) \sim i\omega\ms a(r), \qquad r\to\infty;
		\end{equation}
	\item
		$g$ is regularly varying at $0$ with index $\rho>0$ and the limit
		\begin{equation}\label{Z159}
			\delta \DE \lim_{x\to0}\frac{\int_0^x\sfw(\tau)\RD\tau}{\sqrt{xg(x)}}
		\end{equation}
		exists.
	\end{Enumeratealph}
	Assume that \textup{(a)} and \textup{(b)} are satisfied and let $\hat x(r)$ be
	the unique number that satisfies
	\begin{equation}\label{Z265}
		\hat x(r)g\bigl(\hat x(r)\bigr) = \frac{1}{r^2},
		\qquad r>0.
	\end{equation}
	Then $\alpha=\frac{1-\rho}{1+\rho}$ and
	\begin{equation}\label{Z160}
		\qigs(rz) \sim i\omega_{\alpha,\delta}\Bigl(\frac{z}{i}\Bigr)^\alpha\frac{1}{r\hat x(r)},
		\qquad r\to\infty,
	\end{equation}
	locally uniformly for $z\in\bb C^+$, where $\omega_{\alpha,\delta}$ is as in \eqref{Z98}.
\end{theorem}

\begin{proof}
\begin{Steps}
\item
	As in \cite[Proof of Theorem~6.1]{eckhardt.kostenko:2016} we can associate a canonical system
	with \eqref{Z125}.  Set
	\[
		\zeta(x) \DE g(x)+x, \qquad x\in[0,L),
	\]
	and let $\xi$ be the generalised inverse of $\zeta$:
	\[
		\xi(t) \DE \sup\{x\in[0,L):\zeta(x)\le t\}, \qquad t\in[0,\infty),
	\]
	which is an asymptotic inverse of $\zeta$ in the sense of \cref{Z254}.
	Note that $\zeta$ has a unique generalised inverse since $\zeta$ is strictly increasing.
	Moreover, we have
	\begin{equation}\label{Z127}
		\xi\bigl(\zeta(x)\bigr) = x, \quad x\in[0,L); \qquad
		\zeta\bigl(\xi(t)\bigr) = 
		\begin{cases}
			t, & t\in\ran\zeta,
			\\[0.5ex]
			\sup\{s\in\ran\zeta:s\le t\}, & t\notin\ran\zeta.
		\end{cases}
	\end{equation}
	Since $g$ is non-decreasing, $\xi$ is locally absolutely continuous  and $0\le\xi'(t)\le 1$
	for a.e.\ $t\in[0,\infty)$.  Define the Hamiltonian
	\begin{equation}\label{Z163}
		H(t) \DE
		\begin{pmatrix}
			1-\xi'(t) & \xi'(t)\sfw(\xi(t))
			\\[1ex]
			\xi'(t)\sfw(\xi(t)) & \xi'(t)
		\end{pmatrix},
		\qquad t\in[0,\infty).
	\end{equation}
	By \cite[Proof of Theorem~6.1, in particular, (6.11)]{eckhardt.kostenko:2016} we have $\qigs=q_H$.
	The primitive of $H$ is given by
	\[
		M(t) =
		\begin{pmatrix}
			t-\xi(t) & \int_0^{\xi(t)}\sfw(\tau)\RD\tau
			\\[2ex]
			\int_0^{\xi(t)}\sfw(\tau)\RD\tau & \xi(t)
		\end{pmatrix},
		\qquad t\in[0,\infty).
	\]

\item
	Let us first assume that (a) holds.
	It follows from \cref{Z10} and the assumption $\alpha\in(-1,1)$ that $m_1$ and $m_2$
	are regularly varying with indices $\rho_1,\rho_2>0$ and that
	the limit in \eqref{Z56} exists.
	Since $m_2=\xi$, \cref{Z254} implies that $\zeta$ is regularly varying with index $\frac{1}{\rho_2}$.
	Further, for $x\in[0,L)$ we can write
	\begin{equation}\label{Z126}
		g(x) = \zeta(x)-x = \zeta(x)-\xi\bigl(\zeta(x)\bigr) = m_1\bigl(\zeta(x)\bigr),
	\end{equation}
	which shows that $g$ is regularly varying with index $\frac{\rho_1}{\rho_2}$.
	From \eqref{Z127}, \eqref{Z126} and the existence of the limit in \eqref{Z56} we obtain
	\[
		\frac{\int_0^x\sfw(\tau)\RD\tau}{\sqrt{g(x)x}}
		= \frac{\int_0^{\xi(\zeta(x))}\sfw(\tau)\RD\tau}{\sqrt{m_1(\zeta(x))\xi(\zeta(x))}\,}
		= \frac{m_3(\zeta(x))}{\sqrt{m_1(\zeta(x))m_2(\zeta(x))}\,}
		\to \delta
	\]
	as $x\to0$.  Hence (b) is true.
	Moreover, from \eqref{Z78} we obtain $\alpha=\frac{1-\frac{\rho_1}{\rho_2}}{1+\frac{\rho_2}{\rho_2}}=\frac{1-\rho}{1+\rho}$.

\item
	Let us now consider the converse and assume that (b) holds.
	Then $\zeta$ is regularly varying with positive index, and hence also $m_2=\xi$ 
	is regularly varying by \cref{Z254}.
	To show that $m_1$ is regularly varying, requires more work.
	When $t\in\ran\zeta$, we obtain from \eqref{Z127} that
	\begin{equation}\label{Z156}
		m_1(t) = \zeta\bigl(\xi(t)\bigr)-\xi(t) = g\bigl(\xi(t)\bigr).
	\end{equation}
	Let $\tilde g$ be the right-continuous function
	\[
		\tilde g(x) \DE \int_0^x (\sfw(t))^2\RD t + \upsilon([0,x]),
		\qquad x\in[0,L),
	\]
	and set $\tilde\zeta(x)\DE\tilde g(x)+x$.
	Note that $\tilde g$ and $g$ coincide at all points where $\upsilon$ has no point mass.
	Set $\zeta_{\textsf{max}}\DE\sup\{\zeta(x):x\in[0,L)\}$, which satisfies $\zeta_{\textsf{max}}\in(0,\infty]$
	since $g$ does not vanish identically.
	Now let $(t_1,t_2)$ be a maximal interval in $[0,\zeta_{\textsf{max}})\setminus\ran\zeta$.
	Then $0<t_1<t_2<\zeta_{\textsf{max}}$, and there exists $x_1\in(0,L)$ such that
	\begin{equation}\label{Z157}
		t_1 = \zeta(x_1) \qquad\text{and}\qquad t_2 = \tilde\zeta(x_1).
	\end{equation}
	For all $t\in[t_1,t_2]$ we have $\xi(t)=x_1$ and hence
	\begin{align*}
		m_1(t) &= t-\xi(t) = t-x_1,
		\\
		g\bigl(\xi(t)\bigr) &= g(x_1) = \zeta(x_1)-x_1 = t_1-x_1,
		\\
		\tilde g\bigl(\xi(t)\bigr) &= \tilde g(x_1) = \tilde\zeta(x_1)-x_1 = t_2-x_1
	\end{align*}
	by \eqref{Z157}.  This, together with \eqref{Z156}, shows that
	\begin{equation}\label{Z158}
		g\bigl(\xi(t)\bigr) \le m_1(t) \le \tilde g\bigl(\xi(t)\bigr),
		\qquad t\in[0,\zeta_{\textsf{max}}).
	\end{equation}
	It follows from \cref{Z253} that there exist smoothly varying functions $g_1,g_2$
	such that $g_1(x)\le g(x) \le g_2(x)$ for $x\in[0,L)$ and $g_1(x)\sim g_2(x)$ as $x\to0$.
	Since $g$ and $\tilde g$ differ only where they have jumps, we have 
	$g_1(x) \le g(x) \le \tilde g(x) \le g_2(x)$ for $x\in[0,L)$, which implies 
	that $g(x)\sim\tilde g(x)$ as $x\to0$.  Together with \eqref{Z158}, this shows that
	\begin{equation}\label{Z161}
		m_1(t)\sim g\bigl(\xi(t)\bigr), \qquad t\to0,
	\end{equation}
	and hence $m_1$ is regularly varying.
	Moreover, from \eqref{Z159} we obtain
	\[
		\frac{m_3(t)}{\sqrt{m_1(t)m_2(t)}\,}
		\sim \frac{\int_0^{\xi(t)}\sfw(s)\RD s}{\sqrt{g(\xi(t))\xi(t)}\,}
		\to \delta,
		\qquad t\to0.
	\]
	By \cref{Z10}, the statement in (a) is satisfied.

\item
	In order to show \eqref{Z160}, set $F(x)\DE xg(x)$.
	We obtain from \eqref{Z161} that $(m_1m_2)(t) \sim g(\xi(t))\xi(t) = (F\circ\xi)(t)$
	and hence $\xi\circ(m_1m_2)^{-1}\sim F^{-1}$ by \cref{Z254}.
	Together with \eqref{Z17}, this yields
	\[
		\ms a_H(r) = \frac{1}{rm_2(\mr t(r))} = \frac{1}{r\xi(\mr t(r))}
		= \frac{1}{r\bigl(\xi\circ(m_1m_2)^{-1}\bigr)\bigl(\frac{1}{r^2}\bigr)}
		\sim \frac{1}{rF^{-1}\bigl(\frac{1}{r^2}\bigr)} = \frac{1}{r\hat x(r)}
	\]
	as $r\to\infty$.
\end{Steps}
\end{proof}

\begin{remark}\label{Z162}
\rule{0ex}{1ex}
\begin{Enumerate}
\item
	Since the Hamiltonian in \eqref{Z163} is trace-normalised, we can use \cref{Z124}
	to obtain the following equivalence: $\qigs(ir)\sim cr^\alpha$ as $r\to\infty$ with come constant $c\ne0$
	if and only if $g(x)\sim c'x^\rho$ as $x\to0$ with some constant $c'\ne0$ and the limit in \eqref{Z159} exists.
\item
	In a similar way as in \cref{Z148}, one can show equivalences for the boundary cases $\alpha=\pm1$:
	there exists a regularly varying function $\ms a$ with index $1$ ($-1$, respectively) and
	a constant $\omega\in\bb C\setminus\{0\}$ such that \eqref{Z169} holds if and only if $g$ is slowly varying 
	(rapidly varying with index $\infty$, respectively).
	Only one implication needs a slightly different argument: assume that $g$ is slowly varying;
	then $m_2$ is rapidly varying and, since $m_1(t)+m_2(t)=t$, it follows that $m_1(t)\sim t$
	and hence $\alpha=1$.
\item
	\Cref{Z148} goes far beyond \cite[Theorems~6.1 and 6.2]{eckhardt.kostenko.teschl:2018},
	where only two cases are considered: either $\qigs(ri)\to c_1$ as $r\to\infty$ with $c_1\in\overline{\bb C^+}$
	or $\qigs(ri)\sim ic_2r^\alpha$ with $c_2>0$ and $\alpha\in(0,1)$.
	In the latter case, $c_2>0$ implies that $\delta=0$.
	This is a severe restriction since, e.g.\ when $\sfw$ is regularly varying and 
	the first term on the right-hand side of \eqref{Z164} dominates the second term, then $\delta\ne0$.
\item
	It is shown in \cite[Lemma~7.1]{eckhardt.kostenko:2016} 
	that $\lim_{y\to\infty}\frac{\qigs(ir)}{ir}=\upsilon(\{0\})$. 
	Hence the case $\upsilon(\{0\})>0$ corresponds to maximal growth of $\qigs$.
	The other extreme case discussed in \cref{Z90} can also be characterised:
	we have $\qigs(ir)\sim \frac{i}{cr}$ as $r\to\infty$ for some $c>0$ if and only if $g(x)=0$ for $x\in[0,c]$ 
	and $g(x)>0$ for $x>c$.
\item
	One can recover the asymptotic behaviour of $g$ from the behaviour of $\qigs$ in a similar way as
	in \cref{Z236}: from the argument of $\omega$ in \eqref{Z169} one can obtain $|\omega_{\alpha,\delta}|$ 
	via \eqref{Z79}; from this and \eqref{Z160} one can then derive the asymptotic behaviour of $\hat x(r)$ and 
	hence the asymptotics of $g$ via \eqref{Z265}.
\end{Enumerate}
\end{remark}


\printbibliography

{\footnotesize
\begin{flushleft}
	M.~Langer \\
	Department of Mathematics and Statistics \\
	University of Strathclyde \\
	26 Richmond Street \\
	Glasgow G1 1XH \\
	UNITED KINGDOM \\
	email: \texttt{m.langer@strath.ac.uk} \\[5mm]
\end{flushleft}
\begin{flushleft}
	R.~Pruckner \\
	Institute for Analysis and Scientific Computing \\
	Vienna University of Technology \\
	Wiedner Hauptstra{\ss}e 8--10/101 \\
	1040 Wien \\
	AUSTRIA \\
	email: \texttt{raphael.pruckner@tuwien.ac.at} \\[5mm]
\end{flushleft}
\begin{flushleft}
	H.\,Woracek\\
	Institute for Analysis and Scientific Computing \\
	Vienna University of Technology \\
	Wiedner Hauptstra{\ss}e\ 8--10/101 \\
	1040 Wien \\
	AUSTRIA \\
	email: \texttt{harald.woracek@tuwien.ac.at} \\[5mm]
\end{flushleft}
}

\end{document}